\documentclass[12pt]{amsart}
\usepackage{amsmath}
\usepackage{amssymb}
\newtheorem{theorem}{Theorem}[section]
\newtheorem{lemma}[theorem]{Lemma}

\newtheorem{corollary}[theorem]{Corollary}

\theoremstyle{definition}
\newtheorem{definition}[theorem]{Definition}
\newtheorem{definitions}[theorem]{Definitions}
\newtheorem{example}[theorem]{Example}
\newtheorem{examples}[theorem]{Examples}
\newtheorem{hypotheses}[theorem]{Hypotheses}
\newtheorem{definitions and remarks}[theorem]{Definitions and Remarks}

\theoremstyle{remark}
\newtheorem{remark}[theorem]{Remark}
\newtheorem{remarks}[theorem]{Remarks}
\newtheorem{algorithm}[theorem]{Canonical desingularization algorithm}

\numberwithin{equation}{section}

\newcommand{\exc}{\mathrm{exc}}
\newcommand{\inv}{\mathrm{inv}}
\newcommand{\Sing}{\mathrm{Sing}\,}
\newcommand{\supp}{\mathrm{supp}\,}

\newcommand{\codim}{\mathrm{codim}\,}

\newcommand{\al}{{\alpha}}

\newcommand{\de}{{\delta}}
\newcommand{\D}{{\Delta}}

\newcommand{\la}{{\lambda}}
\newcommand{\Om}{{\Omega}}
\newcommand{\p}{{\partial}}
\newcommand{\s}{{\sigma}}
\newcommand{\Sig}{{\Sigma}}
\newcommand{\vp}{{\varphi}}
\newcommand{\io}{{\iota}}

\newcommand{\IN}{{\mathbb N}}
\newcommand{\IQ}{{\mathbb Q}}
\newcommand{\IA}{{\mathbb A}}

\newcommand{\IC}{{\mathbb C}}
\newcommand{\IK}{{\mathbb K}}

\newcommand{\cC}{{\mathcal C}}
\newcommand{\cD}{{\mathcal D}}
\newcommand{\cE}{{\mathcal E}}
\newcommand{\cF}{{\mathcal F}}
\newcommand{\cG}{{\mathcal G}}
\newcommand{\cH}{{\mathcal H}}
\newcommand{\cI}{{\mathcal I}}
\newcommand{\cJ}{{\mathcal J}}
\newcommand{\cO}{{\mathcal O}}

\newcommand{\um}{\underline{m}}

\newcommand{\tx}{{\tilde x}}

\newcommand{\wcI}{{\widehat \cI}}
\newcommand{\wcO}{{\widehat \cO}}

\newcommand{\ucF}{\underline{\cF}}
\newcommand{\ucG}{\underline{\cG}}
\newcommand{\ucH}{\underline{\cH}}
\newcommand{\ucC}{\underline{\cC}}
\newcommand{\ucD}{\underline{\cD}}

\newcommand{\frN}{{\mathfrak N}}

\newcommand{\lbr}{{[\![}}
\newcommand{\rbr}{{]\!]}}

\begin{document}

\title[Desingularization algorithms I]
{Desingularization algorithms\\ I. Role of 
exceptional divisors}

\author{Edward Bierstone}
\address{Department of Mathematics, University of Toronto, Toronto,
Ontario, Canada M5S 3G3}
\email{bierston@math.toronto.edu}
\thanks{The first author's research was supported in part by NSERC 
grant OGP0009070.}

\author{Pierre D. Milman}
\address{Department of Mathematics, University of Toronto, Toronto,
Ontario, Canada M5S 3G3}
\email{milman@math.toronto.edu}
\thanks{The second author's research was supported in part by NSERC
grant OGP0008949 and the Killam Foundation.}

\subjclass{Primary 14E15, 32S45; Secondary 32S15, 32S20}

\keywords{resolution of singularities, desingularization invariant,
blowing-up, exceptional divisor}

\begin{abstract}
The article is about a ``desingularization principle''
(Theorem 1.14) common to various canonical desingularization 
algorithms in characteristic zero, and the roles played by the
exceptional divisors in the underlying local construction.
We compare algorithms of the authors and of Villamayor and his
collaborators, distinguishing between the fundamental effect of
the way the exceptional divisors are used,
and different theorems obtained because of
flexibility allowed in the choice of ``input data''. We show
how the meaning of ``invariance'' of the desingularization
invariant, and the efficiency of the algorithm depend on the
notion of ``equivalence'' of collections of local data used 
in the inductive construction.
\end{abstract}

\maketitle
\setcounter{tocdepth}{1}
\tableofcontents

\section{Introduction}
This article is about (1) a ``desingularization principle''
(Theorem 1.14 below) that is common to various algorithms for canonical
resolution of singularities in characteristic zero (embedded desingularization,
principalization of an ideal, etc.); (2) the roles played by the exceptional
divisors arising from blowings-up in the local inductive construction at the
heart of the desingularization principle. The exceptional divisors play two
roles -- one local and one in the passage from local to global -- that are
reflected in different ways in the invariant whose maximum loci provide the
centres of blowing up.

This paper is to be followed by ``Desingularization algorithms
II. Locally binomial varieties''.

The main results on resolution of singularities in characteristic zero
originate in the towering work of Hironaka \cite{Hann}. Among our aims
here is to compare the various algorithms for canonical desingularization
developed by the authors \cite{BMihes, BMjams, BMmega, BMinv},
and by Villamayor and his collaborators \cite{V1, V2, 
EVacta, EVweak, BV}. (Lipman \cite{Lip} has raised the question
of such a comparison.)
The approaches have much in common but, apart from different
applications of the general principle to particular results,
there are important differences due to the ways that the
exceptional divisors are used, affecting the invariant, the meaning
of ``invariance'', and the choice of centre of blowing up. 

The general {\it desingularization principle}
(Theorem 1.14) can be stated
roughly as follows: An initial choice of invariant that distinguishes
between ``general'' and ``special'' points and that satisfies certain
basic properties, can be extended to a desingularization invariant
defined over sequences of ``admissible'' blowings-up, satisfying
several simple properties which show that special points can be
eliminated by successively blowing up the maximum loci of the invariant.

The desingularization principle depends on a local inductive construction
(Section 2). We distinguish between the more fundamental effect of
the way that the exceptional divisors are used in the
local construction, and different theorems that can be obtained
because of the flexibility allowed in the choice of ``input data''
and ``when we stop running the algorithm'' (cf. \S1.2). For
example, ``embedded desingularization'' (Example 1.19(2) below, 
\cite[Main Thm. I]{Hann}, \cite
[Thms. 1.6, 11.14]{BMinv}, \cite{V2}) and ``principalization of
an ideal'' (Example 1.19(1) below, \cite[Main Thm. 2]{Hann}, 
\cite[Thm. 1.10]{BMinv})
are both applications of Theorem 1.14, the only differences
being in the notion of transformation used 
and the choice of initial invariant (``strict'' transform and 
the Hilbert-Samuel function for
embedded desingularization, or ``weak'' transform and the order 
of an ideal for principalization; see \S1.4 and Examples 1.8). 
Likewise, the ``weak embedded
desingularization'' theorem of \cite{EVweak} (Example 1.19(4) below)
is obtained by stopping the principalization algorithm early.
(See also \S6.2.) 

The preceding point of view is not always clear in the literature.
It is further developed in Section 6, 
where we treat universal embedded desingularization of (not 
necessarily embedded) spaces (\S6.1), the relationship between weak
and strict transform (\S6.2), and an extension of the general
desingularization principal to
parametrized families (\S6.3; cf. \cite{ENV}).

A second related aim is an understanding of the meaning 
of ``invariance''.
The desingularization invariants are invariants of what? In other words,
on what do the invariants and therefore the algorithms depend? These
questions are closely connected to the notion of ``equivalence'' of
collections of local data used in the inductive construction. The collection
of local data is called a ``presentation'' in \cite{BMinv}. A presentation is
not invariant or canonical. The philosophy behind \cite{BMinv}
is to introduce an equivalence relation on presentations (using sequences
of ``test blowings-up'' that depend on the accumulating exceptional
divisors) so that the corresponding equivalence class is an invariant and
certain natural numerical characters of a presentation depend only on the
equivalence class. A presentation is similar to the notion of ``basic
object'' or ``idealistic space'' used in \cite{V2, EVacta}. Both
originate in Hironaka's idea of an ``idealistic exponent'' \cite{Hid}.
But the equivalence class of a presentation is strictly smaller than
that of an idealistic exponent or basic object used in these articles;
the ``residual multiplicities'' appearing in the desingularization invariants
of \cite{BMinv} or even those of
\cite{V2, EVacta} depend only on the equivalence class
as a presentation, but not only on the equivalence class as a basic object.
(See Remarks 2.5 and \S\S3.5, 5.3 below.)

Our third aim is to show that for certain natural classes of algebraic or
analytic varieties (e.g., locally toric or toric), a combinatorial structure
can be used to simplify one or both of the roles played by the exceptional
divisors. This is the subject of ``Desingularization algorithms II.
Locally binomial varieties''. 
A {\it locally binomial variety} is a variety
defined in local coordinates by systems of binomial equations. (A {\it
binomial} means a difference of two monomials with no common factor.)
A locally toric variety is simply a locally binomial variety that is normal.
Locally binomial varieties are a very natural class, both as a testing-ground
for general conjectures in algebraic or analytic geometry, and because many
general questions and computational problems can be reduced to the binomial
case.

The ideas above will be made more precise in the remainder of this introductory
section. The notion of a presentation and the desingularization algorithms
of \cite{BMinv}, \cite{V2} and \cite{EVacta} are recalled in 
Sections 2 and 3, where we
also outline the proof of the desingularization principle 
Theorem 1.14. Section 4 presents a worked example
(the details of Example 1.2 below), and Section 5 
deals with the idea of equivalence of presentations.
In comparing our approach to
desingularization with that of Villamayor {\it et al}, we nevertheless
use the more analytic language of \cite{BMinv, BMdc}.
We treat equivalence of presentations (Section 5) using
transformation formulas for differential operators developed by
Hironaka \cite[Sect. 8]{Hid}, Giraud \cite{giraud} 
and Encinas and Villamayor \cite[Sect. 4]{EVacta}.

\subsection{Role of exceptional divisors}
The desingularization invariant $\inv(\cdot) \, = \, \inv_X(\cdot)$ 
(or $\inv_\cJ(\cdot)$) will be described precisely below. The invariant
is defined recursively over a sequence of ``admissible'' blowings-up
of a singular space $X$ (or a coherent ideal sheaf $\cJ$). Let $X_j$
(or $\cJ_j$) denote the transform of $X$ (or $\cJ$) in ``year'' $j$
(i.e., after the first $j$ blowings-up; see \S1.5). If $a \in X_j$ (or
$\supp \cJ_j$), then $\inv(a)$ is a finite sequence $(\io(a), s_1(a);
\nu_2(a), s_2(a); \ldots)$, where $\io(\cdot)$ is the initial
invariant and the successive pairs are themselves
defined using data (``presentations'') that involve ``maximal contact''
subspaces $N_{r-1}(a)$ of increasing codimension in an ambient manifold.
For each $r$, the truncated invariant $\inv_{r - 1/2}(a) := (\io(a), s_1(a);
\ldots ; \nu_r(a))$ determines a block $E^r(a)$ of ``old'' exceptional
divisors passing through $a$ (those which do not necessarily have normal
crossings with respect to $N_r(a)$ and which do not appear in a previous
block $E^q(a)$, $1 \leq q < r$; see Definitions 1.15), 
and a block $\cE_r(a)$ of ``new''
exceptional divisors (those exceptional divisors passing through $a$,
not in $E^1(a) \cup \cdots \cup E^r(a)$). The exceptional divisors are
used in two ways; they are:

\begin{enumerate}
\item {\it Counted in.} The term $s_r(a)$ of $\inv(a)$ is simply the number
of elements of $E^r(a)$. This is the local role of the exceptional
divisors -- the old exceptional divisors are counted in by $s_r(a)$ in 
order to guarantee that the centre of blowing up lies inside each of them.

\item {\it Factored out.} The new exceptional divisors (those in $\cE_r(a)$)
are factored out from the resolution data on $N_r(a)$ and the next term
$\nu_{r+1}(a)$ is the ``residual multiplicity'' -- the minimal order
of the resolution data (with respect to suitable weighting) after this
factorization. Factoring of the new exceptional divisors guarantees that 
the centre of blowing up (defined {\it a priori} by a local construction)
extends to a global smooth subspace -- see Corollary 1.17 below.
\end{enumerate}

\begin{remark} {\it The invariants of Bierstone-Milman and Villamayor.}
The fact that the equivalence class of a
presentation is strictly smaller than that as a basic object
has an important consequence for the definition of the invariant
(due to the second role of the exceptional divisors): 
The equivalence relation used in \cite{V2, EVacta,
EVweak} does not identify the powers of the new exceptional
divisors as invariants (\S5.3 below). Thinking of a sequence of
blowings-up as a ``history'', the residual multiplicities are
defined in \cite{BMinv} using only the ``future'' (see \S2.4) and in
the articles of Villamayor {\it et al} by a calculation involving
the ``past'' (Lemma 3.7 below). But the latter applies to only
certain subblocks of the ``more recent'' new exceptional
divisors, so only these subblocks are factored out to define
the residual multiplicities in \cite{V2}, etc. 

Each approach has certain
advantages: The notion of invariance introduced in
\cite{BMinv} is stronger (see Remarks 1.16 and \S3.5) and
the desingularization algorithm is in general faster (cf. Example
1.2). But we need Villamayor's invariant in the desingularization
principle for families (Theorem 6.14), which involves a
comparison of the desingularization invariants for a fibre and
for the total space, inductively over the previous history.
\end{remark}

In ``Desingularization algorithms II'', we will show that, for locally
binomial varieties, the desingularization algorithm can be 
greatly simplified: The local role of the exceptional divisors
(1) above is unnecessary because the successive maximal contact
subspaces and exceptional divisors are coordinate subspaces in
suitable local charts; we can use a simpler invariant of the form
$(\io(a), \nu_2(a), \nu_3(a), \ldots)$. 
For {\it affine} binomial varieties,
there is a purely combinatorial algorithm (cf. \cite[Theorem 1.13]
{BMinv}) in which the second role of the exceptional divisors (2)
above is equally superfluous. The general desingularization 
principle has a combinatorial component that is closely related to
the algorithms of ``Desingularization algorithms II''; we will show in
the latter that the techniques involved lead to several natural
questions about the efficiency of desingularization algorithms,
in general.

\begin{example} 
Consider the hypersurface $X$ in $4$-space
defined by
$$
z^dw^{d-1} - x^{d-1}y^d \ = \ 0 \, , 
$$
where $d$ is a positive integer $\geq 2$. The maximum order 
$2d-1$ is taken only at the origin -- this is the first centre
of blowing up in desingularizing $X$. The strict transform $X_1$
of $X$ is defined by
$$
z^d - x^{d-1}y^d \ = \ 0
$$
(in a local coordinate chart where the blowing-up is given
by substituting $(xw, yw, zw, w)$ for the original variables
$(x, y, z, w)$. See Section 4. 
For simplicity of notation, we use the same
variables before and after blowing up.) In Section 4, we estimate
the number of blowings-up needed to reduce the maximum order $d$
of $X_1$. Let $n($BM$)$, $n($V$)$, $n($EV$)$ and $n($LB$)$ denote the
number of blowings-up prescribed by the algorithms of \cite{BMinv},
\cite{V2}, \cite{EVacta} and the locally binomial algorithm of
``Desingularization algorithms II'', respectively. 
Also write $n($AB$)$
for the number of blowings-up given by the affine binomial 
algorithm. Then
\begin{eqnarray*}
n(BM) & \leq & 2d + j\\
n(V) \ = \ n(EV) & \geq & 9d + k\\
n(LB) & \leq & d + l\\
n(AB) & =  & 1 \, , 
\end{eqnarray*}
where $j$, $k$ and $l$ are independent of $d$. The striking
difference between AB and the other algorithms reflects the
fact that a {\it local} invariant that produces a {\it global}
centre of blowing up, as in BM, V, EV or LB, necessarily has
some inefficiency from a purely local point of view.
\end{example}

\begin{example}
Let $X$ denote the surface
$$
z^2 - x^2y^3 \ = \ 0 \, .
$$
The singular locus of $X$ is the union of the $x$- and $y$-axes.
The algorithms of \cite{BMinv} or \cite{V2} prescribe the origin
as the centre $C_0$ of the first blowing-up $\s_1$. The blowing-up
$\s_1$ is given by the substitution $(xy, y, yz)$ in one of 
three coordinate charts, so that the strict transform
$X_1$ of $X$ by $\s_1$ is given in this chart by the same
equation $z^2 - x^2y^3 = 0$ as before; we seem to have accomplished
nothing! But $y = y_\exc$ now defines the {\it exceptional
hypersurface} $H_1 = \s_1^{-1}(C_0)$. The next blowing-up prescribed
by the algorithms again has centre $C_1 = \{0\}$; $\s_2$ is given
by the substitution $(x, xy, xz)$ in one of the charts, so the
strict transform $X_2$ of $X_1$ in this chart is defined by
the equation $z^2 - x^3y^3 = 0$. Here $x$ and $y$ are both exceptional
divisors: $\{y=0\}$ is the strict transform of $H_1$ above and
$\{x=0\} = H_2 := \s_2^{-1}(C_1)$. The two blowings-up $\s_1$ and
$\s_2$ have not simplified the equation, but serve to {\it
re-mark the variables} $x$ {\it and} $y$ {\it as exceptional
divisors}. The exceptional divisors can be thought of as global
coordinates.
\end{example}

In the general algorithm, when a suitable re-marking of variables
as new exceptional divisors has been completed, the simple
combinatorial part of the algorithm takes over to reduce the orders
of functions that are part of the data in a presentation. In
general, when the re-marking has been completed for the algorithm
of \cite{BMinv}, the exceptional variables may not be new in the
sense of \cite{V2} or \cite{EVacta}, so further blowings-up are
needed before the combinatorial part kicks in. This is the difference
in the algorithms highlighted by the example in Section 4.

\subsection{Input data}
The desingularization invariant is a sequence
$\inv(\cdot) = \inv_X(\cdot)$
(or $\inv_\cJ (\cdot)$) as above, beginning with a local invariant
$\io(\cdot) = \io_X(\cdot)$ (or $\io_\cJ(\cdot)$) of a space $X$
(or an ideal sheaf $\cJ$) that distinguishes between general and
special points (e.g., between smooth and singular points of $X$).
(For simplicity of exposition, we sometimes write $\io(\cdot)$ as
$\nu_1(\cdot)$, though, strictly speaking, we reserve $\nu_1(\cdot)$
for the case that $\io(\cdot)$ is the order of an ideal at a point.)
The residual multiplicity $\nu_{r+1}(\cdot)$, $r \geq 1$, is an
invariant of the equivalence class of a local presentation of
the truncated invariant $\inv_r(\cdot) = (\io(\cdot), s_1(\cdot),
\ldots, \nu_r(\cdot), s_r(\cdot))$. A presentation includes a
collection of data defined on a maximal contact subspace $N_r(\cdot)$
of a certain codimension $q$ -- the {\it codimension} of the 
presentation. The maximum locus of the invariant is a union of
global smooth subspaces having only normal crossings. 
The desingularization algorithm is given by choosing as each
successive centre of blowing up a component of the maximum locus
(or of the maximum locus in a suitable subspace). The local
inductive construction allows us to define $(\nu_{r+1}(\cdot), 
s_{r+1}(\cdot))$ using $\inv_r(\cdot)$ and a local presentation
(of codimension $q$, say), and to pass to a local presentation
of $\inv_{r+1}(\cdot)$ of codimension $q+1$. The possibility
of recognizing invariant characteristics of resolution of
singularities by inductive steps of codimension $+1$ appears
already in \cite{BMihes} and is one of the main features
distinguishing \cite{BMmega,BMinv} and \cite{V1,V2} from the
work of Hironaka and Abhyankar. 

Apart from the way that the exceptional divisors are used in
the local construction, there is great flexibility
in the desingularization algorithm depending on choices of the
following:

\begin{enumerate}
\item A notion of transformation by blowing up -- usually strict
or weak transform. (See \S1.4 below.)

\item An initial invariant $\io(\cdot)$; for example, the order of
an ideal $\cJ$ or a space $X$ at a point $a$, or the Hilbert-Samuel
function of $X$ at $a$ (see Examples 1.8).

\item The codimension of a presentation of $\io(a)$.

\item When we ``stop running'' the algorithm.
\end{enumerate}

This flexibility is illustrated by the results in Section 6.

\subsection{The category of spaces}
$\IN$ denotes the nonnegative integers. Throughout this article,
$\IK$ denotes a field of characteristic zero, and $X$ denotes
an algebraic variety or a scheme of finite type over $\IK$, or
an analytic space over $\IK$ (in the case that $\IK$ is locally
compact). For simplicity, we will always assume that our analytic
spaces are compact, or relatively compact in an ambient space
(so that $\inv_X(\cdot)$ will always have maximum values), but
these assumptions are not necessary (see \cite[Section 13]{BMinv}).
We will usually assume that $X$ is embedded as a closed
subspace of a smooth ambient
space (``manifold'') $M$. (See \S6.1, however.) 
Let $\cO_M$
denote the structure sheaf of $M$; $\cO_M$ is a coherent sheaf of
rings. To be precise, $M$ is a local-ringed space
$(|M|, \cO_M)$, where $|M|$
denotes the underlying topological space of $M$.
A closed subspace $X = (|X|, \cO_X)$ of $M = (|M|, \cO_M)$ is
a local-ringed space corresponding to an ideal (i.e., a sheaf
of ideals) $\cJ = \cI_X$ of finite type in $\cO_M$:
$$
|X| \ = \ \supp \cO_M/\cJ \quad\mbox{and}\quad 
\cO_X \ = \ (\cO_M/\cJ) \big|_{|X|} \, .
$$
We say that $X$ is a {\it hypersurface} if $\cI_X$ is principal.

A function in our category (e.g., an element of $\cO_{X,a}$,
where $a \in X$, or of
$\cO_X(U)$, where $U \subset |X|$ is open) will be called
a ``regular function''.

Our desingularization algorithms hold much more generally than
in the categories above, at least in the hypersurface case
(see \cite{BMinv,BMdc}), but we will not pursue this point here.
We do, however, recommend the proof given in \cite{BMdc} as a
motivation for many of the ideas in this article.

\subsection{Weak and strict transform}
Consider an ideal of finite type $\cJ$ in $\cO_M$.

\begin{definitions}
Let $a \in M$. The {\it order} $\mu_a(\cJ)$ of $\cJ$ at $a$
is defined as
$$
\mu_a(\cJ) \ := \ \max \{\mu \in \IN : \, \cJ \subset 
\um^\mu_{M,a} \} \, ,
$$
where $\um_{M,a}$ denotes the maximal ideal of $\cO_{M,a}$.
($\mu_a(\cJ)$ generalizes the {\it order} of a (germ of a)
function $f \in \cO_{M,a}$,
$$
\mu_a(f) \ := \ \max \{\mu \in \IN : \, f \in \um^\mu_{M,a} \} \, .
\mbox{)}
$$
Let $C$ denote a (smooth) subspace of $M$. The {\it order} of 
$\cJ$ {\it along} $C$ at $a$,
$$
\mu_{C,a}(\cJ) \ := \ \max \{\mu \in \IN : \, \cJ \subset
\cI^\mu_{C,a} \} \, .
$$
Thus, $\mu_{C,a}(\cJ)$ is the generic value of $\mu_x(\cJ)$,
for $x \in C$ near $a$.
\end{definitions}

Let $\s: M' \to M$ denote a blowing-up (or a local blowing-up)
with smooth centre $C$. Given an ideal $\cJ$ in $\cO_M$ (or
a closed subspace $X$ of $M$, with ideal sheaf $\cI = \cJ$),
we define {\it weak} and {\it strict transforms} of $\cJ$
(or $X$) by $\s$ as follows:

\begin{definition}
The {\it weak transform} $\cJ'$ of $\cJ$ by $\s$ is the ideal
sheaf $\cJ' \subset \cO_{M'}$ such that, for all $a' \in M'$,
$\cJ'_{a'}$ is the ideal of $\cO_{M',a'}$ generated by 
$$
\{ y^{-\mu}_\exc f \circ \s : \ f \in \cJ_a \} \, ,
$$
where $a = \s(a')$, $\mu = \mu_{C,a}(\cJ)$ and $y_\exc$
denotes a generator of the principal ideal $\cI_{\s^{-1}(C),a'}$.

The {\it weak transform} $X'$ of $X$ is the subspace of $M'$
defined by the weak transform $\cJ'$ of $\cJ := \cI_X$.
\end{definition}

Clearly, the weak transform $\cJ'$ of $\cJ$ by $\s$ is an
ideal of finite type.

\begin{definition} If $X$ is a hypersurface
(i.e., $\cJ := \cI_X$ is principal), then, by definition, the
{\it strict transform} $X'$ of $X$ (or the {\it strict transform}
$\cJ'$ of $\cJ$) by $\s$ coincides with the weak transform.
(If $f \in \cO_{M,a}$ and $a' \in \s^{-1}(a)$, then we will say
that $f' := y_\exc^{-\mu_{C,a}(f)} f \circ \s$ is the ``strict
transform'' of $f$ at $a'$, though of course this is only
well-defined up to an invertible factor.

Now let $X$ be an arbitrary closed subspace of $M$ (i.e.,
$\cJ = \cI_X$ is an arbitrary ideal of finite type in $\cO_M$).
The {\it strict transform} $X'$ of $X$ by $\s$ is defined
locally, at each $a' \in M'$, as the intersection of the
strict transforms of all hypersurfaces containing $X$ at
$a = \s(a')$ (i.e., $\cJ'_{a'} = \cI_{X',a'} \subset \cO_{M',a'}$
is generated by the strict transforms $f'$ of all $f \in \cJ_a)$.
\end{definition}

It is not trivial that the strict transform $X'$ of $X$ is a closed
subspace of $M'$ (i.e., that $\cJ' \subset \cO_{M'}$ is an ideal
of finite type). This depends on Noetherian and coherence properties
of our category of spaces. (For this reason, it is easier to
prove versions of our desingularization theorems that require only
the weak transform, or the strict transform in the hypersurface
case, in categories more general than schemes of finite type or
analytic spaces; cf. \cite{BMdc}.)

\begin{remark}
Let $\cJ \subset \cO_M$ denote an ideal of finite type. Let
$\cJ^s$, $\cJ^w$ and $\cJ^t$ denote the strict, weak and {\it
total} transforms of $\cJ$ (respectively) by a local blowing-up
$\s$. (By definition, $\cJ^t := \s^{-1}(\cJ)$.) Then
$$
\cJ^s \ \subset \ \cJ^w \ \subset \ \cJ^t
$$
and $\cJ^s = \cJ^t$ if and only if $C \not\subset \supp \cJ$.
\end{remark}

Athough the notions of weak and strict transform both apply to
ideals $\cJ$ or subspaces $X$, we will usually refer to the 
``strict transform $X'$ of $X\, $'' or the ``weak transform $\cJ'$
of $\cJ\, $'' in order to economize notation and to emphasize the
role of the strict transform in embedded resolution of
singularities and of the weak transform
in principalization of an ideal.

\subsection{Desingularization invariants}
Given a closed subspace $X$ of $M$, or an ideal of finite type
$\cJ \subset \cO_M$, we will consider sequences of transformations
\begin{equation}
\begin{array}{rccccccccccl}
\longrightarrow & M_{j+1} & \stackrel{\s_{j+1}}{\longrightarrow}
& M_j & \longrightarrow & \cdots & \longrightarrow & M_1 & \stackrel
{\s_1}{\longrightarrow} & M_0 & = & M \\
& X_{j+1} & & X_j & & & & X_1 & & X_0 & = & X\\
& E_{j+1} & & E_j & & & & E_1 & & E_0 & = &E 
\end{array}
\end{equation}
or sequences of transformations
\begin{equation}
\begin{array}{rccccccccccl}
\longrightarrow & M_{j+1} & \stackrel{\s_{j+1}}{\longrightarrow}
& M_j & \longrightarrow & \cdots & \longrightarrow & M_1 & \stackrel
{\s_1}{\longrightarrow} & M_0 & = & M \\
& \cJ_{j+1} & & \cJ_j & & & & \cJ_1 & & \cJ_0 & = & \cJ\\
& E_{j+1} & & E_j & & & & E_1 & & E_0 & = &E
\end{array}
\end{equation}
where, in each case, $E$ is a finite collection of smooth 
hypersurfaces in $M$ having only normal crossings (usually
$E = \emptyset$) and, for each $j$:
\begin{list}{}{}
\item[$\s_{j+1}$] is a (local) blowing-up of $M_j$ with smooth
centre $C_j$ such that $E_j$ and $C_j$ simultaneously have
only normal crossings.
\item[$E_{j+1}$] 
$:= E'_j \cup \{\s_{j+1}^{-1}(C_j)\}$, where
$E'_j$ denotes the collection of strict transforms $H'$ of all
hypersurfaces $H \in E_j$. (Thus, $E_{j+1}$ has only normal
crossings.)
\item[$X_{j+1}$] denotes the strict transform $X'_j$ of $X_j$
by $\s_{j+1}$\, , or
\item[$\cJ_{j+1}$] denotes the weak transform $\cJ'_j$ of $\cJ_j$
by $\s_{j+1}$.
\end{list}

Let $\Sig$ denote a partially ordered set. Let $\io_X$ denote
a local invariant of $X$ with values in $\Sig$; i.e., a function
$\io_X: X \ni a \mapsto \io_X(a) \in \Sig$ such that $\io_X(a)$
depends only on the local isomorphism class of $X$ at $a$. (Or
let $\io_\cJ: M \to \Sig$ denote a local invariant of $\cJ$. We
can also consider $\io_X$ to be defined on $M$.) Write $\io
= \io_X$ or $\io_\cJ$, to cover both cases.

\begin{examples} (1) $\io_{\cJ}(a) = \mu_a(\cJ)$, 
the order of $\cJ$ at $a$.
\smallskip

(2) $\io_X(a) = \nu_X(a)$, where $\nu_X(a) := \mu_a(\cI_X)$,
the {\it order} of $X$ at $a$.
\smallskip

(3) $\io_X(a) = H_{X,a}$, the {\it Hilbert-Samuel function}
of $X$ at $a$; i.e., the function
$$
H_{X,a}(l) \ = \ \dim_\IK \frac{\cO_{X,a}}{\um_{X,a}^{l+1}},
\quad l \in \IN ,
$$
where $\um_{X,a}$ denotes the maximal ideal of $\cO_{X,a}$.
(In the case of schemes, this definition is correct as stated
only at a $\IK$-rational point $a$; we should otherwise replace
$\dim_\IK$ by length.) The order $\nu_X(a)$ is determined by
$H_{X,a}$; if $X$ is a hypersurface, then $H_{X,a}$ is determined
by $\nu_X(a)$. See \cite[Rmks. 1.3]{BMinv} 
for details of these remarks.
\end{examples}

\begin{definition}
A local blowing-up $\s: M' \to M$ with centre $C$ is 
$\io$-{\it admissible} if $\io$ is locally constant on $C$.
\end{definition}

\begin{hypotheses} We will assume that $\io$ satisfies the
following three properties:
\begin{enumerate}
\item {\it Semicontinuity.}
\begin{enumerate}
\item $\io$ is upper-semicontinuous in the Zariski topology.
\item $\io$ is {\it infinitesimally upper-semicontinuous} in
the sense that, if $\s: M' \to M$ is an $\io$-{admissible}
local blowing-up, then $\io'_{a'} \leq \io_a$, for all $a'
\in \s^{-1}(a)$.
\end{enumerate}
\item {\it Stabilization.} Every decreasing sequence in the
value set of $\io$ stabilizes.
\item $\io$ admits a {\it semicoherent presentation} at every
point, in the sense of 1.12 following.
\end{enumerate}
\end{hypotheses}

\begin{remark} In an analytic category, where we can distinguish
between the classical and Zariski topologies, we will always understand
by (1)(a) above, the following somewhat weaker property: Every
point admits a classical neighbourhood in which $\io$ is
Zariski upper-semicontinuous. Neighbourhoods and germs, in the
analytic case, will be
understood to be in the Zariski sense, but within some classical
coordinate neighbourhood of a given point. This is important,
for example, in 1.12, \S2.5 and Section 3 below (see Remark 3.2), 
but we will not labour the point
here; we refer to \cite{BMinv} for more details.
\end{remark}

\begin{definitions and remarks} Let
$$
S_\io(a) \ := \ \{x \in M: \, \io(x) \geq \io(a)\}
$$
as a germ at $a$ (so that $S_\io(a) := \{x \in M: \, 
\io(x) = \io(a)\}$, as germs, by property (1)(a) above). We
call $S_\io(a)$ (the germ of) the {\it constant locus} 
of $\io$ at $a$. Let
$$
\ucG(a) \ = \ ( N(a), \cG(a)) ,
$$
where $N(a)$ denotes a germ of a submanifold of $M$ at $a$,
of codimension $p$, say, and $\cG(a)$ is a finite collection
of pairs $(g, \mu_g)$, where $g \in \cO_{N(a)}$, $\mu_g \in \IQ$
and $\mu_a(g) \geq \mu_g$. We define the {\it equimultiple locus}
of $\ucG(a)$,
$$
S_{\ucG(a)} \ := \ \{x \in N(a): \ \mu_x(g) \geq \mu_g,
\ \mbox{for all}\ (g, \mu_g) \in \cG(a)\}.
$$
This makes sense as a germ at $a$.

We say that $\ucG(a)$ is a {\it presentation} of $\io$ of
{\it codimension} $p$ at $a$ if:
\begin{enumerate}
\item $S_{\ucG(a)} = S_\io(a)$.
\item If $\s: M' \to M$ is an $\io$-admissible local blowing-up
(with centre $C$) and $a' \in \s^{-1}(a)$, then $\io(a') =
\io(a)$ if and only if $a' \in N(a')$ and 
$$
\mu_{a'}(y_{\exc}^{\mu_g} g \circ \s) \ \geq \ \mu_g, \quad
\mbox{for all}\ (g, \mu_g) \in \cG(a),
$$
where $N(a')$ denotes the germ at $a'$ of the strict transform
$N(a)'$ of $N(a)$, and $y_{\exc}$ denotes a local generator of
the ideal of the {\it exceptional hypersurface} $\s^{-1}(C)$.
\item The preceding properties (1) and (2) are stable (i.e.,
continue to hold) after suitable finite sequences of three
kinds of morphisms (admissible blowing-up, product with a line,
and exceptional blowing-up) and corresponding transformations of
$X$ (or $\cJ$).
\end{enumerate}

The notion of presentation will be made more precise in
Section 2. Exceptional blowings-up are defined in terms
of the new exceptional divisors (\S2.2; cf. \S1.1 above).
Property (3) concerns sequences of ``test blowings-up'' used
to prove invariance of $\inv$, and is used to define equivalence
of presentations (\S2.3). The corresponding stability
property for idealistic exponents or basic objects involves
only admissible blowings-up and product with a line.

We will usually identify a germ with a representative in a
small neighbourhood. Then $\ucG(a) \, = \, (N(a), \cG(a))$
induces a pair $\ucG(x) \, = \, (N(x), \cG(x))$, at any $x \in
S_{\ucG(a)}$ near $a$. A presentation $\ucG(a)$ of $\io$ at $a$
is called {\it semicoherent} if $\ucG(x)$ is a presentation of
$\io$ at $x$, for each $x$ in a neighbourhood of $a$ in $S_\io(a)
= S_{\ucG(a)}$. The notion of semicoherent presentation at $a$
clearly depends only on $\ucG(a)$ and the germ of $X$ (or the
stalk of $\cJ$) at $a$. 
\end{definitions and remarks}

\begin{examples} (1) {\it The order $\mu_{\cdot}(\cJ)$ of an 
ideal $\cJ$ satisfies Hypotheses 1.10:} Properties (1)(a) and (2) 
are obvious. (1)(b) is an elementary Taylor series computation
(cf. \cite[lemma 5.1]{BMinv}). Let $\ucG(a) = \{(g, \mu_g)\}$,
where the $g$ form any finite set of generators of $\cJ_a$,
and $\mu_g = \mu_a(\cJ)$. Then $\ucG(a)$ provides a codimension zero
semicoherent presentation of $\mu_{\cdot}(\cJ)$ at $a$. 
(See \cite[Prop. 6.5]{BMinv}.)
\smallskip

(2) {\it The Hilbert-Samuel function satisfies Hypotheses 1.10:}
Properties (1) and (2) are due to Bennett \cite{Ben}. See
\cite[Thms. 9.2, 7.20]{BMinv} and \cite[Thm. 5.2.1]{BMjams}
for elementary proofs of (1)(a),(b) and (2), respectively.
\cite[Thms. 9.4, 9.6]{BMinv} provide a semicoherent presentation
of the Hilbert-Samuel function.

In Theorem 6.17, we show that the Hilbert-Samuel function
satisfies stronger hypotheses (6.12 below) that are needed
to extend the desingularization principle to parametrized
families.
\end{examples}

The purpose of the local inductive construction is to extend
$\io(a) = \io_X(a)$ or $\io_\cJ(a)$, $a \in M_0 = M$, to an
``invariant'' $\inv(a) = \inv_X(a)$ or $\inv_\cJ(a)$ , where
$a \in M_j$, $j = 0,1,\ldots$, defined recursively over a 
sequence of (local) blowings-up (1.1) or (1.2) provided that,
for all $i \leq j$,
$\s_i$ is $\inv${\it -admissible}; i.e., $\inv(\cdot)$
is locally constant on $C_i$.

In other words, $X$ or $\cJ$ determines $\inv_X(a)$ or $\inv_\cJ(a)$
for $a \in M_0 = M$, and thus the first centre of blowing up 
$C_0 \subset M_0$; then $\inv(a) = \inv_X(a)$ or $\inv_\cJ(a)$
can be defined on $M_1$ and determines $C_1$, etc.

The notation $\inv_X(a)$ or $\inv_\cJ(a)$, where $a \in M_j$,
indicates a dependence on the original space $X_0 = X$ or ideal
$\cJ_0 = \cJ$, and not simply on $X_j$ or $\cJ_j$. Some dependence
on the history of the desingularization process (1.1) or (1.2)
is necessary in order to determine a global centre of blowing up
using a local invariant -- see Example 1.3 above and \cite[Example
1.9]{BMmsri}. Corollary 1.17 below shows how $\inv$ is used to
determine a global centre.

The invariant $\inv(a) = \inv_X(a)$ or $\inv_\cJ(a)$, $a \in M_j$,
is a finite sequence
\begin{equation}
\inv(a) \ = \ (\io(a), s_1(a), \nu_2(a), s_2(a), \ldots , 
\nu_{t+1}(a))
\end{equation}
beginning with $\io(a) = \io_X(a)$ or $\io_\cJ(a)$. The terms
$s_r(a) \in \IN$ are defined in 1.15 below. For each $r = 1,\ldots,$
the residual multiplicity $\nu_{r+1}(a) \in \IQ$ is defined using
a presentation of $\inv_r$ of codimension $q+r$ at $a$ (where we
begin with a presentation of $\inv_{1/2} = \io$ of codimension
$q+1$, say). We can then define a presentation of $\inv_{r+1}$
of codimension $q+r$, and the local inductive contruction allows
us to pass to an equivalent presentation of codimension $q+r+1$,
to complete a cycle in the inductive definition. (See \S\S3.2 
and 6.1 below.)
The construction terminates by exhaustion of variables; if
$n = \dim_a M_j$, then $t \leq n$ and $\nu_{t+1}(a) = 0$ or
$\infty$. Sequences in the value set of $\inv$ can be compared
lexicographically.

Although the residual multiplicities are rational, their 
denominators are controlled in the following way: There is 
$e_1 \in \IN$ (for example, $e_1 = \io(a)$ in the case that
$\io(a)$ is the order of an ideal) such that, for
all $r > 0$, $e_r!\nu_{r+1}(a) \in \IN$, where $e_{r+1} = 
\max\{e_r!, e_r!\nu_{r+1}(a)\}$. (See \S3.3(2).)

\begin{theorem} Let $\io = \io_X$ or $\io_\cJ$ denote an invariant
of $X$ or $\cJ$, satisfying the hypotheses 1.10 above. Then $\io$
extends to an invariant $\inv = \inv_X$ or $\inv_\cJ$ which is
defined over any sequence of tranformations (1.1) or (1.2), where
the successive (local) blowings-up are $\inv$-admissible, having
the following properties:
\begin{enumerate}
\item{\it Semicontinuity.}
\begin{enumerate}
\item
$\inv$ is Zariski upper-semicontinuous.
\item
$\inv$ is infinitesimally upper-semicontinuous; i.e., 
$\inv(a) \leq \inv(\s_j(a))$, for all $a \in M_j$, $j \geq 1$.
\end{enumerate}
\item{\it Stabilization.} 
If $a_j \in M_j$ and $a_j = \s_{j+1}(a_{j+1})$,
$j = 0, 1, \ldots ,$ then there exists $j_0$ such
that $\inv(a_j) = \inv(a_{j+1})$, $j \geq j_0$.
\item{\it Normal crossings.} Let $a \in M_j$. Then $S_\inv (a)$
and $E(a)$ simultaneously have only normal crossings. If 
$\inv(a) = \inv_{t+1/2}(a) = (\ldots, \infty)$, 
then $S_\inv (a)$ is smooth. If $\inv(a) = \inv_{t+1/2}(a)
= (\ldots, 0)$, then each component $Z$ of $S_\inv (a)$
is of the form 
\begin{equation}
Z \ = \ S_\inv (a) \, \cap \, \bigcap\{H \in E(a):\ Z \subset H\}.
\end{equation}
\item{\it Decrease.} Let $a \in M_j$. If $\inv(a) = (\ldots, \infty)$
and $\s$ is the local blowing-up of $M_j$ with centre $S_\inv (a)$,
then $\inv(a') < \inv(a)$ for all $a' \in \s^{-1}(a)$. On the other
hand, suppose that $\inv(a) = \inv_{t+1/2}(a) = (\ldots, 0)$. 
Then there is an
additional invariant $\mu(a) = \mu_X(a)$ or $\mu_\cJ(a) \in \IQ$,
$\mu(a) \geq 1$, such that, if $Z$ denotes any component of 
$S_\inv (a)$ and $\s$ is the local blowing-up with centre $Z$,
then
$$
(\inv(a'), \mu(a')) \ < \ (\inv(a), \mu(a)),
$$
for all $a' \in \s^{-1}(a)$. ($e_t!\mu(a) \in \IN$.)
\end{enumerate}
\end{theorem}

\subsubsection*{Old exceptional divisors.} The terms $s_r(a)$
in $\inv(a)$ can be defined immediately, in an invariant way:
If $a \in M_j$, write $a_i$ to denote the image of $a$ in $M_i$,
where $i \leq j$; i.e., $a_j = a$ and $a_i = (\s_{i+1}\circ\cdots
\circ\s_j)(a)$ if $i < j$.

\begin{definitions}
Let $a \in M_j$. Let $i$ denote the earliest year in the resolution
history (1.1) or (1.2) where $\io(a) = \io(a_i)$ (i.e., $i$ is the
smallest index $i'$ such that $\io(a) = \io(a_{i'})$. We 
sometimes call $i$ the ``year of birth'' of $\io(a)$.) Set
\begin{eqnarray*}
E^1(a) & := & \{H \in E(a):\ H\mbox{ is transformed from }E(a_i)\},\\
s_1(a) & := & \#E^1(a).
\end{eqnarray*}

To define $s_{r+1}(a)$ in general, let $i$ be the earliest year
where $\inv_{r + 1/2}(a) = \inv_{r + 1/2}(a_i)$. Set
\begin{eqnarray*}
E^{r+1}(a) & := & \{H \in E(a)\backslash \left( E^1(a)\cup\cdots\cup
E^r(a) \right) :\ H\mbox{ is} \\
	   &    & \qquad\qquad\qquad \mbox{transformed from }E(a_i)\},\\
s_{r+1}(a) & := & \#E^{r+1}(a).
\end{eqnarray*}
\end{definitions}

\subsubsection*{Meaning of invariance.} $\inv(a)$ depends only
on the local isomorphism class of $(X_{j,a}\ (\mbox{or } \cJ_{j,a}),
E(a), E^1(a),\ldots, E^t(a))$ -- see \S\S3.2, 3.5. 

\begin{remarks}
The proof of \cite{V2} or \cite{EVacta} gives a weaker sense
of invariance; see \S3.5. Although \cite{BMinv} and \cite{V2,EVacta}
count in old exceptional divisors in the same way, the latter
factor out only a subset of the new exceptional divisors $\cE^r(a)$
(those accumulating after the birth of $\inv_r(a)$) to define the
residual multiplicities $\nu_{r+1}(a)$. The resulting difference
in the invariants can first show up in the $\nu_2$-terms. But
$s_2$ depends on $\inv_{3/2}$, $\nu_3$ depends of $\inv_2$ and
the next block of new exceptional divisors, etc., so the difference
is magnified in each of the following terms. 
\end{remarks}

\subsection{Desingularization algorithm}
The following corollary of Theorem 1.14 above captures the
role of the exceptional divisors in the passage from local
to global. (The exceptional divisors involved in (1.4) above
belong to the new block $\cE^t(a)$.)

\begin{corollary}
Let $a \in M_j$, and consider an open neighbourhood $U$ of
$a$ in $M_j$ such that $\inv(a)$ is a maximum value of 
$\inv(\cdot)$ on $U$. Then each component $Z$ of $S_\inv(a)$
extends to a global smooth closed subspace of $U$.
\end{corollary}

\begin{proof} Consider any total order on $\{I:\ I \subset E_j\}$.
Let $a \in M_j$. We label each component $Z$
of $S_\inv(a)$ as $Z_I$, where $I := \{H \in E(a):\ Z \subset H\}$.
Define
\begin{eqnarray*}
J(a) & := & \max\{I:\ Z_I \mbox{ is a component of } S_\inv(a)\},\\
\inv^e(a) & := & (\inv(a), J(a)).
\end{eqnarray*}
It is easy to see that $\inv^e(\cdot)$ is Zariski upper-semicontinuous
on $M_j$ and its maximum locus in any open subspace of $M_j$
is smooth.

Given $a \in M_j$ and a component $Z_I$ of $S_\inv(a)$, we can
choose the order above so that $I = J(a) = \max\{J:\ J \subset E_j\}$.
Then $Z_I$ extends to a smooth closed subspace of the open set
$\{x \in M_j:\ \inv(x) \leq \inv(a)\}$.
\end{proof}

We can order $\{I:\ I \subset E_j\}$ using the resolution
history (1.1) or (1.2) as follows: (Assuming $E_0 = \emptyset$),
write $E_j = \{H^j_1,\ldots,H^j_j\}$, where $H^j_i$ is the
strict transform of $H^{j-1}_i$ by $\s_j$, $i = 1,\ldots,j-1$,
and $H^j_j = \s_j^{-1}(C_{j-1})$. Associate to each $I \subset E_j$
the sequence $(\de_1,\ldots,\de_j)$, where $\de_i = 0$ if $H^j_i
\notin I$ and $\de_i = 1$ if $H^j_i \in I$, and use the lexicographic
ordering of such sequences, for all $j$ and $I \subset E_j$.
Consider the extended invariant $\inv^e(\cdot) := (\inv(\cdot), 
J(\cdot))$ defined using this ordering. 
The desingularization principle Theorem 1.14 above can then 
be applied as follows:

\begin{algorithm}
Blow up with each successive centre given by the maximum locus
of $\inv^e$ (or the maximum locus within a suitable closed
subspace, depending on the problem). Stop when $\inv$ is locally
constant (on a suitable subspace).
\end{algorithm}

\begin{examples}(1) {\bf Principalization of an ideal} $\cJ$ 
(cf. \cite[Main Thm. II]{Hann}, \cite[Thm. 1.10]{BMinv}).
Take $\io(a) = \io_\cJ(a) = \mu_a(\cJ)$. Then there is a finite
sequence of blowings-up (1.2) with $\inv_{\cJ}$-admissible
centres, such that, if $\cJ' \subset \cO_{M'}$ denotes the
final weak transform of $\cJ$, then $\cJ' = \cO_{M'}$ and 
$\s^{-1}(\cJ) = \s^*(\cJ) \cdot \cO_{M'}$ is a normal-crossings
divisor, where $\s:\ M' \to M$ denotes the composite of the
sequence of blowings-up. 
\smallskip

\noindent
{\it Algorithm:} Blow up
with successive centres given by the maximum locus of $\inv_\cJ^e$
in $M_j$. Stop when $\supp \cO_{M_j}/\cJ_j = \emptyset$; i.e.,
$\cJ_j = \cO_{M_j}$.
\medskip

(2) {\bf Embedded desingularization} (cf. \cite[Main Thm. I]{Hann},
\cite[Theorem 1.6]{BMinv}, \cite{V2,EVacta}).
The following ``geometric version'' is meaningful if the set
of smooth points of $X$ is not empty. Take $\io = \io_X = H_{X,\cdot}$.
Then there is a finite sequence of blowings-up (1.1) with
$\inv_X$-admissible centres $C_j$, such that:
\begin{enumerate}
\item[(a)] For each $j$, either $C_j \subset \Sing X_j$ or
$X_j$ is smooth and $C_j \subset X_j \cap E_j$.
\item[(b)] Let $X'$ and $E'$ denote the final strict transform of $X$
and exceptional set, respectively. Then $X'$ is smooth and
$X', E'$ simultaneously have only normal crossings.
\end{enumerate}
\smallskip

\noindent
{\it Algorithm:}
\begin{list}{}{}
\item[{\it Step 1.}] Blow up with successive centres given by the
maximum locus of $\inv_X^e$ in $\Sing X_j$. Stop when $\Sing X_j = 
\emptyset$.
\item[{\it Step 2.}] Continue to blow up with successive centres
given by the maximum locus of $\inv_X^e$ in $S_j := \{x \in M_j:
\ s_1(x) > 0\}$. Stop when $S_j = \emptyset$.
\end{list}
\medskip

(3) Embedded desingularization in the nonreduced case. See
\cite[Section 11]{BMinv}.
\medskip

(4) {\bf Weak embedded desingularization} (cf. \cite{EVweak}).
Take $\io(a) = \io_\cJ(a) = \mu_a(\cJ)$, where $\cJ = \cI_X$.
Then there is a finite sequence of blowings-up (1.2) with
$\inv_{\cJ}$-admissible centres $C_j$, such that:
\begin{enumerate}
\item[(a)] Each $C_j \subset \pi_j^{-1}(\Sing X)$, 
where $\pi_j$ denotes
the composite of the blowings-up to year $j$.
\item[(b)] $X'$ is smooth and $X', E'$ simultaneously have only
normal crossings (in the notation of (2)(b)).
\end{enumerate}

\noindent
{\it Algorithm:}
(Assume that $X$ is pure-dimensional.) Apply
the algorithm for principalization of $\cJ = \cI_X$, but stop
early -- when the maximum value of $\inv_\cJ$ becomes equal
to the generic value of $\inv_\cJ$ on $X_0 = X$. See \S6.2.
\end{examples}

Other applications of the desingularization principle are
given in \cite[Chapter IV]{BMinv} and in Section 6 below.

\begin{remark}
``Weak embedded desingularization'' (4) above is weaker than
``embedded desingularization'' (2) because, in (4), it is not
in general true that the successive centres $C_j$ lie in the
strict transforms $X_j$, nor that the $C_j \cap X_j$ are
smooth. ({\it Any} birational projective morphism of 
quasiprojective varieties is a blowing-up with (not necessarily
smooth) centre.)
\end{remark}

\section{Idea of a presentation}

The desingularization
invariant $\inv = \inv_X$ or $\inv_\cJ$ is to be defined 
recursively over a sequence
of admissible blowings-up (1.1) or (1.2).
The successive pairs $(\nu_r,s_r)$ in the sequence
$\inv$ will themselves be defined inductively.  Given $\inv_{r-1/2}$,
$r \ge 1$, $s_r$ and $\nu_{r+1}$ can be defined recursively over
any sequence of $(r-1/2)$-admissible (i.e., 
$\inv_{r-1/2}$-admissible) (local) blowings-up.

Consider a sequence of $(r-1/2)$-admissible (local) blowings-up
(1.1) or (1.2).
Let $a \in M_j$.
Assume that $\nu_r(a) \ne 0,\infty$.
We have defined $s_r(a)$ in Definitions 1.15.
The following term $\nu_{r+1}(a)$ of $\inv(a)$ can be defined using a
``presentation of $\inv_r$ at $a$''.
A presentation of codimension $p$ involves a collection of regular
functions (i.e., functions in our category) 
with ``assigned multiplicities'' on a ``maximal contact
subspace'' of codimension $p$.  A presentation has an ``equimultiple
locus'' (as a germ at $a$); cf. 1.12.

A presentation is not an invariant.
We introduce a notion of equivalence of presentations;
The equivalence class of a presentation of $\inv_r$ at $a$ does
have an invariant meaning, and $\nu_{r+1}(a)$ depends only on
the equivalence class. See Theorems 2.3, 2.4 and \S3.2.
It is convenient, in fact, to consider two notions of equivalence:
\begin{enumerate}
\item A purely local notion; see \S2.3.
The corresponding equivalence class will be denoted $[\cdot]$.
Two presentations at $a$ are equivalent in this sense if they
have the same equimultiple locus, both at $a$ and also after
certain sequences of transformations (cf. 1.12).

\item A stronger notion of ``semicoherent equivalence''.
Two presentations at $a$ are semicoherent equivalent if they induce
presentations that are equivalent in the sense of (1) at each point
of the equimultiple locus of $a$.
(This again makes sense as a notion about germs.)
See \S2.5. 
The semicoherent equivalence class will be denoted $[\![\cdot]\!]$.
\end{enumerate}

We differentiate between (1) and (2) in order to explain precisely
on what the successive terms of $\inv$ depend.

The idea of a presentation is treated in an abstract way
in this section.
In Section 3 below, we show how the idea is
used to define desingularization invariants and to prove the
desingularization principle Theorem 1.14.

\subsection{Definition of a presentation {\rm (cf. 
\cite[(4.1)]{BMinv}}} Let $M$ denote a manifold and let $a \in M$.
A ({\it local}) {\it presentation} of {\it codimension} $p$ at $a$
is a triple
$$
\ucH(a) \ = \ ( N(a),\cH(a),\cE(a)) ,
$$
where:
\begin{list}{}{}
\item[$N(a)$] is a germ of a submanifold of codimension $p$
at $a$;

\item[$\cH(a)$] $= \{ (h,\mu_h)\}$ is a finite collection of pairs
$(h,\mu_h)$, where $h\in\cO_{N(a)}$, $\mu_h\in\IQ$ and
$\mu_a(h)\ge\mu_h$.

\item[$\cE(a)$] is a finite collection of smooth hypersurfaces
such that $N(a)$, $\cE(a)$ simultaneously have only normal
crossings, and $N(a)\not\subset H$, for all $H\in\cE(a)$.
\end{list}

We define the {\it equimultiple locus} $S_{\ucH(a)}$ of
$\ucH(a)$ (as a germ at $a$) by
$$
S_{\ucH(a)} \ := 
\ \{ x\in N(a):\ \mu_x(h) \ge \mu_h,\ \mbox{for all}
\ (h, \mu_h) \in \cH(a)\} .
$$

\subsection{Transforms of a presentation} Our notion of equivalence 
of presentations is given by stability of the condition that their
equimultiple loci coincide, after sequences of transformations
by three types of morphisms \cite[Section 4]{BMinv}. 
Let $\ucH(a)=( N(a),\cH(a),\cE(a))$ denote a local
presentation of codimension $p$ at $a$. 
We will consider transformations of $\ucH(a)$ by 
morphisms $\s$ of the following three types.
\begin{enumerate}
\item[(i)] {\it Admissible blowing-up}: a local blowing-up $\s$
with centre $C$ such that $C$ and $\cE(a)$ simultaneously have
only normal crossings, and $a\in C\subset S_{\ucH(a)}$.

\item[(ii)] {\it Product with a line}: $\s$ is a projection $M'=W\times\IA^1
\to W\hookrightarrow M$ over a neighbourhood $W$ of $a$.

\item[(iii)] {\it Exceptional blowing-up}: a local blowing-up $\s$ at $a$
with centre $C=H_0\cap H_1$, where $H_0, H_1\in \cE(a)$.
\end{enumerate}

For any of the morphisms $\s$ above, we define a transform
$$
\ucH(a') \ = \ ( N(a'),\cH(a'),\cE(a'))
$$
of $\ucH(a)$ {\it at certain points} $a'\in\s^{-1} (a)$:

\begin{enumerate}
\item[(i)] Admissible blowing-up: Suppose that $a'\in\s^{-1}(a)$
is a point such that
$$
\mu_{a'}(y_\exc^{-\mu_h} h\circ\s) \ \ge \ \mu_h, \quad
\mbox{for all } (h,\mu_h)\in\cH(a)
$$
(where $y_\exc$ denotes a local
generator of the ideal of $\s^{-1}(C)$).
Then we define:
\begin{list}{}{}
\item[$N(a')$] $:=$ germ at $a'$ of the strict transform 
$N(a)'$ of $N(a)$,

\item[$\cH(a')$] $:= \left\{ (h',\mu_{h'}):\ (h,\mu_h)\in\cH(a)\right\}$,
where $h' :=$ germ at $a'$ of $y_\exc^{-\mu_h}h\circ\s$ and
$\mu_{h'} := \mu_h$,

\item[$\cE(a')$] $:= \left\{ H':\ H\in \cE(a),\ a'\in H'\right\}
\cup \left\{\s^{-1} (C)\right\}$ (where $H'$ means the strict transform
of $H$).
\end{list}
\item[(ii)] Product with a line.
Let $a'=(a,0)$.
Then we define:
\begin{list}{}{}
\item[$N(a')$] $:=$ germ at $a'$ of $\s^{-1}\left( N(a)\right)$,

\item[$\cH(a')$] $:= \left\{ (h\circ\s,\mu_h):\ (h,\mu_h)\in\cH(a)\right\}$,

\item[$\cE(a')$] $:= \left\{ \s^{-1} (H):\ H\in\cE(a)\right\} \cup \{ W\times 0\}$.
\end{list}
\item[(iii)] Exceptional blowing-up.
Let $a'$ be the unique point of $\s^{-1}(a)\cap H'_1$.
Then we define $N(a')$, $\cH(a')$ as in (ii), and $\cE(a')$
as in (i).
\end{enumerate}

\subsection{Equivalence of presentations {\rm \cite[Section 4]{BMinv}}}

\begin{definitions}
Two presentations
\begin{eqnarray*}
\ucH(a) & = & \left( N=N(a),\, \cH(a),\, \cE(a)\right) \\
\ucF(a) & = & \left( P=P(a),\, \cF(a),\, \cE(a)\right) ,
\end{eqnarray*}
perhaps of different codimension but with common $\cE(a)$,
are {\it equivalent with respect to transformations of types}
(i), (ii) if:
\begin{enumerate}
\item $S_{\ucH(a)} = S_{\ucF(a)}$.

\item If $\s$ is an admissible blowing-up and $a'\in \s^{-1}(a)$,
then $a'\in N'$ and $\mu_{a'}(y_\exc^{-\mu_h} h\circ\s)\ge\mu_h$,
for all $(h,\mu_h)\in\cH(a)$, if and only if $a'\in P'$ and
$\mu_{a'}(y_\exc^{-\mu_f}f\circ\s)\ge\mu_f$, for all $(f,\mu_f)\in
\cF(a)$.

\item Conditions (1) and (2) continue to hold for 
$\ucH(a')$, $\ucF(a')$
obtained by any sequence of transformations by morphisms of types
(i), (ii).
\end{enumerate}

We likewise define {\it equivalence with respect to transformations
of types} (i), (ii), (iii) (by simply using all three types of
morphisms in (3) above).
We write $[\ucH(a)]_{{\rm (i,ii)}}$ and 
$[\ucH(a)]_{{\rm (i,ii,iii)}}$
for the corresponding equivalence classes.
\end{definitions}

\begin{definition}
We introduce an intermediate notion of equivalence by allowing
only certain sequences of morphisms (i), (ii) and (iii) in
Definitions 2.1 above; namely,
\begin{equation*}
\begin{array}{rccccccccccl}
\longrightarrow & M_j & \stackrel{\s_j}{\longrightarrow}
& \cdots & \longrightarrow &
M_i & \longrightarrow & \cdots & \longrightarrow & M_0 & = & M \\
& \cE(a_j) &  & & & \cE(a_i) & & & 
& \cE(a_0) & = & \cE(a) 
\end{array}
\end{equation*}
where, if $\s_{i+1},\ldots,\s_j$ are exceptional blowings-up,
then $i\ge 1$ and $\s_i$ is either of type (iii) or (ii).
In the latter case, $\s_i$: $M_i=M_{i-1}\times\IA^1\to M_{i-1}$ is
the projection and each $\s_{k+1}$, $k=i,\ldots,j-1$, is the blowing-up
with centre $C_k = H_0^k \cap H_1^k$, where $H_0^k$, $H_1^k\in
\cE(a_k)$, $a_{k+1}=\s_{k+1}^{-1} (a_k)\cap H_1^{k+1}$, and the
$H_0^k$, $H_1^k$ are determined by some fixed $H\in\cE(a_{i-1})$
inductively in the following way:
$H_0^i := M_{i-1}\times \{0\}$, $H_1^i := \s_i^{-1} (H)$,
and, for each $k=i+1, \ldots, j-1$ $H_0^k := \s_k^{-1} (C_{k-1})$,
$H_1^k :=$ the strict transform of $H_1^{k-1}$ by $\s_k$.
\end{definition}

Let $[\ucH(a)]$ denote the equivalence class corresponding to
Definition 2.2.
Clearly,
\begin{equation}
[\ucH(a)]_{{\rm (i,ii,iii)}} \ \subset \ [\ucH(a)] \ \subset 
\ [\ucH(a)]_{{\rm (i,ii)}} .
\end{equation}

\subsection{Invariants of a presentation}
Let $\ucH(a)=( N(a),\cH(a),\cE(a))$ be a 
presentation at $a$;
say $\cH(a)=\{ (h,\mu_h)\}$.
We define
\begin{equation*}
\begin{array}{rcccl}
\mu(a) & = & \mu_{\ucH(a)} & := & \min_{\cH(a)}
\displaystyle\frac{\mu_a(h)}{\mu_h} ,\\
\mu_H(a) & = & \mu_{\ucH(a),H} & := & \min_{\cH(a)}
\displaystyle\frac{\mu_{H,a}(h)}{\mu_h} , \quad H\in\cE(a),\\
\nu(a) & = & \nu_{\ucH(a)} & := & \mu_{\ucH(a)} - 
\displaystyle\sum_{H \in \cE(a)}
\mu_{\ucH(a),H}. 
\end{array}
\end{equation*}
(Recall Definitions 1.4.)
In particular, $0\le\nu(a)<\infty$ if $\mu(a)<\infty$.
(We set $\nu(a)=\infty$ if $\mu(a)=\infty$.)

\begin{theorem}
Let $\ucH(a)$ and $\ucF(a)$ denote presentations 
that are equivalent with respect to transformations
of types (i) and (ii). If $\ucH(a)$ and $\ucF(a)$ have
the same codimension, then
$$
\mu_{\ucH(a)} \ = \ \mu_{\ucF(a)} .
$$
If $\codim \ucH(a) > \codim \ucF(a)$, then $\mu_{\ucF(a)} = 1$.
\end{theorem}

\begin{theorem}
If $\ucH(a)$ and $\ucF(a)$ are presentations of the same
codimension that are equivalent (i.e., $[\ucH(a)]=[\ucF(a)]$),
then, for all $H\in \cE(a)$,
$$
\mu_{\ucH(a),H} \ = \ \mu_{\ucF(a),H} .
$$
\end{theorem}

See \cite[Propositions 4.8, 4.11]{BMinv} or \cite[Propositions
4.4, 4.6]{BMmsri}
for proofs of these assertions. (The second assertion of
Theorem 2.3 is not stated explicitly in these references,
but is clear from the proof of \cite[Proposition 4.8]{BMinv}
or \cite[Proposition 4.4]{BMmsri}, and is worth noting 
-- see \S6.1.1 below.)

\begin{remarks}
As we have remarked in Section 1, our presentations $\ucH(a)
= ( N(a),\cH(a),\cE(a))$ are similar to Villamayor's 
basic objects. In the latter $\cH(a)$ is replaced by an 
{\it idealistic exponent} $(J,b)$ in the sense of Hironaka
\cite{Hid}: $J$ is an ideal in $\cO_{N(a)}$ and $b\in\IN$.
For example, choose $q\in\IN$ such that $q\cdot\mu_h\in\IN$,
for all $(h,\mu_h)\in\cH(a)$. Then we can take $b=\max(q\mu_h)!$
and $J=$ ideal generated by the $h^{b/\mu_h}$, for all
$(h,\mu_h)\in\cH(a)$. (Each $b/\mu_h\in\IN$.)

Our notion of equivalence of presentations, however, is stronger
than the notions of equivalence of idealistic exponents or of
basic objects used by Hironaka \cite{Hid} and Villamayor
\cite{V2, EVacta}.
The latter involve stability under transformations of types
(i), (ii) alone, so the corresponding equivalence classes 
(essentially $[\cdot]_{{\rm (i,ii)}}$) are larger.

Example 5.14 below shows that the conclusion of Theorem 2.4
is not necessarily true if $\ucH(a)$ and $\ucF(a)$ are merely
equivalent with respect to transformations of types (i) and (ii);
in particular, $[\cdot]$ is in general a {\it strictly smaller}
class of equivalence than $[\cdot]_{{\rm (i,ii)}}$. Example 5.14
shows that even the variant of $\nu(a)$ used by Villamayor is
not an invariant of the equivalence class of an idealistic exponent.
It is for this reason that the definitions of $\inv$,
the resolution algorithms and the meanings of
invariance are not the same in \cite{BMinv} and \cite{V2,EVacta}.
The proofs in the latter show that the underlying invariant
$\inv(a)$, $a\in M_j$ depends on the previous history
of the resolution process, but does not show that it depends
only on $X_{j,a}$, $E(a)$, and the $E^q(a)$.

On the other hand, we need to use $[\cdot]$ rather than the
smaller equivalence class $[\cdot]_{{\rm (i,ii,iii)}}$ 
-- see Remark 2.7 below.
\end{remarks}

\subsection{Semicoherent equivalence}
Let $\ucH(a)=( N(a),\cH(a),\cE(a))$ be a 
presentation.
Say $\cH(a)=\{ (h,\mu_h)\}$.
We identify $N(a)$ (respectively, each $h$) with a submanifold
(respectively, a function) in some neighbourhood of $a$.
Let $\ucH(x)=( N(x),\cH(x),\cE(x))$ denote the
presentation induced by $\ucH(a)$ at each $x\in N(a)$.
($\cE(x) := \{ H\in\cE(a):\ x\in H\}$.)

Let $\ucF(a)$ be another presentation at $a$.
We say that $\ucF(a)$ and $\ucH(a)$ are {\it semicoherent
equivalent} if $\ucF(x)$ and $\ucH(x)$ are equivalent
(in the sense of Definition 2.2) at each $x$ in a neighbourhood
of $a$ in $S_{\ucF(a)} = S_{\ucH(a)}$. 
This notion of semicoherent equivalence clearly depends
only on $\ucF(a)$ and $\ucH(a)$.
Write $[\![\ucH(a)]\!]$ for the semicoherent equivalence class.

The stronger notion of semicoherent equivalence is important
in the local construction \S2.6 below that will be 
used in Section 3 in the inductive definition
of $\inv$. The following theorem is the basis of the induction
on the codimension of a presentation.

\begin{theorem}
Let $\ucG(a)=( N(a),\cG(a),\emptyset)$ be a presentation
of codimension $p$. If $\mu_{\ucG(a)} = 1$, then $\ucG(a)$ is
semicoherent equivalent to a presentation $\ucC(a)=( N_{+1}(a),
\cC(a),\emptyset)$ of codimension $p+1$.
\end{theorem}

We will prove Theorem 2.6 in Section 5 below; see also 
\cite[Proposition 4.12]{BMinv}.

\subsection{Local construction} 
We use the notation of \S2.4. 
Suppose that $\mu(a)<\infty$.
Define
$$
D(a) \ = \ D_{\ucH(a)} \ := \ \prod_{H \in \cE(a)} x_H^{\mu_H(a)} ,
$$
where $x_H\in\cO_{N(a)}$ denotes a generator of the ideal of
$N(a)\cap H$, for each $H\in\cE(a)$.
If $\nu(a)=0$, we define
$$
\cG(a) \ = \ \cG_{\ucH(a)} \ := \ \left\{ ( D(a),1)\right\} .
$$
If $0<\nu(a)<\infty$, then, for each $(h,\mu_h)\in\cH(a)$,
we can write
\begin{equation}
h \ = \ D(a)^{\mu_h} \cdot g_h
\end{equation}
(see Remark 2.7 following), and we define
\begin{equation}
\begin{array}{rcl}
\cG(a) & = & \cG_{\ucH(a)} \\ 
& := & \left\{ ( g_h, \mu_h\nu(a)):
\ (h,\mu_h) \in \cH(a)\right\}
\cup \left\{ ( D(a),1-\nu(a))\right\} .
\end{array}
\end{equation}
(The element $( D(a),1-\nu(a))$ plays no part
and can be deleted unless $\nu(a)<1$.)
Set $\ucG(a) = \ucG_{\ucH(a)} = ( N(a),\cG(a),\cE(a))$.
Then $\mu_{\ucG(a)}=1$ and
\begin{equation}
S_{\ucG(a)} \ = \ \{ x \in S_{\ucH(a)} :\ \nu(x) \ge
\nu(a) \} ,
\end{equation}
where $\nu(x)=\nu_{\ucH(x)} = \min_{\ucH(a)} \mu_x (h) / \mu_h$.
((2.4) makes sense as an equality of germs at $a$.)

\begin{remark}          
Using (2.4), it is easy to see that $[\![\ucG(a)]\!]$ depends
only on $[\![\ucH(a)]\!]$ \cite[Proposition 4.24]{BMinv}.
We need to use $[\cdot]$ rather than the smaller equivalence
class $[\cdot]_{{\rm (i,ii,iii)}}$ because it is not {\it a
priori} true that the semicoherent class 
$[\![\ucG(a)]\!]_{{\rm (i,ii,iii)}}$ corresponding to 
$[\cdot]_{{\rm (i,ii,iii)}}$ depends only on 
$[\![\ucH(a)]\!]_{{\rm (i,ii,iii)}}$. We do not know whether
$[\![\cdot]\!] = [\cdot]$.
\end{remark}

\begin{remark}
The factors appearing in (2.2) are (perhaps rational) powers of
elements of $\cO_{N(a)}$.
We can avoid non-integral powers by replacing $\ucH(a)$ by
an equivalent presentation where all $\mu_h=d$, for
some $d\in\IN$: Choose $q\in\IN$ such that $q\mu_h\in\IN$,
for all $h$; let $d=\max(q\mu_h)!$ and replace each
$(h,\mu_h)$ by $(h^{d/\mu_h} , d)$ to get an equivalent
presentation as claimed.
If $\cH(a)=\{ (h,d)\}$ with common $d\in\IN$, then $D(a)^d$
is a monomial in $x_H$, $H\in\cE(a)$ (with integral powers)
and (2.2) becomes $h=D(a)^d g_h$, so that each $g_h \in \cO_{N(a)}$; 
$D(a)^d$ is the greatest common factor of the $h$
which is monomial in $x_H$, $H\in\cE(a)$.
\end{remark}

\section{The invariant and the desingularization principle}

In this section, we give the local inductive construction
needed to define the desingularization invariant $\inv = 
\inv_X$ or $\inv_\cJ$ (\S3.2), and to prove the desingularization
principle, Theorem 1.14 (see \S3.3). We compare our local construction
with that used in Villamayor's algorithm \cite{V2} (\S3.4), as well
as with the variant of Encinas and Villamayor \cite{EVacta} (\S3.6).
We answer the question, ``$\inv$ is an invariant of what?''; see
\S3.5.

\subsection{Presentation of a local invariant} 
Let $\io = \io_X$ (or $\io_\cJ$) denote a local invariant of
spaces $X$ (or ideals of finite type $\cJ$). (We assume that
$X$ is a closed subspace of a manifold $M$, or that $\cJ$ is
a subsheaf of $\cO_M$.) 

Let us first be more precise about the Definitions 1.12
used in the Hypotheses 1.10 that we will impose on $\io$.
Assume that $\io$ satisfies the semicontinuity hypotheses
1.10(1). We will use the following transforms of $X$ (or
$\cJ$) by the three types of morphisms listed in \S2.2: If
$\s: \ M' \to M$ is an $\io$-admissible (local) blowing-up
with smooth centre $C$, we consider the strict transform
$X'$ of $X$ (or the weak transform $\cJ'$ of $\cJ$). On the
other hand, if $\s: \ M' \to M$ is a morphism of either type
(ii) (product with a line) or type (iii) (exceptional
blowing-up, in the presence of exceptional divisors), we transform
$X$ (or $\cJ$) simply by inverse image: $X' := \s^{-1}(X)$
(or $\cJ' := \s^{-1}(\cJ)$).

Let $\ucG(a)=( N(a),\cG(a),\cE(a))$ denote a presentation
at $a$, as in \S2.1. Of course, if $S_{\ucG(a)} = S_\io(a)$, then
a local blowing-up at $a$ is $\io$-admissible if and only if
it is admissible for $\ucG(a)$. (See Definitions 1.9 and \S2.2.)

\begin{definitions}
$\ucG(a)$ is a {\it presentation} of $\io$ at $a$ if conditions
(1)-(3) of Definitions 1.12 hold, where (3) refers to any finite
sequence of morphisms allowed by Definition 2.2.

A presentation $\ucG(a)$ of $\io$ at $a$ is {\it semicoherent}
if it induces a presentation of $\io$ at each $x$ in a 
neighbourhood of $a$ in $S_\io(a) = S_{\ucG(a)}$ (cf. \S2.5).
\end{definitions}

The Hypothesis 1.10(3) means that $\io$ admits a semicoherent
presentation $\ucG(a)=( N(a),\cG(a),\emptyset)$
at each point $a$.

\begin{remark}
In general (e.g., in an analytic category), it is necessary to
be somewhat more precise about the meaning of a presentation
associated to an invariant (as above or as in
\S3.2 below): We assume that $M$ can be covered by coordinate
charts $U$ such that, for each $a \in U$, the functions involved
in the presentation at $a$ (e.g., the functions in $\cG(a)$ 
above) are quotients of elements of $\cO(U)$ with
denominators not vanishing at $a$ (likewise for a collection
of functions defining the maximal contact submanifold,
e.g., $N(a)$ above); cf. \cite[Definition and
remarks 4.14]{BMinv}. (The resulting definition is identical
to the preceding in the case of algebraic varieties or
schemes.) The equimultiple locus of a presentation at $a \in U$,
and the ideas of semicoherent presentation of an invariant
or of semicoherent equivalence of presentations
(\S2.5) involve germs with respect to the Zariski topology
of $U$. Presentations in this more precise sense exist in Examples
1.13 (according to the references given), and are needed to prove
upper semicontinuity of $\inv$ in Theorem 1.14 (\S3.3 below). 
\end{remark}

\begin{remarks}
(1) Suppose that $\ucG(a)=( N(a),\cG(a),\cE(a))$ is
a (semicoherent) presentation of $\io$ at $a$. If 
$\ucH(a)=( P(a),\cH(a),\cE(a))$ is a presentation
at $a$, then $\ucH(a)$ is a (semicoherent) presentation of
$\io$ at $a$ if and only if $\ucH(a)$ is (semicoherent)
equivalent to $\ucG(a)$. 
\smallskip

(2) Consider any sequence of $\io$-admissible transformations
(1.1) (or (1.2)). Let $a \in M_j$ and let $a_i$ denote the
image of $a$ in $M_i$, for all $i \leq j$. Suppose that 
$\io(a) = \io(a_0)$. If $\ucG(a_0)=( N(a_0),\cG(a_0),
\emptyset)$ is a (semicoherent) presentation of $\io$
at $a_0$, then we can consider the successive transforms
$\ucG(a_i)=( N(a_i),\cG(a_i),\cE(a_i))$ of $\ucG(a_0)$
at $a_i$, $i = 1,\ldots,j$. It follows from Definitions 3.1
that $\ucG(a)$ is a (semicoherent) presentation of $\io$
at $a$.
\smallskip

(3) Suppose that $\io$ satisfies the Hypotheses 1.10(1) and (3)
for any $X$ (or for any $\cJ$). Consider any sequence of
$\io$-admissible transformations (1.1) (or (1.2)). It follows
that, for all $j$ and all $a \in M_j$, there is a semicoherent
presentation $\ucG(a)=( N(a),\cG(a),\cE_1(a))$ of 
$\io$ at $a$, where $\cE_1(a) := E(a)\backslash E^1(a)$.
(Recall Definitions 1.15. We obtain $\ucG(a)$ simply by
transforming a presentation $\ucG(a_i)=( N(a_i),\cG(a_i),
\emptyset)$ of $\io$ at $a_i$, where $i \leq j$ is the
year of birth of $\io(a)$.)

By (1) above, if $\ucG(a)=( N(a),\cG(a),\cE_1(a))$ is a presentation
(respectively, semicoherent presentation) of $\io$ at $a$, then
the equivalence class $[\ucG(a)]$ (respectively, $[\![\ucG(a)]\!]$)
depends only on
$X_{j,a}$ (or $\cJ_{j,a}$) and $\cE_1(a)$.
\end{remarks}

The following is a corollary of Theorem 2.6.

\begin{corollary}
Suppose that $\io$ satisfies the Hypotheses 1.10(1) and (3).
Let $\ucG(a_0) = ( N(a_0),\cG(a_0),\emptyset)$ be a (semicoherent) 
presentation of $\io$ at $a_0 \in M$, of codimension $p(a_0)$, say.
Suppose that $\mu_{\ucG(a_0)} = 1$. Consider any
sequence of $\io$-admissible transformations (1.1) (or (1.2)).
Then, for all $j$ and all $a \in M_j$ lying over $a_0$, 
if $\io(a) = \io(a_0)$, $\io$ admits a
(semicoherent) presentation $\ucC(a)=( N_{+1}(a),\cC(a),
\cE_1(a))$ of codimension $p(a_0)+1$ at $a$.
\end{corollary}

\subsection{The local inductive construction and the
desingularization invariant}
Consider $\io = \io_X$ (or $\io_\cJ$) satisfying the Hypotheses
1.10. Our aim is to extend $\io$ to an invariant $\inv = \inv_X$
(or $\inv_\cJ$) defined recursively over a sequence of admissible
transformations (1.1) (or (1.2)). (See the introductory paragraph
of Section 2 above.)

For simplicity in this section, we will assume that $\io$ admits
a semicoherent presentation $\ucG(a)= ( N(a),\cG(a),
\emptyset)$ of codimension $0$ at every point $a$, where
$\mu_{\ucG(a)} = 1$. (For example, if $\io(a) = \io_\cJ(a)$ denotes
the order $\mu_a(\cJ)$ of an ideal of finite type 
$\cJ \subset \cO_M$ at $a \in
M$, then we can take $N(a) \ =$ the germ of $M$ at $a$, and
$\cG(a) = \{(g, \mu_a(\cJ)\}$, where the $g$ form any finite set of
generators of $\cJ_a$ -- cf. Examples 1.13, 1.19.) But
this simplifying assumption is not necessary; see \S6.1 below 
(in particular, for the Hilbert-Samuel function).

First consider a sequence of $\io$-admissible transformations
(1.1) (or (1.2)). Recall that, if $a \in M_j$, then $\io(a)$
denotes $\io_{X_j}(a)$ (or $\io_{\cJ_j}(a)$), and $E(a)$ denotes
$\{H \in E_j: \ a \in H\}$. We write $\inv_{1/2} := \io$.

Let $a \in M_j$. Define $E^1(a)$ as in Definitions 1.15. Set
$\cE_1(a) := \cE(a)\backslash E^1(a)$. By Corollary 3.4, 
$\inv_{1/2} = \io$ admits a semicoherent presentation
$$
\ucC_1(a) \ = \ ( N_1(a), \cC_1(a), \cE_1(a))
$$
of codimension $1$ at $a$. By Remarks 3.3, the equivalence
classes $[\ucC_1(a)]$
and $[\![\ucC_1(a)]\!]$
depend only on $X_{j,a}$ (or $\cJ_{j,a}$) and $\cE_1(a)$.
Define $\inv_1(a) := (\io(a), s_1(a) )$,
where $s_1(a) := \# E^1(a)$, and set
$$
\ucH_1(a) \ = \ ( N_1(a), \cH_1(a), \cE_1(a)) ,
$$
where $\cH_1(a) := \cC_1(a) \cup \left(E^1(a)\big|_{N_1(a)} , 1 \right)$ and
$$
\left(E^1(a)\big|_{N_1(a)} , 1\right) \ := \ \left\{(x_H\big|_{N_1(a)} , 1): \ 
H \in E^1(a)\right\}
$$
(and $x_H$ denotes a generator of $\cI_{H,a}$).

Clearly, $\ucH_1(a)$ is a {\it semicoherent presentation} of
$\inv_1$ at $a$, and
the equivalence class $[\ucH_1(a)]$ (respectively,
$[\![\ucH_1(a)]\!]$) depends only on $[\ucC_1(a)]$ and 
$E^1(a)$ (respectively, on $[\![\ucC_1(a)]\!]$ and $E^1(a)$).

\begin{remark}
By ``semicoherent presentation of $\inv_1$ at $a$'',
we mean the analogue of Definitions 3.1 for $\inv_1$,
but where the condition (3) (of Definitions 1.12) is replaced
by the weaker condition of stability after finite sequences of
transformations of type (i) (admissible blowings-up) only.

We use this weaker version of semicoherent presentation of an
invariant here (and below) because, together with the fact that
$[\![\ucH_1(a)]\!]$ depends only on $[\![\ucC_1(a)]\!]$ and $E^1(a)$,
it suffices to continue the induction, and avoids the necessity
of defining $\inv$ over sequences of all three types of tranformations
in \S2.2. The stronger version Definition 3.1 for $\inv_{1/2} = 
\io$ is needed only to start the induction.
\end{remark}

Define $\nu_2(a) := \nu_{\ucH_1 (a)}$ and $\inv_{3/2}(a)
:= ( \inv_1(a); \nu_2(a))$. 
Then $\nu_2 (a)$ depends only on $[\ucH_1(a)]$
(hence only on $[\![ \ucH_1(a)]\!]$);
thus $\inv_{3/2}(a)$
depends only on $X_{j,a}$, $E(a)$ and $E^1(a)$.

If $\nu_2(a)=\infty$, we write $\inv(a) := (
\nu_1(a), s_1(a);\infty)$; then $S_{\inv}(a) =
N_1(a)$.
($\nu_2(a)=\infty$ only if $s_1(a)=0$.)
If $\nu_2 (a)=0$, we write $\inv(a) := (\nu_1(a),
s_1(a); 0)$ and define $\ucG_2 (a) :=
\ucG_{\ucH_1 (a)} = ( N_1(a), \cG_2(a), \cE_1(a))$,
where $\cG_2(a)=\{ ( D_2(a),1)\}$,
and $D_2(a) := D_{\ucH_1 (a)}$ (see \S2.6).
Then $\ucG_2 (a)$ is a codimension $1$ presentation
of $\inv$ at $a$, and
$$
S_{\inv}(a) \ = \ \left\{ x\in N_1(a): \ \mu_a( D_2(a)
) \ge 1\right\} .
$$

On the other hand, suppose that $0<\nu_2(a)<\infty$. Then
$$
\ucG_{\ucH_1 (a)} \ = \ \left( N_1(a), \cG_2(a), \cE_1(a)\right),
$$
where $\cG_2(a) :=  \cG_{\ucH_1 (a)}$ is given by (2.3).
Then $\ucG_{\ucH_1 (a)}$ is a codimension $1$ semicoherent
presentation of $\inv_{3/2}$ at $a$, the semicoherent
equivalence class $[\![ \ucG_{\ucH_1 (a)} ]\!]$ depends only on
$[\![ \ucH_1 (a) ]\!]$ (Remark 2.7),
and $\mu_{\ucG_{\ucH_1 (a)}} = 1$.
Define $\cE_2(a) := \cE_1(a) \backslash E^2(a)$ (see
Definitions 1.15) and put
$$
\ucG_2(a) \ := \ ( N_1(a), \cG_2(a), \cE_2(a)).
$$
Then $\ucG_2(a)$ is a semicoherent presentation of 
$\inv_{3/2}$ at $a$, and $[\ucG_2(a)]$ (respectively,
$[\![ \ucG_2(a) ]\!]$) depends only on $[\ucG_{\ucH_1 (a)}]$
and $E^2(a)$ (respectively, on $[\![ \ucG_{\ucH_1 (a)} ]\!]$      
and $E^2(a)$). It follows from Theorem 2.6 (cf. Corollary 3.4)
that $\ucG_2(a)$ is semicoherent equivalent to a presentation
$$
\ucC_2(a) \ = \ ( N_2(a), \cC_2(a), \cE_2(a))
$$ of codimension $2$ (where $N_2(a)$ is a submanifold of
$N_1(a)$). This completes one cycle in the inductive definition
of $\inv$.

In general, let $r\geq 1$ and suppose that we have
introduced $\inv_{r-1/2} (a) = ( \inv_{r-1}(a);
\nu_r (a))$. Consider a sequence of $(r-1/2)$-admissible
transformations (1.1) (or (1.2)). Let $a \in M_j$. Assume
that $0<\nu_r(a)<\infty$. By induction, we can assume that
$\inv_{r-1/2}$ admits a semicoherent presentation 
$$
\ucC_r(a) \ = \ ( N_r(a), \cC_r(a), \cE_r(a))
$$
of codimension $r$ at $a$, where $\cE_r(a) := \cE_{r-1}(a)
\backslash E^r(a)$, and that the semicoherent equivalence
class of $\ucC_r(a)$ depends only on $X_{j,a}$ (or $\cJ_{j,a}$),
$E(a), E^1(a), \ldots, E^r(a)$. (We can assume, inductively,
that the semicoherence class of $\ucC_r(a)$ depends only on
$\cE_r(a)$ and the semicoherence class of $\ucH_{r-1}(a)$.)

Define $\inv_r(a) := ( \inv_{r-1/2}(a), s_r(a))$,
where $s_r(a) := \# E^r(a)$, and set
$$
\ucH_r(a) \ := \ ( N_r(a), \cH_r(a), \cE_r(a)),
$$
where $\cH_r(a) := \cC_r(a) \cup \left(E^r(a)\big|_{N_r(a)}, 1 \right)$.
Then $\ucH_r(a)$ is a codimension $r$ semicoherent presentation
of $\inv_r$ at $a$, and its equivalence class (or semicoherent
equivalence class) depends only on $E^r(a)$ and that of $\ucC_r(a)$.

Define $\mu_{r+1} (a) := \mu_{\ucH_r (a)}$,
$\mu_{r+1,H}(a) := \mu_{\ucH_r (a),H}$, for all
$H\in\cE_r (a)$, and
\begin{equation}
\nu_{r+1} (a) \ := \ \nu_{\ucH_r (a)} \ = \ \mu_{r+1}
(a) - \sum_{H\in\cE_r (a)} \mu_{r+1,H} (a).
\end{equation}
If $0\leq \nu_{r+1} (a) < \infty$, define
$$
D_{r+1} (a) \ := \ D_{\ucH_r (a)} \ = \ \prod_{H \in \cE_r (a)}
x_H^{\mu_{r+1,H} (a)} ,
$$
and introduce 
$\ucG_{\ucH_r (a)}$ as in \S2.6.
Then $\ucG_{\ucH_r (a)} = ( N_r(a), \cG_{r+1}(a), \cE_r(a))$ 
is a semicoherent
codimension $r$ presentation of $\inv_{r+1/2} :=
(\inv_r;\nu_{r+1})$ at $a$, and
$[\![ \ucG_{\ucH_r (a)}]\!]$ depends
only on $[\![\ucH_r (a)]\!]$. If $0 < \nu_{r+1} (a) < \infty$,
then 
$$
\ucG_{r+1}(a) \ := \ ( N_r(a), \cG_{r+1}(a), \cE_{r+1}(a)),
$$
where $\cE_{r+1}(a) := \cE_r(a) \backslash E^{r+1}(a)$, is
semicoherent equivalent to a presentation
$$
\ucC_{r+1}(a) \ := \ ( N_{r+1}(a), \cC_{r+1}(a), \cE_{r+1}(a))
$$
of codimension $r+1$. Clearly $[\![\ucC_{r+1}(a)]\!]$ depends
only on $[\![\ucH_r (a)]\!]$ and $\cE_{r+1}(a)$. Etc.

Eventually, we find $t\leq n := \dim_a M_j$ such that $\nu_{t+1}(a)=0$
or $\infty$, and we set $\inv(a) := \inv_{t+1/2} (a)$.
Suppose that $\nu_{t+1}(a)=\infty$.
Then $S_{\inv}(a) = N_t(a)$.
On the other hand, if $\nu_{t+1} (a) =0$, then
\begin{equation}
S_{\inv} (a) \ = \ \left\{ x\in N_t(a) :\ \mu_x
( D_{t+1} (a)) \ge 1\right\} ,
\end{equation}
where $D_{t+1} (a) = \prod_{\cE_t(a)} x_H^{\mu_{t+1,H}(a)}$. 

\begin{remark}
In the local inductive construction above, we pass from $\ucG_r(a)
= ( N_{r-1}(a), \cG_r(a), \cE_r(a))$
to an equivalent presentation $\ucC_r(a) = ( N_r(a), 
\cC_r(a), \cE_r(a))$
in codimension $+1$, and
then to $\ucH_r(a)$ by adjoining $\left(E^r(a)\big|_{N_r(a)} , 1 \right)$.
The construction of \cite{BMinv} follows a slightly different
route -- from $\ucG_r(a)$ to $\ucF_r(a) = ( N_{r-1}(a), \cF_r(a),
\cE_r(a))$, where $\cF_r(a) = \cG_r(a) \cup 
\left(E^r(a)\big|_{N_{r-1}(a)}, 1\right)$,
and then to an equivalent presentation $\ucH_r(a)$ in codimension
$+1$. The latter gives a little more flexibility in the
choice of the maximal contact submanifold $N_r(a)$, but the
construction above makes for a clearer parallel treatment of
alternative initial invariants $\io$ -- see, for example, \S6.1.
\end{remark}

\subsection{Proof of the desingularization principle Theorem 1.14}

(1) Semicontinuity. The semicontinuity properties (a) and (b) can
be proved for the truncated invariants $\inv_{r-1/2}$ and $\inv_r$
by induction on $r$. By the hypothesis 1.10(1), $\inv_{1/2}$ satisfies
(a) and (b). Assume that $\inv_{r-1/2}$ satisfies (a) and (b), where
$r \geq 1$. The properties (a) and (b) for $\inv_r$ are then
consequences of the following semicontinuity assertion for 
$E^r(\cdot)$:\ If $a \in M_j$, then $E^r(x) = E(x) \cap E^r(a)$,
where $x \in S_{\inv_{r-1/2}}(a)$ \cite[Proposition 6.6]{BMinv}. 
To prove (a) and (b) for 
$\inv_{r+1/2}$, consider a semicoherent presentation 
$\ucH_r(a) \ := \ ( N_r(a), \cH_r(a), \cE_r(a))$ of 
$\inv_r$ at $a \in M_j$; say $\cH_r(a) = \{(h, \mu_h)\}$. If
$x \in S_{\inv_r}(a)$, then 
$$
\mu_h \nu_{r+1}(x) \ = \ \min_{\cH_r(a)}\mu_x\left(\displaystyle{
\frac{h}{D_{r+1}(a)^{\mu_h}}}\right) .
$$
(See \S\S2.4, 2.6.) The semicontinuity properties for $\inv_{r+1/2}$
are consequences of the analogous properties for the order of an
element $g = h/D_{r+1}(a)^{\mu_h}$ such that $\mu_a(g) = 
\mu_h \nu_{r+1}(a)$. (This is where Remark 3.2 is relevant, in
general.)
\smallskip

(2) Stabilization. 
Suppose that $\ucC_1(a)$ is a presentation of 
$\inv_{1/2} = \io$ at $a \in M_j$, as above, and that
$\cC_1(a) = \{(c, \mu_c)\}$. Choose $q \in \IN$
such that $q \cdot \mu_c \in \IN$, for all
$(c, \mu_c)$. Let $e_1 = e_1(a) := \max q \cdot \mu_c\, $.
Then, for all $r > 0$, $e_r! \nu_{r+1}(a) \in \IN$,
where $e_{r+1} = \max\{ e_r!, e_r! \nu_{r+1}(a) \}$.
The assertion follows from Hypotheses 1.10
for $\io$ (using infinitesimal semicontinuity
((1) above) and the stabilization property 1.12(3) of 
a presentation of $\io$ at $a$).
\smallskip

(3) Normal crossings. Say $\inv(a) = \inv_{t+1/2}(a)$. The
assertion is an immediate consequence of the fact that 
$S_{\inv}(a) = N_t(a)$ in the case that $\nu_{t+1}(a) = 
\infty$, and of (3.2) in the case that $\nu_{t+1}(a) = 0$.
(See also (4) below).
\smallskip

(4) Decrease. Say $\inv(a) = \inv_{t+1/2}(a)$. First suppose
that $\nu_{t+1}(a) = \infty$. Then $S_{\inv}(a) = N_t(a)$.
If $\s$ is the local blowing-up with centre $N_t(a)$, then
the strict transform $N_t(a)' = \emptyset$, so that $\inv(a')
< \inv(a)$, for all $a' \in \s^{-1}(a)$. (See Definitions 1.12
and 3.1.)

On the other hand, suppose that $\nu_{t+1}(a) = 0$. Then
$h = D_{t+1}(a)^{\mu_h}$, for some $(h, \mu_h) \in \cH_t(a)$,
and $S_\inv(a) = S_{\inv_t}(a) = \left\{x \in N_t(a):\ 
\mu_x(D_{t+1}(a)) \geq 1\right\}$. (We can choose coordinates
$x = (x_1,\ldots, x_{n-t})$ for $N_t(a)$ such that $D_{t+1}(a)$
is a monomial $x_1^{\Om_1}\cdots x_{n-t}^{\Om_{n-t}}$ with
rational exponents and, if $\Om_l \neq 0$, then $x_l = x_H$,
for some $H \in \cE_t(a)$, and $\Om_l = \mu_{t+1,H}(a)$. Thus
$\mu_x(D_{t+1}(a))$ makes sense as a rational number.) Therefore,
$S_{\inv_t}(a)$ is a union of smooth components $\bigcup_I Z_I$,
where $Z_I = \{x \in N_t(a):\ x_l = 0, \, l \in I\}$ and the
union is over the minimal subsets $I$ of $\{1,\ldots, n-t\}$
such that $\sum_{l \in I} \Om_l \geq 1$; equivalently, over
the subsets $I$ such that 
$$
0 \ \leq \ \sum_{k \in I} \Om_k - 1 \ \leq \ \Om_l, \quad
\mbox{for all } l \in I .
$$
Set $\mu (a) := \mu_{t+1}(a)$. Consider a local blowing-up
$\s$ with centre $Z_I$, for some $I$ as above. Suppose that
$a' \in \s^{-1}(a)$ and $\inv_t(a') = \inv_t(a)$. Then $a'
\in N(a)'$. $N(a)'$ is a union of coordinate charts 
$\bigcup_{l \in I} U_l'$ such that $\s \big|_{U_l'}$ is given by the
substitution $x_l = y_l$, $x_k = y_l y_k$ if $k \in I\backslash
\{l\}$, and $x_k = y_k$ if $k \notin I$. Consider 
$h = D_{t+1}(a)^{\mu_h}$.
Suppose that $a' \in U_l'$. Write $d = \mu_h$. 
Then $(h', d) \in \cH_t(a')$,
where $h' := y_l^{-d}\left( D_{t+1}(a)^d \circ \s\right)
= \left(y_1^{\Om_1'} \cdots
y_{n-t}^{\Om_{n-t}'}\right)^d$, and
$$
\Om_k' \ = \ \Om_k, \ \ k \neq l; 
\ \quad \Om_l' \ = \ \sum_{k \in I} \Om_k - 1 .
$$
Therefore, 
$$
1 \ \leq \ \mu_{t+1}(a') \ \leq \ \sum_{k = 1}^{n-t} \Om_k' \ 
< \ \sum_{k = 1}^{n-t} \Om_k \ = \ \mu_{t+1}(a) ,
$$
as required.

\subsection{Villamayor's algorithm}
Villamayor's algorithm can be described in the framework of
\S3.2 above, by making a simple modification in the construction
and the corresponding invariant.
This modification corresponds to using only a certain subset
$\cE'_r(a)\subset\cE_r(a)$ when we define $\nu_{r+1}(a)$
according to the formula (3.1); equivalently, to factoring
from $\cH_r(a)$ only a part of the exceptional monomial $D_{r+1}(a)$
in order to define $\cG_{r+1}(a)$.
Such a modification for given $r$ will, in general,
change all subsequent terms of the invariant and the corresponding
presentations.

Consider $r \geq 1$. Suppose (inductively) that we have introduced
$\inv_{r-1/2}(a) = ( \inv_{r-1}(a); \nu_r(a))$. As before,
we define $E^r(a) := \{ H\in\cE_{r-1}(a):\ H$ is transformed from
$E(a_i)\}$, where $i$ denotes the year of birth of $\inv_{r-1/2}(a)$
(see Definitions 1.15), and we set $s_r(a) := \#E^r(a)$, $\cE_r(a)
:= \cE_{r-1}(a) \backslash E^r(a)$. But now let $\cE'_r(a) \subset
\cE_r(a)$ denote the subset consisting of only those exceptional
hypersurfaces passing through $a$ that have accumulated since the
year of birth of $\inv_r(a) = (\inv_{r-1}(a); \nu_r(a), 
s_r(a))$. In other words, if $i'$ denotes the year of birth
of $\inv_r(a)$, then $\cE'_r(a) := \cE_{r-1}(a) \backslash E^{r'}(a)$,
where $E^{r'}(a) := \{H\in\cE_{r-1}(a):\ H$ is transformed from
$E(a'_i)\}$.

We define $\inv_1$ exactly as before, but we use $\ucH_1(a) \ = 
\ ( N_1(a), \cH_1(a), \cE'_1(a))$ as an associated
presentation (where $N_1(a)$, $\cH_1(a)$ have the same meaning
as before), and we define
$$
\nu_2(a) \ := \ \mu_2(a) - \sum_{H \in \cE'_2(a)} \mu_{2,H}(a),
$$
where $\mu_2(a) := \mu_{\ucH_1(a)}$ and $\mu_{2,H}(a) :=
\mu_{\ucH_1(a),H}$ for all $H \in \cE'_2(a)$. (See \S2.4.)
Then $\inv_{3/2} = (\inv_1; \nu_2)$ has a codimension
$1$ presentation 
$$
\ucG_2(a) \ = \ ( N_1(a), \cG_2(a), \cE'_2(a)),
$$
where $\cG_2(a) = \cG_{\ucH_1(a)}$. We pass to an equivalent
presentation 
$$
\ucC_2(a) \ = \ ( N_2(a), \cC_2(a), \cE'_2(a))
$$
in codimension $2$, and adjoin $\left(E^2(a)\big|_{N_2(a)}, 1\right)$
to get a presentation of $\inv_2 = (\inv_{3/2}, s_2)$.

(Although we are using the same notation for presentations  
as before),  $\nu_2$ and all subsequent terms of the invariant
will, in general, be
different, as will the corresponding presentations:
The change in the local construction is repeated for each
successive $r$; when we pass from
from $\ucH_r(a)$ to 
$\ucG_{\ucH_r(a)} = ( N_r(a),\cG_{r+1}(a), \cE'_r(a))$
in order to obtain a presentation of $\inv_{r+1/2}$ at $a$, 
we factor from $\cH_r(a)$ not
$D_{r+1}(a) = \prod_{\cE_r(a)} x_H^{\mu_{r+1,H}(a)}$ as before,
but only the product
$$
D'_{r+1} (a) \ := \ \prod_{H\in\cE'_r (a)} \mu_H^{\mu_{r+1,H}(a)}
$$
(i.e., the product over those $H$ accumulated since the birth
of $\inv_r (a)$).
This change will be magnified in all subsequent
terms of $\inv$ because $s_{r+1}$ depends on $\inv_{r+1/2}$,
etc.

\subsection{Meaning of invariance} 
Consider $\inv = \inv_X$ (or $\inv_\cJ$) defined over a 
sequence of admissible blowings-up (1.1) (or (1.2)) according
to the construction of \S3.2 above. We have shown in \S3.2
that, if $a \in M_j$, then $\inv(a)$ depends only on the 
local isomorphism class of $\left(X_{j,a}(a)\ (\mbox{or } \cJ_{j,a}), 
E(a), E^1(a), E^2(a),\ldots\right)$. The key point is that, for each $r$,
$\nu_{r+1}(a)$ depends only on the equivalence class $[\ucH_r(a)]$
(by Theorems 2.3 and 2.4).

According to these theorems, the version of $\nu_{r+1}(a)$ 
defined in \S3.4 is also an invariant of the corresponding 
equivalence class $[\ucH_r(a)]$. 
{\it But it is not an invariant with respect to the weaker
notion of equivalence of idealistic exponents} \cite{Hid} 
{\it or equivalence of basic objects} \cite{V2, EVacta}
-- see Example 5.14 below.

This is the reason for our introducing the idea
of equivalence of presentations involving exceptional blowings-up.
The result is not only a proof of invariance in a stronger
sense, but also an invariant whose maximum locus, in general,
provides a larger centre of blowing up because Theorem 2.4
shows that the $\mu_{r+1,H}(a)$ are invariants for the
larger block $\cE_r(a)$ of exceptional divisors $H$.

The difference in the algorithms is thus not accidental:
for the smaller block $\cE'_r(a)$
of exceptional divisors $H$ used by Villamayor, there is
another approach to invariance of $\mu_{r+1,H}(a)$
in terms of the previous history:

\begin{lemma}
Suppose that $a\in M_j$ and that $\inv_r(a)=\inv_r(a_{j-1})$.
Let $H\in \cE'_r(a)$.
(Note that $\cE'_r(a) =\emptyset$ unless $\inv_r(a)=\inv_r(a_{j-1})$.)
If $H=\s_j^{-1} (C_j)$, then
$$
\mu_{r+1,H}(a) \ = \ \sum_{K\in J} \mu_{r+1,K} (a_{j-1}) 
+ \nu_{r+1}(a_{j-1}) - 1 ,
$$
where $J := \{ K\in\cE'_r (a_{j-1}) :\ C_j\subset K\}$.
Otherwise, $H$ is the strict transform 
of an element $K\in\cE'_r(a_{j-1})$,
and
$$
\mu_{r+1,H} (a) \ = \ \mu_{r+1,K} (a_{j-1}) .
$$
\end{lemma}

The proof is a simple calculation.
(Cf. \cite[\S5.4]{V2}.)   
The transformation formulas of Lemma 3.7 are in general
not valid for the exceptional divisors $H \in \cE_r(a) \backslash
\cE'_r(a)$.

Desingularization as realized by either the algorithm of
the authors or that of Villamayor is canonical; in fact,
local isomorphisms between open subsets of $X$ lift throughout
the sequence of blowings-up determined by the algorithm.
But the notion of invariance provided by Lemma 3.7 is weaker
than that determined by equivalence of presentations because
it depends in a stronger way on the history.
Of course, the fact that $\nu_{r+1}(a)$ as defined in \S3.4
is also an invariant of the corresponding equivalence class
$[\ucH_r(a)]$ shows that Villamayor's version of $\inv_X$  
is an invariant in a stronger
sense than shown by Lemma 3.7; at $a\in X_j$, it depends on
$X_{j,a}$, $E(a)$, the $E^r(a)$ and the $E^{r'}(a)$.

\subsection{The variant of Encinas and Villamayor}
Encinas and Villamayor \cite{EVacta} have modified Villamayor's algorithm
to avoid blowings-up that are superfluous in a certain situation
based on Abkyankar's idea of ``good points'' \cite{Ab}. 
\cite{EVacta} suggests a substantial modification of
Villamayor's construction; in particular, each pair $(\nu_r,s_r)$
in the invariant as described above is replaced by a triple.
But the inductive construction described in \S3.2 above is 
flexible, and \cite{EVacta} 
can be understood in the following way.

We make a modification in the definition of $\inv_r$ by induction
on $r$:
Suppose that we have defined $\inv_{r-1}$ and an associated
codimension $r-1$ presentation $\ucH_{r-1}(a) =
( N_{r-1}(a), \cH_{r-1}(a), \cE_{r-1}(a))$ at $a\in M_j$,
as above.
(We are using the notation above.
The modification of \cite{EVacta} can be applied equally to 
the construction
of the authors (\S3.2 above) or that of Villamayor (\S3.4), so that
$\inv_{r-1}$ and $\ucH_{r-1}(a)$ here mean the notions defined
for either algorithm.
For Villamayor's version, therefore, $\cE_{r-1}(a)$ here is to be
understood as $\cE'_{r-1}(a)$ in the notation of \S3.4.)

Define $\mu_r(a)$, $\mu_{rH}(a)$ for all $H\in \cE_{r-1}(a)$,
$\nu_r(a)$, $D_r(a) = \prod x_H^{\mu_{rH}(a)}$ and $\cG_r(a)$
as before.
(In the case of Villamayor's construction again, $D_r(a)$ here means
$D'_r(a)$ in the notation of \S3.4.)
Of course,
$$
\mu_r(a) \ = \ \sum_{H\in\cE_r(a)} \mu_{rH} (a) + \nu_r(a) \ \ge \ 1 .
$$
But now, we make a modification of the definition of $\nu_r$
in the case that
$$
\sum \left( \mu_{rH} (a) - [\mu_{rH}(a)]\right) + \nu_r(a) \ < \ 1
$$
(where $[\cdot]$ denotes the integral part):
In this case, we redefine $\nu_r(a)$ as $\nu_r^*(a) := 0$,
and $\cG_r(a) := \{ (D_r(a),1)\}$.
This is the only change.

Like $\nu_r(a)$, the modified version $\nu_r^*(a)$ depends only
on the equivalence class $[\ucH_{r-1}(a)]$ since it is defined
in terms of $\mu_r(a)$ and the $\mu_{rH}(a)$.

It follows that the corresponding modified invariant $\inv^*$ 
can be defined
inductively over a sequence of transformations 
(1.1) (or (1.2)), provided we assume that each
successive centre of blowing up is a component of the maximum
locus of $\inv^*$.
Theorem 1.14 holds for such sequences; our proof as sketched above
applies with no change.
(When $\nu_r^*(a)=0$, the $\inv^*$-stratum of $a$ coincides
with the $\inv_{r-1}$-stratum and has only normal crossings
as in Theorem 1.14(3), where each component has codimension $1$
in $N_{r-1}(a)$ -- this special situation is analogous to the
idea of a ``good point''.
The blowing-up of $N_{r-1}(a)$ with centre such a component is
the identity.)
This is the variant of Encinas and Villamayor; the situation in
which it applies does not occur in Example 1.2 -- see Section 4
following.

\section{Worked example}

We illustrate the desingularization principle in this section
by working out Example 1.2 above (except for the estimates
on $n(AB)$ and $n(LB)$ given by the affine and locally binomial
algorithms. These will be computed in Desingularization 
algorithms II.) Let
$X$ denote the hypersurface in 4-dimensional affine space $M$ 
given by $g_0 = 0$, where
\begin{equation}
g_0(x,y,z,w) \ = \ z^d w^{d-1} - x^{d-1} y^d ,
\end{equation}
for any natural number $d\ge 2$.
The hypersurface $X$ takes its maximum order $2d-1$ precisely at
the origin, so in any case we take $C_0 = \{ 0\}$ as the
centre of the first blowing-up $\s_1$: $M_1\to M_0 =M$.
Then $M_1$ can be covered by four affine coordinate charts,
corresponding to the four variables.
For example, $\s_1$ is given in the ``$w$-chart'' $U_w$ by
substituting $(xw,yw,zw,w)$ for the original variables
$(x,y,z,w)$.
The strict transform $X_1=X'_0$ of $X_0=X$ is given in $U_w$
by $g_1=0$, where $g_1=w^{-(2d-1)} g_0\circ \s_1$; i.e.,
$$
g_1(x,y,z,w) \ = \ z^d - x^{d-1} y^d .
$$
(For economy of notation, we are using the same letters for
the variables before and after blowing up, and we are ``bookkeeping''
by writing $U_w$ for the ``$w$-chart'' of $M_1$; $U_w$ is
the complement in $M_1$ of the strict transform by $\s_1$
of the coordinate subspace $\{ w=0\}$ of $M_0=M$.)

We will estimate the number of blowings-up needed to reduce
the maximum order $d$ of $X_1$.
Let $n(BM)$ and $n(V)$ denote the number of blowings-up
determined by the algorithms of \cite{BMinv} and \cite{V2}, 
respectively.
We will show that
\begin{eqnarray*}
n(BM) &\le& 2d+j ,\\
n(V) &\ge& 9d+k ,
\end{eqnarray*}
where $j$, $k$ are independent of $d$.
It is easy to check in the calculations below that the variant
of Encinas-Villamayor \cite{EVacta} (see \S3.6) makes no difference
to either algorithm when applied to (4.1).

\subsection{Year one}
$X_1\cap U_w$ is given by $g_1=0$ as above.
(Note that the order at $0$ of a binomial majorizes its
orders at points of the chart.
In particular, $X_1$ has order $<d$ throughout the charts
$U_z$, $U_y$.)
Let $a_1=0$.
Then $E(a_1) = \{ H_1\}$, where $H_1$ denotes the exceptional
hypersurface $\{ w=0\}$.
We have $\nu_1(a_1)=d$, $E^1(a_1)=E(a_1)=\{ H_1\}$, $s_1(a_1)=1$,
and $\cE_1(a_1) := E(a_1)\backslash E^1(a_1)=\emptyset$.
Then $\inv_{1/2}=\nu_1$ has a (semicoherent) codimension $0$
presentation at $a_1$ given by $\cG_1(a_1)=\{ (g_1,d)\}$; therefore,
$\inv_{1/2}$ and $\inv_1 = (\nu_1,s_1)$ have (semicoherent)
codimension $1$ presentations
\begin{eqnarray*}
\ucC_1(a_1) &=& ( N_1(a_1), \cC_1(a_1), \cE_1(a_1)) ,\\
\ucH_1(a_1) &=& ( N_1(a_1), \cH_1(a_1), \cE_1(a_1)) ,
\end{eqnarray*}
(respectively), 
where $N_1(a_1) = \{z=0\}$, 
$\cC_1(a_1) = \{ (x^{d-1} y^d,d)\}$ and
$$
\cH_1(a_1) \ = \ \{ (x^{d-1} y^d,d),\ (w,1)\} .
$$
Therefore,
$$
\nu_2 (a_1) \ = \ \mu_2 (a_1) \ := \ \min_{(h,\mu_h)\in\cH_1(a_1)}
\frac{\mu_a(h)}{\mu_h} \ = \ 1
$$
and $s_2(a_1) = 0$. Then $\inv_{3/2}$ has presentations
\begin{eqnarray*}
\ucG_2(a_1) &=& ( N_1(a_1),\, \cG_2(a_1) = \cH_1(a_1),
\, \emptyset), \\
\ucC_2(a_1) &=& ( N_2(a_1), \cC_2(a_1), \emptyset),
\end{eqnarray*}
of codimensions $1$ and $2$ (respectively), where 
$N_2(a_1) = \{ z=w=0\}$ and $\cC_2(a_1) = \{ (x^{d-1} y^d,d)\}$,
and $\inv_2$ has a codimension $2$ presentation
$$
\ucH_2(a_1) \ = \ ( N_2(a_1),\, \cH_2(a_1) = \cC_2(a_1),
\, \emptyset) .
$$
(There is no change from $\cH_1(a_1)$ to $\cG_2(a_1)$ because
$\nu_2(a_1)=1$, and no change from $\cC_2(a_1)$ to
$\cH_2(a_1)$ because $s_2(a_1)=0$.)
Therefore,
$$
\nu_3(a_1) \ = \ \mu_3(a_1) \ = \ \frac{2d-1}{d} 
$$
and $s_3(a_1)=0$. Then $\inv_{5/2}$ admits 
a codimension $2$ presentation at $a_1$,
$$
\ucG_3(a_1) \ = \ ( N_2(a_1), \cG_3(a_1), \emptyset) ,
$$
where $\cG_3(a_1) = \{ (x^{d-1} y^d , 2d-1) \}$, and we can replace
the latter by
$$
\cG_3(a_1) \ = \ \{ (x,1),(y,1)\}
$$
to get an equivalent presentation (cf. Corollary 5.10 
below). Therefore,
$\inv_{5/2}$ has a codimension $3$ presentation at $a_1$,
$$
\ucC_3(a_1) \ = \ ( N_3(a_1), \cC_3(a_1), \emptyset) ,
$$
where $N_3(a_1)=\{z=w=y=0\}$ and $\cC_3(a_1)=\{(x,1)\}$,
and $\inv_3$ has the same codimension $3$ presentation
$\ucH_3(a_1) = \ucC_3(a_1)$. 
Hence $\nu_4(a_1)=1$, $s_4(a_1)=0$, and $\inv_4$ has a codimension 4
presentation $\ucH_4(a_1)=( N_4(a_1),\emptyset,\emptyset)$ at
$a_1$, where $N_4(a_1) = \{z=w=y=x=0\} = \{0\}$.
Thus
$$
\inv_X(a_1) \ = \ \left( d,1;\, 1,0;\, \frac{2d-1}{d}, 0;
\, 1,0;\, \infty \right) ;
$$
the preceding presentation of $\inv_4$ is also a presentation
of $\inv_X$, and, in the chart $U_w$, the centre of the
next blowing-up $\s_2$ is $C_1 = N_4(a_1)=\{0\}$.
By symmetry, as a global centre of the blowing-up $\s_2$:
$M_2\to M_1$ we would take $C_2 :=$ union of the origins
in the charts $U_w$, $U_x$.

\subsection{Year two} 
Let $X_2$ denote the strict transform $X'_1$ of $X_1$ by $\s_2$.
Consider the chart $U_{wy}$ of $M_2$, in which $\s_2$ is
given by the substitution $(xy,y,yz,yw)$.
Then $X_2\cap U_{wy}=\{ g_2=0\}$ where
$$
g_2(x,y,z,w) \ = \ z^d - x^{d-1} y^{d-1} .
$$
Let $a_2=0$.
Then $E(a_2)=\{ H_1,H_2\}$, where $H_1$ and $H_2$
denote the exceptional hypersurfaces $\{ w=0\}$ and
$\{ y=0\}$, respectively; $H_2=\s_2^{-1} (C_1)$ and $H_1$
is the strict transform of the hypersurface in year one
that we also denoted $H_1$ (to economize notation).
We have $\nu_1(a_2)=d$, $E^1(a_2)=\{ H_1\}$, $s_1(a_2)=1$,
and $\cE_1(a_2) := E(a_2) \backslash E^1(a_2) = \{ H_2\}$.
At $a_2$, $\inv_1$ admits a codimension $1$ presentation
$$
\ucH_1(a_2) \ = \ ( N_1(a_2), \cH_1(a_2), \cE_1(a_2)) ,
$$
where $N_1(a_2) = \{ z=0\}$ and
$$
\cH_1(a_2) \ = \ \{ (x^{d-1} y^{d-1},d),\ (w,1)\} .
$$
We compute
\begin{eqnarray*}
\mu_2(a_2) & := & \min_{\cH_1(a_2)}\frac{\mu_{a_2}(h)}{\mu_h} \ 
= \ 1 ,\\
\mu_{2H_2}(a_2) & := & \min_{\cH_1(a_2)}
\frac{\mu_{H_2,a_2}(h)}{\mu_h} \ = \ 0 ,\\
\nu_2(a_2) & := & \mu_2(a_2) -\sum_{H\in \cE_1(a_2)} \mu_{2H}(a_2) \ 
= \ 1 - 0 \ = \ 1 .
\end{eqnarray*}
Therefore, $E^2(a_2)=\emptyset$ and $s_2(a_2)=0$; $\inv_2$ has a 
codimension $2$ presentation
$$
\ucH_2(a_2) \ = \ ( N_2(a_2), \cH_2(a_2), \cE_2(a_2)) ,
$$
where $N_2(a_2) = \{ z=w=0\}$ and $\cH_2(a_2) 
= \{ (x^{d-1} y^{d-1},d)\}$.
Therefore,
\begin{eqnarray*}
\mu_3(a_2) & = & \frac{2d-2}{d}, \quad 
\mu_{3H_2}(a_2) \ = \ \frac{d-1}{d} ,\\
\nu_3 (a_2) & = & \mu_3(a_2) - \mu_{3H_2}(a_2) \ = \ \frac{d-1}{d} ,
\end{eqnarray*}
and $s_3(a_2)=1$.
At $a_2$, $\inv_{5/2}$ has a codimension 2 presentation
$$
\ucG_3(a_2) \ = \ ( N_2(a_2), \cG_3(a_2), \cE_3(a_2)) ,
$$
where $\cE_3(a_2)=\emptyset$ and 
$$
\cG_3(a_2)\ = \ \left\{ (x^{d-1},d-1),\ 
\left( y^{\frac{\scriptstyle d-1}{\scriptstyle d}} ,\, 
1 - \frac{d-1}{d}\right) \right\}
$$
or, equivalently,
$$
\cG_3(a_2) \ = \ \left\{ (x,1),\ \left( y,\, \frac{1}{d-1}\right)\right\} ;
$$
$\inv_{5/2}$ has a codimension $3$ presentation
$$
\ucC_3(a_2) \ = \ ( N_3(a_2), \cC_3(a_2), \cE_3(a_2)) ,
$$
where $N_3(a_2)=\{ z=w=x=0\}$ and $\cC_3(a_2) = ( y,\, 
1/(d-1) )$. Then $\inv_3$ has a codimension $3$ presentation
$$
\ucH_3(a_2) \ = \ ( N_3(a_2), \cH_3(a_2), \cE_3(a_2)) ,
$$
where $\cH_3(a_2)=\{ (y,1)\}$.
Clearly,
$$
\inv_X(a_2) \ = \ \left( d,1;\, 1,0;\, \frac{d-1}{d},1;
\, 1,0;\,\infty\right)
$$
and the centre of the next blowing-up $\s_3$ is
$C_2=\{ 0\}$, in this chart.

Of course, over $U_w$, $X_2$ lies in the union
of the following charts, for each of which we give a
defining equation $g_2=0$ and equations of the exceptional
hypersurfaces $H_1$, $H_2$:
\medskip

\begin{tabular}{cll}
Chart & \quad\ Equation of $X_2$ & \ Exceptional hypersurfaces\\[0.7ex] 
$U_{wy}$ & \ $z^d-x^{d-1}y^{d-1} \ = \ 0$ & \quad $H_1: w=0,\ H_2: y=0$\\
$U_{wx}$ & \ $z^d - x^{d-1} y^d \ = \ 0$  & \quad $H_1: w=0,\ H_2: x=0$\\ 
$U_{ww}$ & \ $z^d - x^{d-1}y^d w^{d-1} \ = \ 0$ & \quad $H_1=\emptyset,
\ H_2: w= 0$ 
\end{tabular}
\medskip

A calculation parallel to that above shows that, in each of these
charts, the centre of the next blowing-up is the origin.

To estimate $n(BM)$, we have to follow the branching of
the coordinate charts, as above, after each subsequent blowing-up.
We give the calculation in detail for two particular branches
and leave it to the reader to check (exploiting the
similar form of the data in the various charts) that
for no branch do we need a number of blowings-up $> 2d+j$.

\subsection{First branch} 
Charts $U_w$, $U_{wy}$, $U_{wyx}$, $U_{wyxw}$, $\ldots$.

\subsection*{Year three}  Let $X_3$ denote the strict transform
$X'_2$ of $X_2$ by $\s_3$.
In the chart $U_{wyx}$ of $M_3$, $\s_3$ is given by
the substitution $(x,xy,xz,xw)$.
Therefore $X_3\cap U_{wyx} = \{ g_3=0\}$, where
$$
g_3 (x,y,z,w) \ = \ z^d - x^{d-2} y^{d-1} .
$$
Let $a_3=0$.
Then $E(a_3)=\{ H_1,H_2,H_3\}$, where
$$
H_1:\ w=0,\quad
H_2:\ y=0,\quad
H_3:\ x=0
$$
(indexing the $H_i$ according to our previous convention).
We claim that
$$
\inv_X(a_3) \ = \ (d,1;\, 1,0;\, 0) .
$$
The first two pairs can be computed as before;
$\inv_1$ at $a_3$ has a codimension $1$ presentation
$\ucH_1 (a_3) = ( N_1(a_3),\cH_1(a_3),\cE_1(a_3))$,
where $N_1(a_3) = \{z=0\}$, $\cH_1(a_3) = \{ (x^{d-2}
y^{d-1},d),\, (w,1)\}$, $\cE_1(a_3)=\{ H_2,H_3\}$,
and $\inv_2$ has a codimension $2$ presentation
$\ucH_2(a_3)$ given by $N_2(a_3)=\{ z=w=0\}$,
$\cH_2(a_3)=\{ (x^{d-2}y^{d-1},d)\}$, $\cE_2(a_3)=
\cE_1(a_3) = \{ H_2,H_3\}$.
Therefore,
$$
\nu_3(a_3) \ = \ \mu_3(a_3) - \mu_{3H_2} (a_3) - \mu_{3H_3}
(a_3) \ = \ 0 ,
$$
and $\ucG_3(a_3) = \ucH_2(a_3)$ is a codimension $2$ presentation
of $\inv_{5/2}$ or $\inv_X$ at $a_3$; in particular, $D_3(a_3)^d=
x^{d-2}y^{d-1}$, in the notation of \S3.2.
Since $D_3(a_3)^d$ has order at least $d$ only at the origin,
the centre of the next blowing-up $\s_4$ in this chart is
$C_3=\{ 0\}$.

\subsection*{Year four}  Let $X_4 := X'_3$.
In the chart $V := U_{wyxw}$, $\s_4$ is given by the substitution
$(xw,yw,zw,w)$.
Therefore, $X_4\cap V$ is defined by
\begin{equation}
g_4(x,y,z,w) \ = \ z^d - x^{d-2} y^{d-1} w^{d-3} .
\end{equation}
Let $a_4=0$.  Then $E(a_4)=\{ H_2,H_3,H_4\}$, where
\begin{equation}
H_2:\ y=0,\quad
H_3:\ x=0,\quad
H_4:\ w=0 .
\end{equation}
($H_1\cap V=\emptyset$; the variable $w$ has been ``re-marked''
as $H_4$!).
Now $\inv_1(a_4)=(d,0)$ and $\inv_1$ at $a_4$ has a codimension
$1$ presentation $\ucH_1(a_4)$ given by $N_1(a_4)=\{z=0\}$,
$\cH_1(a_4) = \{ (x^{d-2}y^{d-1}w^{d-3}, d)\}$ and
$\cE_1(a_4) = E(a_4) = \{ H_2,H_3,H_4\}$.
Therefore,
$$
\inv_X(a_4) \ = \ (d,0;\, 0)
$$
and $\ucG_2(a_4) = \ucH_1(a_4)$ is a codimension 1 presentation
of $\inv_{3/2}$ or of $\inv_X$ at $a_4$; in particular,
$\cG_2(a_4) = ( D_2 (a_4)^d,d)$, where $D_2(a_4)^d=
x^{d-2} y^{d-1} w^{d-3}$.
Then $S_{\inv_X} (a_4)$ has three components as follows;
we order these components using the lexicographic ordering of
the quadruples shown (cf. \S1.6):
$$
\begin{array}{rcl}
\{ z = w = x = 0\} & \hphantom{x} & (0,0,1,1)\\
\{ z = w = y = 0\} & \hphantom{x} & (0,1,0,1)\\
\{ z = y = x = 0\} & \hphantom{x} & (0,1,1,0)
\end{array}
$$
The centre $C_4$ of the next blowing-up is given by
the maximum order; i.e., $C_5=\{ z=y=x=0\}$.
We are now in the combinatorial situation of \S3.3(4).
We need a number of blowings-up of the order of $2d$
(i.e., $2d+j$ blowings-up, where $j$ is independent of
$d$) to decrease the order $d$ of the strict transform
of $X$ along any branch of coordinate charts over $V$;
for example:
\medskip

\renewcommand{\arraystretch}{1.3}
\noindent
{\scriptsize
\begin{tabular}{ccclc}
{\small Year $j$} & {\small Chart} & 
{\small  Strict transform $g_j$} & 
{\small Exceptional divisors} & 
{\small Centre $C_j$}\\[0.7ex]
$5$ & $V_x$ & $z^d-x^{d-3}y^{d-1}w^{d-3}$ 
& $x: H_5,\ y: H_2,\ w: H_4$ 
& $\{z=w=y=0\}$\\
$6$ & $V_{xw}$ & $z^d-x^{d-3}y^{d-1}w^{d-4}$
& $x: H_5,\ y: H_2,\ w: H_6$ 
& $\{z=y=x=0\}$\\
$7$ & $V_{xwx}$ & $z^d-x^{d-4}y^{d-1}w^{d-4}$
& $x: H_7,\ y: H_2,\ w: H_6$
& $\{z=w=y=0\}$\\
etc. &  &  &  &
\end{tabular} }

\renewcommand{\arraystretch}{1}

\subsection{Second branch}
We will briefly make a calculation analogous to the preceding
or the branch of charts $U_w$, $U_{wy}$, $U_{wyw}$,
$U_{wywy}$, $\ldots$.
Years one and two are the same as in the branch above.

\subsection*{Year three}  In the chart $U_{wyw}$, $\s_3$ is given
by the substitution $(xw,yw,zw,w)$.
Therefore, $X_3\cap U_{wyw}$ is defined by $g_3=z^d-x^{d-1}
y^{d-1} w^{d-2}$.
Let $a_3=0$.
Then $E(a_3)=\{ H_2,H_3\}$ where $H_2$ and $H_3$ are defined
by $y$ and $w$, respectively ($H_1\cap U_{wyw}=\emptyset$.)
Then $\inv_1(a_3)=(d,0)$.
We calculate $N_1(a_3)=\{ z=0\}$, $\cH_1(a_3) = \{ (x^{d-1}
y^{d-1} w^{d-2},d)\}$, $\nu_2(a_3)=(d-1)/d$,
$s_2(a_3)=2$, $N_2(a_3)=\{ z=x=0\}$,
and $\cH_2(a_3) = \{ (y,1), (w,1)\}$.
Therefore,
$$
\inv_X(a_3) \ = \ \left( d , 0;\, 
\frac{d-1}{d},2;\, 1,0;\, 1,0;\, \infty\right)
$$
and $C_3=\{ 0\}$.

\subsection*{Year four}  
In $U_{wywy}$, $\s_4$ is given by $(xy,y,yz,yw)$;
therefore $X_4\cap U_{wywy}$ is defined by
$g_4 = z^d-x^{d-1} y^{2d-4} w^{d-2}$.
Let $a_4=0$.
Then $H_3$, $H_4$ are defined by $w$, $y$, respectively
($H_1$ and $H_2$ do not intersect this chart).
We calculate
$$
\inv_X (a_4) \ = \ \left( d,0;\, \frac{d-1}{d},1;\, 1,0;\, \infty\right) ,
$$
together with $N_1(a_4)=\{z=0\}$, $\cH_1(a_4)=\{ (x^{d-1} y^{2d-4}
w^{d-2},d)\}$, 
$N_2(a_4)=\{ z=x=0\}$,
and $\cH_2(a_4)=\{ (w,1)\}$.
Therefore, $C_4=\{ z=w=x=0\}$.

\subsection*{Year five}  In $U_{wywyw}$, $\s_5=(xw,y,zw,w)$.
Therefore $X_5\cap U_{wywyw}$ is given by $g_5=z^d-x^{d-1}
y^{2d-4} w^{d-3}$ and $y$, $w$ are exceptional divisors
representing $H_4$, $H_5$ (respectively).
At $a_5=0$, we calculate
$$
\inv_X(a_5) \ = \ \left( d,0;\, \frac{d-1}{d},0;\, 0\right) ,
$$
together with $N_1(a_5)=\{ z=0\}$, $\cH_1(a_5)=\{ (x^{d-1}
y^{2d-4} w^{d-3}, d)\}$, 
$N_2 (a_5) = \{ z=x=0\}$, $\cH_2 (a_5) = \{ (y^{2d-4} w^{d-2},1)\}=\cG_3
(a_5)$.
Therefore, $S_{\inv_X}(a_5)$ has two components, ordered as follows:
$$
\begin{array}{rcl}
\{z=y=x=0\} & \hphantom{x} & (0,0,0,1,0)\\
\{z=w=x=0\} &              & (0,0,0,0,1)\\
\end{array}
$$
(cf. \S4.3, year four).
The next blowing-up $\s_6$ has centre $C_5=\{ z=y=x=0\}$.
\medskip

Over $U_{wywyw}$, the strict transform $X_6$ of $X_5$ by $\s_6$
lies entirely in two charts, given as follows:
$$
\begin{array}{rcccl}
U_{wywywx} & \hphantom{x} & (x,xy,xz,w) & \hphantom{x} & 
g_6=z^d-x^{2d-5} y^{2d-4} w^{d-3} \\
U_{wywywy} & \hphantom{x} & (xy,y,yz,w) & \hphantom{x} &
g_6=z^d-x^{d-1} y^{2d-5} w^{d-3}
\end{array}
$$
It is easy to check that, following either branch, the order $d$
is decreased by a number of blowings-up of the order of $2d$.

\subsection{Villamayor's algorithm} 
We now reconsider our example (4.1) above using Villamayor's algorithm.
We will exhibit a sequence of points over which a number
of blowings-up of the order $9d$ is needed to decrease the order
$d$ of $X_1$ at $a_1$:
Following the first branch above, there is no change until year
four, so we reconsider our calculation from that point:
In each year $j$ below, $a_j$ denotes the origin of the chart
we study.

\subsection*{Year four}  Let $V := U_{wyxw}$ as in \S4.3, year
four, above.
Then $X_4\cap V$ is defined by 
$g_4 (x,y,z,w) = z^d - x^{d-2} y^{d-1} w^{d-3}$ (4.2),
and $E(a_4) = \{ H_2,H_3,H_4\}$, where
$H_2$, $H_3$ and $H_4$ are defined by $y$, $x$ and $w$,
respectively (4.3) (and $H_1\cap V=\emptyset$).
As before, we have $\inv_1^* (a_4) = (d,0)$.
(We use an asterisk to distinguish the invariant corresponding
to Villamayor's algorithm.)
But, following Villamayor's algorithm, we take
$$
\ucH_1(a_4) = ( N_1(a_4), \cH_1(a_4), \cE'_1(a_4)) ,
$$
where $N_1(a_4)=\{ z=0\}$, $\cH_1(a_4) = \{ (x^{d-2} y^{d-1} w^{d-3},d
)\}$, and $\cE'_1(a_4)=\emptyset$ (the latter because this year four
is the year of birth of the value $(d,0)$ of $\inv_1^*(a_4)$).
Then $E^2(a_4) = E(a_4) = \{ H_2,H_3,H_4\}$ and
$s_2(a_4)=3$.
Therefore, we take $\ucH_2(a_4) = ( N_2(a_4), \cH_2(a_4), 
\emptyset )$, where $N_2(a_4)=\{ z=w=0\}$, $\cH_2(a_4)
= \{ ((x,1), (y,1)\}$, and get
$$
\inv_X^*(a_4) = \left( d,0;\, \frac{3d-6}{d},3;\, 1,0;\,
1,0;\, \infty\right) ;
$$
the centre of the next blowing-up $\s_5$ in this chart is
$C_4=\{ 0\}$.

\subsection*{Year five}  In $V_w$, $g_5=z^d-x^{d-2}y^{d-1}w^{2d-6}$.
In a fashion similar to year four, we compute
$$
\inv_X^* (a_5) = \left( d,0;\, \frac{2d-3}{d},3;\,
1,0;\, 1,0;\, \infty\right) .
$$
(Here $\cE'_1(a_5)=\{ H_5\}$, where $H_5=\{w=0\}$, so we
factor $w^{2d-6}$ from $\cH_1(a_5)$ to get $\nu_2(a_5)=
(2d-3)/d$.)
Therefore, $C_5=\{0\}$.

\subsection*{Year six}  In $V_{ww}$, $g_6=z^d-x^{d-2}y^{d-1}w^{3d-9}$.
We compute
$$
\inv_X^* (a_6) = \left( d,0;\, \frac{2d-3}{d},2;\,
1,1;\, 1,0;\, \infty\right) .
$$
(where $s_2(a_6)=2$ counts $H_2$: $y=0$ and $H_3$: $x=0$,
and $s_3(a_6)=1$ counts $H_6$: $w=0$), together with
$N_1(a_6)=\{z=0\}$, $\cH_1(a_6)=\{ (x^{d-2} y^{d-1} w^{3d-9},d)\}$,
$N_2(a_6)=\{z=y=0\}$, $\cH_2(a_6)=\{ (x,1)\}$, $N_3(a_6)=\{z=y=x=0\}$,
and $\cH_3(a_6)=\{ (w,1)\}$. In particular, $C_6=\{ 0\}$.

\subsection*{Year seven}  
In $V_{www}$, $g_7=z^d-x^{d-2}y^{d-1}w^{4d-12}$.
Then
$$
\inv_X^* (a_7) = \left( d,0;\, \frac{2d-3}{d},2;\,
1,0;\, \infty\right) ,
$$
and $N_2(a_7)=\{ z=y=0\}$, $\cG_3(a_7)=\cH_2(a_7)
= \{ (x,1)\}$, so that $C_7= \{ z=y=x=0\}$.
\medskip

We go on as follows. (In this table, ``chart $x$'' means $U_x$,
where $U$ denotes the chart in the previous year. For example,
the chart in year eight is $V_{wwwx}$.)
\medskip

\renewcommand{\arraystretch}{1.5}
\noindent
{\tiny
\begin{tabular}{ccclc}
{\footnotesize Year $j$} & {\footnotesize Chart} & 
{\footnotesize Strict transform $g_j$} & 
{\footnotesize \qquad\quad $\inv_X^*(a_j)$} & 
{\footnotesize Centre $C_j$}\\[0.7ex]
$8$ & $x$ & $z^d-x^{d-3}y^{d-1}w^{4d-12}$ 
& $(d,0;\, \frac{d-1}{d},3;\, 1,0;\, 1,0;\, \infty)$
& $\{0\}$\\
$9$ & $w$ & $z^d-x^{d-3}y^{d-1}w^{5d-16}$
& $(d,0;\, \frac{d-1}{d},2;\, 1,1;\, 1,0;\, \infty)$
& $\{0\}$\\
$10$ & $w$ & $z^d-x^{d-3}y^{d-1}w^{6d-20}$
& $(d,0;\, \frac{d-1}{d},2;\, 1,0;\, \infty)$
& $\{z=y=x=0\}$\\
$11$ & $x$ & $z^d-x^{d-4}y^{d-1}w^{6d-20}$
& $(d,0;\, \frac{d-1}{d},1;\, 7d-24,2;\, 1,0;\, \infty)$
& $\{0\}$\\
$12$ & $w$ & $z^d-x^{d-4}y^{d-1}w^{7d-25}$
& $(d,0;\, \frac{d-1}{d},1;\, d-4,2;\, 1,0;\, \infty)$
& $\{z=y=x=0\}$\\
$13$ & $w$ & $z^d-x^{d-4}y^{d-1}w^{8d-30}$
& $(d,0;\, \frac{d-1}{d},1;\, d-4,1;\, \infty)$
& $\{z=y=x=0\}$\\
$14$ & $x$ & $z^d - x^{d-5} y^{d-1} w^{8d-30}$
& $(d,0;\, \frac{d-1}{d},1;\, 0)$ & $\{z=w=y=0\}$
\end{tabular} }
\medskip

\renewcommand{\arraystretch}{1}

As an aid to computing the entries in the above table, we
note that, in year eleven, we can take 
$\cG_2(a_{11})=\{ (y,1),\, (x^{d-4}w^{6d-20},1)\}$,
and in year thirteen, we can take $N_2(a_{13})=\{ z=y=0\}$,
$\cH_2(a_{13})=\{ (x^{d-4} w^{8d-30},\, 1)\}$.

In year fourteen, we can take $N_2(a_{14})=\{ z=y=0\}$ 
and $\cH_2(a_{14})=\{ (x^{d-5} w^{8d-30},1)\}$;
$S_{\inv_X}(a_{14})$ has two components $\{ z=y=x=0\}$ and
$\{ z=w=y=0\}$; $C_{14}$ is the latter.

We now follow alternately the ``$w$-chart'' or ``$x$-chart'',
finding alternately $C_{\cdot} = \{z=y=x=0\}$ or
$C_{\cdot}=\{ z=w=y=0\}$ until, in year $2d+4$, we have
$g_{2d+4}=z^d-y^{d-1}w^{7d-25}$ and $C_{2d+4}=
\{ z=w=y=0\}$.
Now following the ``$w$=charts'', we reduce the order $d$
after $7d-25$ further blowings-up (i.e., in year $9d-21$).
(Throughout the calculation above, we assume that $d\ge 5$.)

\section{Equivalence of presentations}

Our purpose in this section is to elucidate the ideas of
equivalence of presentations introduced in Section 2.
In particular, Corollaries 5.12 and 5.13 below imply Theorem 2.6,
which is used in the inductive
definition of $\inv$ in \S3.2, to pass from  
a codimension $r$ presentation of $\inv_{r+1/2}$ to a 
presentation in codimension $r+1$. 
It is convenient to use transformation formulas for differential
operators introduced by Hironaka \cite[Sect.~8]{Hid} and developed
by Giraud \cite{giraud} and Villamayor and Encinas \cite{V2,EVacta}. 
(An alternative approach is given in
\cite[Propositions 4.12, 4.19]{BMinv}.) The formulas describe
the way that the partial derivatives of a regular function
transform by admissible blowings-up. The treatment below differs from 
that of \cite{V2,EVacta} in our use of coordinate charts (as
in \cite{BMinv, BMsam}) and of transformation formulas that account
also for the effect of exceptional blowings-up.

\subsection{Transformation of differential operators}
Let $U$ denote a coordinate chart in a manifold $M$,
with coordinate system
$(x_1,\ldots,x_n)$. (The $x_i$ are regular functions on
$U$. In the case of schemes of finite type, ``coordinate
chart'' means ``\'etale (or regular) coordinate chart'',
as defined in \cite[\S3]{BMinv}. See also \cite[\S2]{BMsam}
for coordinate charts, more generally.)
Let $a=0$.
Consider a blowing-up $\s$ with centre
$$
C \ = \ Z_I \ := \ \{ x_i = 0,\ i\in I\} ,
$$
where $I\subset \{ 1,\ldots,n\}$.
If $i\in I$, then $\s$ is given in the chart $U_i=U_{x_i}\subset
\s^{-1} (U)$ (cf. Section 4) by the formulas $x_i=y_i=y_\exc$,
$x_j=y_iy_j$ if $j\in I\backslash \{ i\}$, and $x_j=y_j$ if
$j\not\in I$.
(For $j\in I\backslash \{ i\}$, $y_j$ does not necessarily
vanish at $a'\in \s^{-1} (a)$.)
The following lemma is a simple calculation. 

\begin{lemma} 
Let $f\in\cO_a = \cO_{M,a}$. Suppose that $d \leq \mu_{C,a}(f)$.
Let $i \in I$. Then 
\begin{eqnarray*}
\frac{1}{y_i^{d-1}} \left( \frac{\p f}{\p x_j} \circ\s \right) & =
& y_i \frac{\p}{\p y_j} \left( \frac{f\circ \s}{y_i^d}\right),
\qquad \mbox{if }\, j\not\in I;\\
\frac{1}{y_i^{d-1}} \left( \frac{\p f}{\p x_j} \circ\s\right) & =
& \frac{\p}{\p y_j} \left( \frac{f\circ\s}{y_i^d}\right),
\qquad \mbox{if }\, j\in I\backslash \{ i\};\\ 
\frac{1}{y_i^{d-1}} \left( \frac{\p f}{\p x_i}\circ\s \right) & =
& d \frac{f\circ\s}{y_i^d} + y_i \frac{\p}{\p y_i}
\left( \frac{f\circ\s}{y_i^d}\right) - \sum_{j\in I\backslash\{ i\}}
y_j \frac{\p}{\p y_j} \left( \frac{f\circ\sigma}{y_i^d}\right) .
\end{eqnarray*}
\end{lemma}

Let $\cE(a)$ denote a collection of smooth hypersurfaces
passing through $a$.
Suppose that $\cE(a)$ and $C$ simultaneously have only normal
crossings.
Then we can choose coordinates so that $C$ has the form $Z_I$
as above, and also, for each $H\in\cE(a)$, $\cI_{H,a}$ is
generated by $x_j$, for some $j$; we will write $x_j=x_H$.

The following lemma will be used to study the effect 
of an exceptional blowing-up.
Let $H_1,H_2\in\cE(a)$.
Suppose $x_1=x_{H_1}$, $x_2=x_{H_2}$.
Let $\s$ denote the blowing-up with centre $C=H_1\cap H_2$.
In the chart $U_1=U_{x_1}$, $\s$ is given by $x_1=y_1$,
$x_2=y_1y_2$, and $x_j=y_j$ if $j>2$.

\begin{lemma}
Let $f\in\cO_a$, $\mu_a(f)\ge 1$.
Then
\begin{eqnarray*}
\left( x_2 \frac{\p f}{\p x_2}\right)\circ \s & = &
y_2 \frac{\p(f\circ\s)}{\p y_2} ;\\
\left( x_1 \frac{\p f}{\p x_1}\right)\circ \s & = &
y_1 \frac{\p(f\circ\s)}{\p y_1} - y_2 \frac{\p(f\circ\s)}{\p y_2} ;\\ 
\frac{\p f}{\p x_j} \circ\s & = & \frac{\p (f\circ\s)}{\p y_j} ,
\qquad \mbox{if }\, j>2.
\end{eqnarray*}
\end{lemma}

Consider a coordinate system of the form
\begin{equation}
x \ = \ (\xi,u) \ = \ (\xi_1,\ldots,\xi_r,\, u_1,\ldots,u_s) ,
\end{equation}
where the $\xi_i$ are precisely the $x_H$, $H\in\cE(a)$.
We will say that the variables $u_1,\ldots,u_s$ are
{\it complementary to} $\cE(a)$.

\begin{lemma}
Let $f\in\cO_a$. Then the  ideal generated by
\begin{eqnarray*}
x_H \frac{\p f}{\p x_H} , & \hphantom{=} & H\in\cE(a) ,\\
\frac{\p f}{\p u_j} ,     & \hphantom{=} & 1\le j\le s ,
\end{eqnarray*}
\smallskip
is independent of the choice of generators $x_H$ of
the ideals of $H\in\cE(a)$ and the choice of complementary
variables $u$.
\end{lemma}

\begin{proof}
If we change the generators $\xi_i=x_H$ of the ideals of the
$H\in\cE(a)$, say to $\eta_i$, then we have $\eta_i=\la_i\xi_i$,
where each $\la_i=\la_i(\xi,u)$ is a unit.
Consider such a change of generators, as well as new complementary
variables $v=(v_1,\ldots,v_s)$.
Then we can write
$$
(\eta,v) \ = \ \left( \la_1(\xi,u)\xi_1,\ldots,\la_r(\xi,u)\xi_r,\, 
v_1(\xi,u),\ldots,v_s(\xi,u)\right) .
$$
Therefore, for each $i=1,\ldots,r$,
$$
\xi_i\frac{\p}{\p \xi_i} \ =
\ \frac{1}{\la_i} \frac{\p\eta_i}{\p\xi_i} \eta_i \frac{\p}{\p\eta_i} +
\sum_{k\ne i} \frac{\xi_i}{\la_k} \frac{\p\la_k}{\p\xi_i} 
\eta_k \frac{\p}{\p\eta_k} +
\sum_{\ell=1}^s \xi_i \frac{\p v_\ell}{\p\xi_i} \frac{\p}{\p v_\ell} ,
$$
and, for each $j=1,\ldots,s$,
$$
\frac{\p}{\p u_j} \ = \ \sum_{k=1}^r \frac{1}{\la_k} 
\frac{\p\la_k}{\p u_j}
\eta_k \frac{\p}{\p \eta_k} + \sum_{\ell=1}^s \frac{\p v_\ell}{\p u_j}
\frac{\p}{\p v_\ell} .
$$
\end{proof}

\subsection{Passage to codimension $+1$}
Let $M$ denote a manifold and let $a \in M$. Let $N(a)$
denote a germ of a submanifold of $M$ at $a$; say $p =
\codim \, N(a)$. 

\begin{lemma}
Let $z_1,\ldots,z_p \in \cO_a = \cO_{M,a}$ denote generators
of $\cI_{N(a)}$ (so that $z_1,\ldots,z_p$ have
linearly independent gradients). Let $f \in \cO_a$ and let
$d \in \IN$. Then $\mu_a(f) \geq d$ if and only if
$$
\mu_a \left( \frac{\p^{|\beta|} f}{\p z^\beta} \Big|_{N(a)} \right)
\ \geq \ d - |\beta|, \quad \mbox{for all }\, \beta \in \IN^p, \, 
|\beta| \leq d - 1.
$$
\end{lemma}

Let $\cE(a)$ denote a collection of smooth hypersurfaces
passing through $a$, such that $N(a)$ and $\cE(a)$ simultaneously
have only normal crossings, and $N(a) \not\subset H$, for all
$H \in \cE(a)$. Then we can choose a coordinate system of the
form (5.1) for $N(a)$.

\begin{definitions}
Choose a coordinate system $x=(\xi,u)$ for $N(a)$ as in (5.1).
Let $f \in \cO_{N(a)}$ and let $\mu_f \in \IQ$ such that
$\mu_a(f) \geq \mu_f$. Set
\begin{eqnarray*}
\D(f,\mu_f) & := & \bigg\{ (f,\,\mu_f-1),\ 
\left( x_H \frac{\p f}{\p x_H},\, \mu_f-1 \right),\, \mbox{for all }\,
H\in\cE(a),\\
	    &    & 
\quad \left( \frac{\p f}{\p u_j},\, \mu_f-1 \right),
\, \mbox{for all } j=1,\ldots,s 
\bigg\},\\
\left(\D (f,\mu_f)\right) & := & \mbox{the ideal in } \cO_N(a)
\mbox{ generated by }\\
  &    & \quad \left\{ h:\ (h,\mu_h) \in \D(f,\mu_f) \right\}, 
\\
\cD(f,\mu_f) & := & \left\{ (f,\mu_f)\right\} \cup \D(f,\mu_f);\\
S_{(f,\mu_f)} & := & \left\{ x\in N(a):\ \mu_x(f)\ge \mu_f\right\},
\\
S_{\cD(f,\mu_f)} & := & \left\{ x\in N(a):\ \mu_x(h)\ge \mu_h,\ 
\mbox{for all } (h,\mu_h)\in\cD(f,\mu_f)\right\}
\end{eqnarray*}
($S_{(f,\mu_f)}$ and $S_{\cD(f,\mu_f)}$ make sense as germs at $a$).
Note that an element $(h,\mu_h)$ in $\cD(f,\mu_f)$ with
$\mu_h\le 0$ imposes no condition in $S_{\cD(f,\mu_f)}$. 
\end{definitions}

\begin{corollary}
Let $f \in \cO_{N(a)}$ and $\mu_f \in \IQ$ such that
$\mu_a(f) \geq \mu_f$. Then
$$
S_{(f,\mu_f)} \ = \ S_{\cD(f,\mu_f)} .
$$
\end{corollary}

Now let
$$
\ucF(a) = \left( N(a), \, \cF(a),\, \cE(a)\right)
$$
denote a presentation (of codimension $p$) at $a$; say
$\cF(a) = \{ (f,\mu_f)\}$.
{\it In the case that all} $\mu_f$ {\it are equal}, we will write
$\big( \cF(a)\big)$ for the ideal in $\cO_{N(a)}$ generated by
$\big\{ f:\ (f,\mu_f)\in\cF(a)\big\}$.

\begin{definitions}
Choose coordinates as in Definitions 5.5 above. Set
\begin{eqnarray*}
\D\left( \cF(a)\right) & := & \bigcup_{\cF(a)} \D(f,\mu_f) ,\\ 
\cD\left( \cF(a)\right) & := & \cF(a) \cup \D(\cF(a)) \ 
= \ \bigcup_{\cF(a)} \cD (f,\mu_f) ,\\ 
\ucD\left(\ucF(a)\right) & := & 
\left( N(a), \cD\left( \cF(a)\right) , \cE(a)\right) .
\end{eqnarray*}
\end{definitions}

Let $\s$ be a morphism of type (i), (ii) or (iii) (cf. \S2.2)
and let $\ucF(a')=( N(a'),\cF(a'),\cE(a'))$ or
$\ucF(a)' = ( N(a)',\cF(a)',\cE(a)')$ denote the
transform of $\ucF(a)$ by $\s$ at a point $a'\in\s^{-1}(a)$ as
in \S2.2.
(It will be convenient to use both notations for the transform.)
We will also write $\D(f,\mu_f)'$, $\D(\cF(a))'$ and
$\cD( \cF(a))'$ for the analogous transforms of
$\D(f,\mu_f)$, $\D( \cF(a))$ and $\cD( \cF(a))$,
respectively; for example,
$$
\D(f,\mu_f)' \ := \ \left\{ (h',\mu_{h'}):\ (h,\mu_h)\in
\D(f,\mu_f)\right\} ,
$$
where $(h',\mu_{h'})$ is given by the transformation rules in \S2.2.

\begin{lemma}
$\D(f,\mu_f)'\, \subset \, ( \D(f',\mu_{f'}))$.
\end{lemma}

(This is a minor abuse of notation; we mean that,
for all $(h',\mu_{h'})\in\D (f,\mu_f)'$, $\, h'\in
(\D (f',\mu_{f'}))$.)

\begin{proof}
For transformations of types (i) or (iii) (admissible or
exceptional blowings-up), this follows from the formulas
in Lemmas 5.1 and 5.2.
For a transformation of type (ii) (product with a line),
it is trivial.
\end{proof}

\begin{theorem}
Let $(f,\mu_f)\in\cF(a)$.
Then
$$
S_{(f',\mu_{f'})} \ = \ S_{\cD(f,\mu_f)'}
$$
after any transformation of type (i), (ii) or (iii) (and, in fact,
after any sequence of transformations 
of types (i), (ii) and (iii)).
\end{theorem}

\begin{proof}
For a given sequence of transformations of types (i), (ii) and
(iii), write $f^{(0)}=f$, $\D(f,\mu_f)^{(0)}=\D(f,\mu_f)$ and
$\cD(f,\mu_f)^{(0)}=\cD(f,\mu_f)$, and recursively define
$f^{(k+1)} := (f^{(k)})'$, $\mu_{f^{(k+1)}}=\mu_{(f^{(k)})'}=\mu_f$,
$\D(f,\mu_f)^{(k+1)} := (\D(f,\mu_f)^{(k)})'$ and
$\cD(f,\mu_f)^{(k+1)} := ( \cD(f,\mu_f)^{(k)})'$, for
all $k\ge 0$.

We have $(f,\mu_f)\in\cD(f,\mu_f)$, so that for all $k$,
$(f^{(k+1)},\mu_{f^{(k+1)}}) \in\cD (f,\mu_f)^{(k+1)}$ and
$$
S_{\cD(f,\mu_f)^{(k+1)}} \ \subset
\ S_{(f^{(k+1)},\mu_{f^{(k+1)}})} .
$$

By Lemma 5.8, $\D(f,\mu_f)'\subset( \D(f',\mu_{f'}))$.
Assume that $\D(f,\mu_f)^{(k)}\subset(\D(f^{(k)},\mu_{f^{(k)}}))$,
by induction.
Consider $(h,\mu_h)\in\D (f,\mu_f)^{(k)}$.
Then $h'\in (\D(f^{(k+1)},\mu_{f^{(k+1)}}))$, again by
Lemma 5.8.
Therefore, for all $k$,
$$
\D(f,\mu_f)^{(k+1)} \ \subset \ \left( \D(f^{(k+1)},\mu_{f^{(k+1)}}\right)
$$
and hence
$$
S_{\cD(f^{(k+1)},\mu_{f^{(k+1)}})} \ \subset
\ S_{\cD(f,\mu_f)^{(k+1)}} .
$$

But
$$
S_{\cD(f^{(k+1)},\mu_{f^{(k+1)}})} \ = 
\ S_{(f^{(k+1)},\mu_{f^{(k+1)}})},
$$
by Corollary 5.6, so the result follows.
\end{proof}

\begin{corollary}
$\ucF(a)$ and $\ucD(\ucF(a))$ are equivalent with
respect to transformations of types (i), (ii) and (iii).
(See Definitions 2.1.)
\end{corollary}

\begin{definitions}
Write $\cD^0 (\cF(a)) = \D^0
(\cF(a)) = \cF(a)$, and, for all $d=0,1,2,\ldots$, set
\begin{eqnarray*}
\D^{d+1} (\cF(a)) & := & \D\left(\D^d\left(\cF(a)\right)\right) ,\\
\cD^{d+1} ( \cF(a)) & := & 
\bigcup_{q=0}^{d+1}\D^q(\cF(a)) \ =
\ \cD\left(\cD^d\left(\cF(a)\right)\right) .
\end{eqnarray*}
\end{definitions}

\begin{corollary}
Let $z\in\cO_{N(a)}$ and let $d$ be a positive integer.
Suppose that $\mu_a(z)=1$ and 
that $z\in\left(\D^{d-1}\left(\cF(a)\right)\right)$.
Then, after any sequence of transformations of types (i), (ii) or
(iii), $z'\in\left(\D^{d-1}\left(\cF(a')\right)\right)$, and
$$
S_{\ucF(a')} \ \subset \ V(z') \ \subset \ N(a') .
$$
\end{corollary}

The first assertion is a consequence of Lemma 5.8 and the
second is a consequence of Corollary 5.10.
From Corollaries 5.10 and 5.12, we deduce:

\begin{corollary}
Under the hypotheses of Corollary 5.12, set
\begin{eqnarray*}
N_{+1}(a) & := & V(z) \ \subset \ N(a) ,\\
\cH(a) & := & \cD^{d-1} ( \cF(a)) |_{N_{+1}(a)} , \\
\ucH(a) & := & ( N_{+1}(a),\cH(a),\cE(a)) .
\end{eqnarray*}
Then the presentations $\ucF(a)$ and $\ucH(a)$ are equivalent
with respect to transformations of types (i), (ii) and (iii).
\end{corollary}

It is clear that $\ucF(a)$ and $\ucH(a)$ are, in fact,
semicoherent equivalent. (See \S2.5). Theorem 2.6 follows
from Corollaries 5.12 and 5.13 --
this is the basis of our constructive definition of
$\inv(\cdot)$ by induction on codimension (\S3.2).

\subsection{On the notion of equivalence}
We conclude this section by showing that, in contrast
to Theorem 2.4, it is not true that the
$\mu_{r+1,H}(a)$ and $\nu_{r+1}(a)$, $r>0$ (even
as occuring in Villamayor's invariant) in general
depend only on the equivalence class of a presentation
of $\inv$ at $a$ with respect to transformations of
types (i) and (ii). (See Remarks 2.5.)

\begin{example}
Let $X$ denote the surface in affine 3-space $U$ defined by
$$
z^d - x^{d-1} y^d \ = \ 0 ,
$$
where $d\ge 2$.
We will use the notational conventions of Section 4.
Let $\s_1$ denote the blowing-up of $U$ with centre
$C_0=\{ 0\}$.
In the chart $U_x$, the strict transform $X_1$ of $X$ is given
by $g_1=0$, where
$$
g_1(x,y,z) \ = \ z^d - x^{d-1} y^d .
$$
Let $\s_2$ be the blowing-up of $U_x$ with centre $C_1=\{0\}$.
In the chart $U_{xy}$, the strict transform $X_2$ of $X_1$
is given by $g_2=0$, where
$$
g_2(x,y,z) \ = \ z^d - x^{d-1} y^{d-1} .
$$

Let $a=0$ in $U_{xy}$.
Following either the algorithm of the authors or that of
Villamayor,
we have $E^1(a)=\emptyset$ and $\cE_1(a)=E(a)=\{H_1,H_2\}$,
where the exceptional hyperplanes $H_1$ and $H_2$ are given
by $x=0$ and $y=0$, respectively.
Then
$$
\mu_2(a) = \frac{2d-2}{d},\quad
\mu_{2H_1}(a) = \frac{d-1}{d},\quad
\mu_{2H_2}(a) = \frac{d-1}{d} ,
$$
and
$$
\inv_X (a) = (d,0;0) .
$$

At $a$, $\inv_1$ has a 
codimension $1$ presentation $\ucH_1(a)=
(N_1(a), \cH_1(a),\cE_1(a))$, where $N_1(a)=\{ z=0\}$
and $\cH_1(a) = \{ (x^{d-1} y^{d-1},d)\}$.
By Theorem 2.4, $\mu_{2H_1}(a)$ and $\mu_{2H_2}(a)$
depend only on $[\ucH_1(a)]$.
But it is not true that they depend only on $[\ucH_1(a)]_{{\rm (i,ii)}}$:

Recall that, for each $H\in\cE_1(a)$, $\mu_{2H}(a)=
\mu_{\ucH_1 (a),H} = \min_{\cH_1(a)}(\mu_{H,a}(h)/\mu_h)$.
It follows from Lemmas 5.1 and 5.8 above
that $\ucH_1(a)$ is equivalent with respect to transformations
of types (i) and (ii) to a presentation $\ucC(\ucH_1(a)) =
( N_1(a),\cC(\cH_1(a)),\cE_1(a))$, where
$$
\cC(\cH_1(a)) \ := \ \cH_1(a)\bigcup 
\left\{ \left(\frac{\p h}{\p x_i},\,
\mu_h-1\right):\ (h,\mu_h)\in\cH_1(a),\ i=1,\ldots,m\right\},
$$
where $(x_1,\ldots,x_m)$ denotes the coordinates of $N_1(a)$.
In our example,
$$
\cC( \cH_1(a)) \ = 
\ \big\{ (x^{d-1} y^{d-1},d),\ (x^{d-2}y^{d-1},d-1),\ 
(x^{d-1} y^{d-2}, d-1)\big\} ,
$$
so that
$$
\mu_{\ucC(\ucH_1(a)),H_1} \ = \ \frac{d-2}{d-1} 
\ = \ \mu_{\ucC(\ucH_1(a)),H_2} .
$$
\end{example}

\section{Applications of the desingularization principle}

The main topics of this section -- universal embedded
desingularization (of spaces that are not necessarily
embedded), comparison of weak and strict transforms, and
simultaneous desingularization of parametrized families 
-- are introduced individually in the subsections below.
The point of view is that of the general desingularization
principle, Theorem 1.14. Particular results concerning 
embedded desingularization (e.g., Theorems 6.1 and 6.8)
depend on the fact that the Hilbert-Samuel function satisfies
the Hypotheses 1.10 that are needed
to apply the desingularization principle (see Example 1.13(2)). 
In Theorem 6.18, we show that the Hilbert-Samuel function also
satisfies the stronger Hypotheses 6.13 below that are needed to 
extend the desingularization principle to parametrized 
families (Theorem 6.15). Theorems 6.8 and 6.15 thus generalize
results of \cite{BV} and \cite{ENV}, respectively, which are
tied to the weak embeddeded desingularization algorithm of
\cite{EVweak} (Example 1.19(4)). 

\subsection{Universal desingularization}
Let $\io = \io_X$ (or $\io_\cJ$) denote a local invariant of
spaces $X$ (or ideals of finite type $\cJ$; see \S1.3) satisfying the
Hypotheses 1.10. The invariant $\inv(\cdot) = \inv_X(\cdot)$ 
or $\inv_\cJ(\cdot)$, and therefore the desingularization
algorithm depend {\it a priori} on the codimension of a
presentation of $\inv_{1/2} = \io$ (see Hypotheses 1.10(3)) and,
in particular, on the dimension of the ambient manifold $M$.
We avoided these issues in \S3.2 by making the following
simplifying assumptions in Hypothesis 1.10(3): 
(1) $\io$ admits a semicoherent 
presentation $\ucG(a)= ( N(a),\cG(a),\emptyset)$ of 
codimension $0$ at every point $a \in M$; (2) $\mu_{\ucG(a)} = 1$.
(To begin the inductive construction, we then showed 
that, for any sequence of $\io$-admissible 
transformations (1.1) (or (1.2)), and all $a \in M_j$,
$j = 0, 1, \dots$, the invariant $\io$ 
admits a semicoherent presentation
$\ucC(a)=( N_1(a),\cC(a),\cE_1(a))$ of codimension $1$ at $a$,
where $\mu_{\ucC(a)} \geq 1$ (Corollary 3.4).) 

In \S6.1.1 below, we show how to modify 
the local inductive construction
(\S3.2) so that it applies without the
simplifying assumptions above. We assume only that $\io$ admits
a semicoherent presentation $\ucC(a)=( N(a),\cC(a),\emptyset)$
at every point $a \in M$, of codimension $p(a) \geq 1$
(Hypothesis 1.10(3); see Remark 6.2 below).

For example, \cite[Theorems 9.4, 9.6]{BMinv} provide such
a semicoherent presentation of the Hilbert-Samuel function
with variable codimension $p(a)$. The largest codimension of a
presentation of the Hilbert-Samuel function at $a$ is not
determined uniquely by $H_{X,a}$ \cite[Remarks 9.15(1)]{BMinv}.
This is why it is important to modify the local inductive
construction to ensure that $\inv(a)$ does not depend
on the codimension of a presentation of $\io$ at $a$.

The invariant $\inv(\cdot)$ and therefore the desingularization
algorithm, as described in \S6.1.1, nevertheless still depend
on the dimension of the ambient manifold $M$. In \S6.1.2, we show
that $\inv_X(\cdot)$ can be made independent of the embedding
space of $X$, by a simple variation in the definition. 
Given a scheme of finite type or an analytic space $X$ (not
necessarily globally embedded), then $\inv_X(\cdot)$ can be
defined in this way using any local embedding $X \big| U 
\hookrightarrow M$ over an open subset $U$ of $X$. As a result,
we obtain the following universal embedded desingularization theorem 
(cf. \cite[Theorem 13.2]{BMinv}) for (not necessarily embedded) spaces
$X$. We include the argument here both as an illustration of the
desingularization principle and to correct an error in 
\cite[Remarks 9.15(3)]{BMinv} that is illustrated by an example
of Encinas \cite{Eemb}. 

\begin{theorem} {\rm (1)} There is a finite sequence of blowings-up
$\s_{j+1}:\ X_{j+1} \to X_j$, where $X_0 = X$, such that, for any
local embedding $X|_U \hookrightarrow M$ (over an open subset $U$
of $|X|$), the sequence of blowings-up $\s_{j+1}$ restricted to
the inverse images of $U$ is induced by embedded desingularization
of $X|_U$ in the sense of Example 1.19(2).

{\rm (2)} The desingularization is {\rm universal} in the sense that,
to each $X$ we associate a morphism $\s_X:\ X' \to X$ such that
\begin{enumerate}
\item[(i)] $\s_X$ is a composite of a finite sequence of
blowings-up as in (1).
\item[(ii)] If $\vp:\ X\big|_U \to Y\big|_V$ 
is an isomorphism over
open subsets $U, V$ of two spaces $X, Y$ (respectively),
then there is an isomorphism $\vp':\ X'\big|_{\s_X^{-1}(U)} \to
Y'\big|_{\s_Y^{-1}(V)}$ such that $\s_Y \circ \vp' = \vp \circ \s_X$.
(The lifting $\vp'$ of $\vp$ is necessarily unique.) In fact,
$\vp$ lifts to isomorphisms throughout the desingularization
towers.
\end{enumerate}
\end{theorem}

\subsubsection{Independence of the codimension of a presentation}
Assume that $\io = \io_X$ (or $\io_\cJ$) admits a semicoherent
presentation $\ucC(a)=( N(a),\cC(a),\emptyset)$ at every point
$a \in M$, of codimension $p(a) \geq 1$. Consider any sequence
of $\io$-admissible transformations (1.1) (or (1.2)). We use
the notation of \S3.2. Let $a \in M_j$. Then 
$\inv_{1/2} = \io$ admits a semicoherent presentation
$$
\ucC_1(a) \ = \ ( N_1(a), \cC_1(a), \cE_1(a))
$$
at $a$, of codimension $p(a) \geq 1$. As before, we define
$\inv_1(a) := (\io(a), s_1(a) )$ and set 
$$
\ucH_1(a) \ = \ ( N_1(a), \cH_1(a), \cE_1(a)) ,
$$
where 
$$
\cH_1(a) := \cC_1(a) \cup 
\left(E^1(a)\big|_{N_1(a)} , 1 \right). 
$$
Then $\ucH_1(a)$ is a 
semicoherent presentation of $\inv_1$ at $a$, of codimension $p(a)$.

If $p(a) = 1$, we continue the inductive definition of $\inv(a)$
exactly as in \S3.2. If $p(a) > 1$, however, we consider the
following variation of the definition given in \S3.2. Recursively,
for each $r = 1, \ldots, p(a)-1$, we define
\begin{eqnarray*}
\nu_{r+1} (a) &:=& 1,\\
\inv_{r+1/2}(a) &:=& ( \inv_r(a); \nu_{r+1}(a)),
\end{eqnarray*}
$E^{r+1}(a)$ and $\cE_{r+1}(a)$ as before, $s_{r+1}(a) := 
\# E^{r+1}(a)$, and
$$
\inv_{r+1}(a) \ := \ ( \inv_{r+1/2}(a), s_{r+1}(a)).
$$
We set $N_{r+1}(a) := N_r(a) = N_1(a)$, 
$\cC_{r+1}(a) := \cH_r(a)$,
\begin{eqnarray*}
\cH_{r+1}(a) &:=& \cC_{r+1}(a) \cup 
\left(E^{r+1}(a)\big|_{N_1(a)}, 1 \right),\\
\ucC_{r+1}(a) &=& ( N_{r+1}(a), \cC_{r+1}(a), \cE_{r+1}(a)),\\
\ucH_{r+1}(a) &=& ( N_{r+1}(a), \cH_{r+1}(a), \cE_{r+1}(a));
\end{eqnarray*}
then $\ucC_{r+1}(a)$ (respectively, $\ucH_{r+1}(a)$) is a
codimension $p(a)$ presentation of $\inv_{r+1/2}$ (respectively,
of $\inv_{r+1}$ at $a$. Finally, 
$$
\inv_{p(a)}(a) \ = \ ( \io(a), s_1(a); 1, s_2(a);
\ldots ; 1, s_{p(a)}(a) )
$$
and 
$$
\ucH_{p(a)}(a) \ = \ ( N_{p(a)}(a), \cH_{p(a)}(a), \cE_{p(a)}(a) )
$$
is a codimension $p(a)$ presentation of $\inv_{p(a)}$ at $a$,
where $N_{p(a)}(a) = N_1(a)$,
\begin{eqnarray*}
\cH_{p(a)}(a) &=& \cC_1(a) \cup
\bigcup_{r=1}^{p(a)} \left(E^r(a)\big|_{N_1(a)}, 1 \right),\\
\cE_{p(a)}(a) &=& E(a) \backslash \bigcup_{r=1}^{p(a)} E^r(a).
\end{eqnarray*}

The definitions of $\inv(a)$ and an associated presentation
now proceed as in \S3.2. The resulting definition of $\inv(a)$
is independent of the choice of a presentation of $\io$ at $a$
(in particular, independent of its codimension $p(a)$ -- see
Theorem 2.3). We thus
obtain the desingularization principle Theorem 1.14 in the
more general setting.

\begin{remark}
We have assumed that $\io$ admits a semicoherent presenation
of codimension $p(a) \geq 1$ at each point $a$ in order to
make the above generalization of \S3.2 consistent with the
latter (and because we know of no interesting example where 
we need to use a presentation $\ucG(a)$ with $p(a) = 
\codim \ucG(a) = 0$ and $\mu_{\ucG(a)} > 1$).
The same local construction can, of course, be used
if we merely assume that $p(a) \geq 0$, but the codimension
of the presentation of $\inv_{r+1}$ will be shifted by $1$.
($\inv_{r+1}$ will have a presentation $\ucH_{r+1}(a)$ of
codimension $r$ at $a$, when $r > p(a)$.)
\end{remark}

\subsubsection{Independence of the embedding dimension}
Consider $X \subset M$ and $\inv = \inv_X$, where $\io_X(a)$
is the Hilbert-Samuel function $H_{X,a}$. If $X \subset M
\subset M'$, where $\dim M' > \dim M$, then $\inv_X$ as
defined above would not be the same for $X$ as a subspace
of $M$ or $M'$ (cf. examples of \cite{Eemb}). We can resolve
this problem by a simple variation of \S6.1.1:

The Hilbert-Samuel function $H_{X,a}$ determines the 
{\it minimal embedding dimension} $e_{X,a}$ of $X$ at
$a \in M$:
$$
e_{X,a} \ = \ H_{X,a}(1) - 1\ .
$$
Let $n = \dim_a M$ and $e(a) = e_{X,a}$. There is a 
semicoherent presentation $\ucC(a)=( N(a),\cC(a),\emptyset)$
of $\io_X(\cdot) = H_{X,\cdot}$ at $a$, where $N(a)$ lies in
a minimal embedding submanifold for
$X$ at $a$; in particular, $\ucC(a)$ has codimension $p(a) \geq
n - e(a)$. Consider a sequence of $\io_X$-admissible transformations
(1.1). (We use the notation above.) Let $a \in M_j$. Then 
$\inv_{1/2} = \io_X$ has a semicoherent presentation 
$\ucC_1(a) = ( N_1(a), \cC_1(a), \cE_1(a))$ at $a$, of codimension
$p(a) \geq n - e(a)$, where $e(a) = e_{X_j, a}$. 
 
We use the construction of \S6.1.1, but where {\it we now think of}
$\ucC_1(a)$ {\it as a presentation of codimension} $e(a) - (n - p(a))$
{\it in a minimal embedding space for} $X_j$ {\it at} $a$. So, if 
$p(a) > n - e(a)$ and $q(a)$ denotes $e(a)-(n-p(a))$, then 
$$
\inv_{q(a)}(a) \ = \ ( H_{X_j,a}, s_1(a); 1, s_2(a);
\ldots ; 1, s_{q(a)}(a) )\ ,
$$
with a semicoherent presentation 
$$
\ucH_{q(a)}(a) \ = \ ( N_{q(a)}(a), 
\cH_{q(a)}(a), \cE_{q(a)}(a) )
$$
at $a$, of codimension $p(a)$ in $M_j$, or of codimension
$q(a)$ in a minimal embedding submanifold, and we
now proceed as in \S3.2.

If $p(a) = n - e(a)$, then $\inv_1(a) = ( H_{X_j,a}, s_1(a) )$
has a presentation $\ucH_1(a) = ( N_1(a), \cH_1(a), \cE_1(a) )$
at $a$, where $\cH_1(a) = \cC_1(a) \cup 
\left( E^1(a)\big|_{N_1(a)}, 1 \right)$, of codimension $p(a)$
in $M_j$, or of codimension $0$ in a minimal embedding
submanifold. In this case, it is easy to see that 
$\mu_{\ucH(a)} = 1$, so we can find an equivalent presentation
in codimension $+1$, and proceed as usual. For example, $X_j$
is smooth at $a$ if and only if $\dim_a X_j = e(a)$; in this
case, $\inv_1(a) = (H_{e(a)}, s_1(a))$, where $H_e$ denotes the
Hilbert-Samuel function of a smooth space of dimension $e$:
$$
H_e(k) \ = \ {e+k \choose e}\ , \quad k \in \IN \ .
$$

The resulting desingularization invariant $\inv_X$ 
is independent of the dimension of a smooth embedding space
$M$ for $X$, and the desingularization principle Theorem
1.14 applies to give the universal embedded desingularization
Theorem 6.1. The embedded desingularization
algorithm stops over a neighbourhood of $a$ when 
$\inv_X(a) = (H_{e(a)}, 0; \infty)$.

\subsection{Comparison of weak and strict transforms}
Let $X$ denote a closed subspace of $M$ and let
$\cJ = \cI_X \subset \cO_M$. (See \S1.3.) The purpose of
this subsection is to study the strict transforms of $X$
by the sequence of blowings-up involved in principalization
of $\cJ$ according to the desingularization principle
Theorem 1.14. (See Example 1.19(1).) We show, in particular,
that the theorems of \cite{EVweak} (cf. Example 1.19(4))
and \cite{BV} are direct consequences of the desingularization
principle (Corollaries 6.5 and 6.7 below). Theorem 6.8
following generalizes these results by using the 
desingularization principle in a novel way.

Consider an $\inv_{\cJ}$-admissible sequence of transformations
(1.2), where $\inv_{1/2}(a) = \mu_a(\cJ_j)$, $a \in M_j$.
(See \S3.2.) The following lemma is the key point of this
subsection.

\begin{lemma}
Let $a \in M_j$. Then 
\begin{equation}
\inv_{\cJ}(a) \ = \ (1, 0; \ldots ; 1, 0; \infty)
\end{equation}
(where, let us say, there are $t$ pairs $(1, 0)$) if and
only if there are local coordinates $(x_1, \ldots, x_{n-t},
x_{n-t+1}, \ldots, x_n)$ for $M_j$ at $a = 0$ in which:
\begin{enumerate}
\item
$E(a) = \{ H \}$, where each $x_H = x_i$, for some 
$i = 1, \ldots, n-t$.
\item
Set $\tx = (x_1, \ldots, x_{n-t})$. Then $\cJ_{j,a}$ is 
generated by 
$$
x_n\ , \ \tx^{\theta_1}x_{n-1}\ , \ldots, \ \tx^{\theta_{t-1}}x_{n-t+1}\ ,
$$
where each $\tx^{\theta_k}$ denotes a monomial
$$
\tx^{\theta_k} \ = \ x_1^{\theta_{k1}} \cdots x_{n-t}^{\theta_{k,n-t}}\ ,
$$
and
\begin{enumerate}
\item
each $\tx^{\theta_k}$ is a monomial in the $x_H$, $H \in E(a)$;
i.e., $\theta_{ki} = 0$ unless $x_i = x_H$, for some $H$;
\item
$\theta_1 \leq \theta_2 \leq \cdots \leq \theta_{t-1}$ (where
$\leq$ means componentwise inequality).
\end{enumerate}
\end{enumerate}
\end{lemma}

\begin{proof}
The assertion can be seen by following the local construction in \S3.2.
In the language of the latter, $\tx^{\theta_1} = D_2(a)$ and
$$
\tx^{\theta_{k+1}-\theta_{k}} \ = \ D_{k+2}(a) , \quad
k \ = \ 1, \ldots , t-2 \ .
$$
(This is where (2)(b) comes from.) Moreover, $\{ x_n = \cdots =
x_{n-t+1} = 0 \}$ is a maximal contact subspace $N_t(a)$
(in the notation of \S3.2), and $\tx$ restricts to a coordinate
system on $N_t(a)$.
\end{proof}

Now set $X_0 = X$ and, for each $j$, let $X_{j+1}$ denote
the strict transform of $X$ by the blowing-up $\s_{j+1}$.
(We are not assuming that the blowings-up $\s_{j+1}$ are in
any sense ``admissible'' for the strict transforms.)

\begin{corollary}
If $a \in M_j$ and $\inv_{\cJ}(a) = (1, 0; \ldots ; 1, 0; \infty)$,
then $a \in X_j$.
\end{corollary}

\begin{proof}
Let $\pi_j : \ M_j \to M_0 = M$ denote the composite of the
blowings-up $\s_1, \ldots, \s_j$. By Lemma 6.3, 
$\pi_j \big|_{N_t(a) \backslash \bigcup \{ H \in E(a) \}}$ is an
isomorphism of $N_t(a) \backslash \bigcup \{ H \in E(a) \}$
with an open set of smooth points of $X$. The claim follows.
\end{proof}

\begin{corollary}[{\rm \cite{EVweak}; see Example 1.19(4)}]
Suppose that $X$ is pure-dimensional (of dimension $n-t$, say,
where $n = \dim M$). Let $\cJ = \cI_X$. Then there is a finite
sequence of $\inv_{\cJ}$-admissible blowings-up (1.2) with
smooth centres in the successive inverse images of $\Sing X$,
such that:
\begin{enumerate}
\item
The final strict transform $X' = X_{j_0}$ is smooth.
\item
$X'$ and $E' = E_{j_0}$ simultaneously have only normal
crossings.
\end{enumerate}
\end{corollary}

\begin{proof}
We simply apply the algorithm for principalization of 
$\cJ = \cI_X$, stopping when the maximum value of the
$\inv_{\cJ}(\cdot)$ becomes $(1,0;\ldots;1,0; \infty)$
(where there are $t$ pairs $(1,0)$).
\end{proof}

The precise statement of Lemma 6.3 immediately provides
stronger versions of this result. Corollary 6.7 below,
for example, is the theorem of \cite{BV}. 
The point is that,
assuming (6.1), we can use Lemma 6.3 to completely describe
the stratification by values of $\inv_{\cJ}$, in a neighbourhood
of $a$. We will need to work only with the largest value of
$\inv_{\cJ}$, apart from $\inv_{\cJ}(a)$. 

\begin{lemma}
Assume (6.1). Suppose that 
$$
\theta_{t-q-1} \ < \ \theta_{t-q} \ = \ \cdots \ = \ \theta_{t-1} \ ,
$$
where $q \leq t-1$ (and where $\theta_0$ means $0 \in \IN^{n-t}$;
we are using the notation of Lemma 6.3). Let $\tau(a)$ denote
the largest value of $\inv_{\cJ}$, apart from $\inv_{\cJ}(a)$,
in a small neighbourhood of $a$. Then 
\begin{equation}
\tau(a) \ = \ (1, 0; \ldots ; 1, 0; 0)\ ,
\end{equation}
where there are $t-q$ pairs $(1,0)$ in (6.2). The locus of
points $x$ near $a$ where $\inv_{\cJ}(x) = \tau(a)$ is 
given by
\begin{equation}
\begin{array}{l}
x_n \ = \ x_{n-1} \ = \ \cdots \ = \ x_{n-t+q+1} \ = \ 0 \ ,\\
x_{n-k} \ \neq \ 0 \ , \quad \mbox{ for some } 
k \ = \ t-q, \ldots, t-1 \ ,\\ 
D_{t-q+1}(a) \ = \ \tx^{\theta_{t-q} - \theta_{t-q-1}} \ = \ 0 \ .
\end{array}
\end{equation}
\end{lemma}

This is a simple consequence of Lemma 6.3.
Let $\Sig = \{ x: \inv_{\cJ}(x) \geq \tau(a) \}$ (i.e.,
the closure of the locus (6.3)). Then $\Sig$ is given by
\begin{equation*}
\begin{array}{l}
x_n \ = \ x_{n-1} \ = \ \cdots \ = \ x_{n-t+q+1} \ = \ 0 \ ,\\
D_{t-q+1}(a) \ = \ 0 \ .
\end{array}
\end{equation*} 
Consider the local blowing-up $\s$ with centre given by any
component of $\Sig$; i.e., with centre
\begin{equation*}
\begin{array}{l}
x_n \ = \ x_{n-1} \ = \ \cdots \ = \ x_{n-t+q+1} \ = \ 0 \ ,\\
x_i \ = \ 0 \ , \mbox{ for some } x_i \mbox{ occuring in }
D_{t-q+1}(a) \ .
\end{array}
\end{equation*}
In the ``$x_i$-chart'' $U_{x_i}$ of this blowing-up, the
weak transform $\cJ_j'$ of $\cJ_j$ is generated by
$$
x_n \ , \ \tx^{\theta_1}x_{n-1} \ ,\ldots, 
\ \tx^{\theta_{t-q-1}}x_{n-t+q+1} \ , 
\ \tx^{\theta_{t-q}-(i)}x_{n-t+q} \ ,\ldots,
\ \tx^{\theta_{t-q}-(i)}x_{n-t} \ ,
$$
where $(i)$ denotes the multiindex of length $n-t$ with
$1$ in the $i$'th place and $0$ elsewhere.
The strict transform $X_j'$ of $X_j$ is still given by
$$
x_n \ = \ x_{n-1} \ = \ \cdots \ = \ x_{n-t+1} \ = \ 0 \ 
$$
in the chart $U_{x_i}$.

\begin{corollary}[{\rm \cite{BV}}]
Suppose that $X$ is a closed subspace of $M$. Let
$\cJ = \cI_X$. Then there is a finite
sequence of $\inv_{\cJ}$-admissible blowings-up (1.2) with
smooth centres in the successive inverse images of $\Sing X$,
such that:
\begin{enumerate}
\item
The final strict transform $X' = X_{j_1}$ is smooth.
\item
$X'$ and $E' = E_{j_1}$ simultaneously have only normal
crossings.
\item
The final total transform $\cJ_{\pi_{j_1}^{-1}(X)}$ is
the product of $\cI_{X_{j_1}}$ with a normal crossings
divisor supported on the exceptional locus. ($\pi_{j_1}$
denotes the composite of the sequence of blowings-up.)
\end{enumerate}
\end{corollary}

\begin{proof}
Let us first suppose that $X$ is pure-dimensional (of
dimension $n-t$, say, where $n = \dim M$. We apply the
algorithm for principalization of $\cJ$ as in Corollary
6.5, stopping (in year $j_0$, say) when the maximum value
of $\inv_{\cJ}(\cdot)$ becomes $(1, 0; \ldots; 1, 0; \infty)$
(where there are $n$ pairs $(1, 0)$).

We can now simply continue to blow up with centre determined by
the maximum value of the extended invariant $\inv_{\cJ}^e$
(cf. proof of Corollary 1.17) {\it outside} the strict
transform of $X$. (For example, if the maximum stratum in
$M_{j_0} \backslash X_{j_0}$ is not closed in $M_{j_0}$, then
its closure is smooth and has only normal crossings with
respect to $X_{j_0}$, by Lemma 6.6.) We continue to blow up
with centre given by the closure of the locus of maximum values of
$\inv_{\cJ}^e$ outside the strict transform of $X$, until the
support of the weak transform equals the support of the strict
transform, say in year $j_1$; i.e., until at any point $a \in X_{j_1}$,
$\cJ_{j_1,a}$ is generated by
$$
x_n, \ x_{n-1}\ , \ldots, \ x_{n-t+1}
$$
(in suitable coordinates as in Lemma 6.3). In particular,
$\cJ_{j_1} = \cI_{X_{j_1}}$ and the assertion follows.

We can of course continue in this way to prove Corollary 6.7
(and also Corollary 6.5) without the assumption of 
pure-dimensionality. Suppose that $X$ has smooth parts of
codimensions $t_1 , t_2, \ldots, t_q$. We blow up first
until the maximum value of the invariant is $(1, 0; \ldots; 1, 0; 
\infty)$ (with $t_1$ pairs $(1, 0)$), say on $Z(t_1)$. We
then continue as above until the maximum value of the invariant
{\it outside} the strict transform of $Z(t_1)$ is 
$(1, 0; \ldots; 1, 0; \infty)$ (with $t_2$ pairs $(1, 0)$),
say on $Z(t_2)$. We now continue using the maximum value of
the (extended) invariant on the complement of the strict
transforms of $Z(t_1)$ and $Z(t_2)$, etc.
\end{proof}

The assertion of Corollary 6.7 is not restricted to the
weak form of desingularization given by Corollary 6.5.
Beginning with the conclusion of the embedded desingularization
theorem (Example 1.19(1) above),
we can transform to the product condition (3) (in Corollary 6.7)
using the following theorem (applied with $X$ $=$ the final strict
transform of our original closed subspace ($Y$, say) and with
$\cJ$ $=$ the final weak transform of $\cI_Y$).

\begin{theorem}
Let $X$ denote a smooth closed subspace of $M$, and let
$E$ be a collection of smooth hypersurfaces in $M$ such
that $X$ and $E$ simultaneously have only normal crossings.
Let $\cJ$ be an ideal of finite type in $\cO_M$ such that
$\supp \cO_M / \cJ \subset X \cup E$ and $\cJ = \cI_X$
on $X \backslash E$. Then there exists a finite sequence of
transformations
\begin{equation}
\begin{array}{rccccccccl}
\longrightarrow & M_{j+1} & \stackrel{\s_{j+1}}{\longrightarrow}
& M_j & \longrightarrow & \cdots & \longrightarrow 
& M_0 & = & M \\
& X_{j+1} & & X_j & & & & X_0 & = & X\\
& \cJ_{j+1} & & \cJ_j & & & & \cJ_0 & = & \cJ\\
& E_{j+1} & & E_j & & & & E_0 & = &E 
\end{array}
\end{equation}
where, for each $j$,
\begin{list}{}{}
\item[$\s_{j+1}$] is a blowing-up of $M_j$ with smooth centre
$C_j$ that is $\inv_{\cJ}$-admissible, is supported in $E_j$,
and simultaneously has only normal crossings with respect to
$X_j$ and $E_j$,
\item[$X_{j+1}$] denotes the strict transform of $X_j$,
\item[$\cJ_{j+1}$] denotes the weak transform of $\cJ_j$,
\end{list}
\noindent
and such that the final transforms $X'$ and $\cJ'$ satisfy
$$
\cJ' \ = \ \cI_{X'} \ ;
$$
thus, the final total transform $\pi^{-1}(\cI_{X})$ is the
product of $\cI_{X'}$ with a normal crossings divisor
supported in $E'$.
\end{theorem}

\begin{proof}
We define a new invariant $\inv_{\cJ,X}(\cdot)$ inductively
over a sequence (6.4), in the following way. Let $a \in M_j$.
Set $\inv_{1/2}(a) = \io(a)$, where
$$
\io(a) \ := \ (\inv_{\cJ}(a), \, \de_{X_j}(a)) \ ,
$$
and
\[
\de_{X_j}(a) \ = \ \left\{ \begin{array}{r@{\ , \quad}l}
1 & a \in X_j \\ 0 & a \notin X_j \ .
\end{array} \right.
\]
The $\inv_{1/2} = \io$ satisfies Hypotheses 1.10(1) and (2), as
well as (a slight variant of) (3): If $\ucH(a) = ( N(a),
\cH(a), \cE(a) )$ is a semicoherent presentation of 
$\inv_{\cJ}$ at $a$ (the term $\cE(a)$ is irrelevant here),
then $\inv_{1/2}$ has a semicoherent 
presentation at $a$ given by
$$
( N(a), \ \cH(a) \cup \{ (g_k, 1) \}, \ \cE_1(a) )\ ,
$$
where $\{ g_k \}$ is the set of restrictions to $N(a)$
of a finite set of generators of $\cI_{X_j,a}$ with
linearly independent gradients (and $\cE_1(a) =
E(a) \backslash E^1(a)$, where $E^1(a)$ is determined
using $\inv_{1/2}$ according to Definitions 1.15, as usual).
(Although this semicoherent presentation of $\inv_{1/2}$
is in the weaker sense of Remark 3.5, it suffices to begin
the inductive construction as in \S3.2 because its
semicoherent equivalence class depends only on $\cJ_{j,a}$,
$X_{j,a}$ and the various blocks of exceptional divisors 
involved in defining $\inv_{\cJ}(a)$.)

Then $\inv_{1/2}$ extends to an invariant $\inv_{\cJ,X}$
defined inductively over a sequence of blowings-up (6.4) with
$\inv_{\cJ,X}$-admissible centres, according to the desingularization
principle Theorem 1.14. If $a \in X_j \backslash E_j$, then
$$
\inv_{\cJ}(a) \ = \ (1, 0; \ldots; 1, 0; \infty)
$$
($(1, 0)$ occuring $t$ times, where $t = \codim_a X_j$), and
$$
\inv_{\cJ,X}(a) \ = \ ( (\inv_{\cJ}(a), 1), 0; \infty ) \ .
$$
On the other hand, if $a \in X_j$ and 
$$
\inv_{\cJ,X}(a) \ = \ (((1,0;\ldots;1,0;\infty), 1), 0 ; \infty) \ ,
$$
then $\inv_{\cJ}(a) \ = \ (1, 0; \ldots; 1, 0; \infty)$, so that
Lemma 6.3 applies. Hence we can blow up the maximum locus of
$\inv_{\cJ,X}^e$ until $\inv_{\cJ,X}$ is constant on $X_j$
(at least when $X$ is pure-dimensional), and then argue as in
the proof of Corollary 6.7; the general (not necessarily 
pure-dimensional) case follows as in the latter.
\end{proof}

\subsection{Simultaneous desingularization of parametrized
families}
A {\it family of spaces parametrized by a space} $T$ means
a morphism $p: \ X \to T$ (e.g. a morphism of schemes of
finite type over a field of characteristic zero, or a morphism
of analytic spaces). The fibres $X_t$, $t \in T$, form
a family of closed subspaces of $X$.

\begin{definitions} Let $p: M \to T$ be a morphism of smooth
spaces (manifolds). We say that $p$ is {\it smooth at} $a \in M$
if there is a system of local coordinates $(x_1, \ldots, x_n)$
in a neighbourhood of $a$, in which $p$ is a projection
onto a coordinate subspace $(x_1, \ldots, x_n) \mapsto
(x_1, \ldots, x_k)$ (where $k \leq n$). We say that $p$ is 
{\it smooth} if it is smooth at every point of $M$.

Let $E$ denote a configuration (finite collection) of smooth
closed subspaces of $M$ that simultaneously have only normal
crossings. We say that the morphism $p$ restricts to a {\it
smooth projection} from $E$ to $T$ if $p$ restricts to a smooth
morphism of every element of $E$ and of every intersection of
elements of $E$. (There is a more general notion of ``smooth
morphism'' between spaces that are not necessarily smooth
\cite[Chapt. III, Sect. 10]{Hart} which we will not need;
it is equivalent to the preceding when the target space is
smooth and the source has only normal crossings.)
\end{definitions}

\begin{definitions} {\it Simultaneous resolution of singularities
of a family of spaces} $X \to T$ means (some version of)
resolution of singularities of $X$ by blowings-up with smooth
centres such that all centres and also the final strict 
transform of $X$ project smoothly to $T$.

{\it Simultaneous embedded desingularization} means simultaneous
resolution of singularities with the additional property that,
at each step, the entire configuration of exceptional divisors
together with the centre of blowing up (or together with the
final strict transform of $X$) projects smoothly to $T$.
\end{definitions}

A simultaneous resolution of singularities (respectively,
simultaneous embedded resolution of singularities) of
$X \to T$ restricts
to a resolution of singularities (respectively, embedded
resolution of singularities) of every fibre $X_t$.

Any theorem of resolution of singularities by blowings-up
with smooth centres immediately implies that the parameter-
space of a proper family can be stratified so that the
fibres can be simultaneously desingularized over every stratum
(cf. \cite[Thm. 4, Ref. BM99]{PS}, \cite[Sect. 4]{ENV}):

\begin{theorem}
Let $p: \ X \to T$ denote a proper morphism. Then there
is a finite filtration by closed subsets,
$$
T \ = \ T_0 \ \supset \ T_1 \ \supset \ \cdots \ \supset T_l
\ = \ \emptyset \ ,
$$
such that, for all $k$, $U_k := T_k \backslash T_{k+1}$ is smooth
and the family $X_k \big|_{U_k} \to U_k$ admits simultaneous
(embedded) desingularization.
\end{theorem}

\begin{proof}
We can assume the $p$ is surjective.
Consider (embedded) resolution of singularities of $X$.
By the ``generic smoothness theorem'' (cf. \cite[Cor. 10.7]{Hart}),
there is a proper closed subspace $T_1$ of $T$ such that 
$T \backslash T_1$ is smooth and (over $T \backslash T_1$)
all centres of blowing up, as
well as the final strict transform of $X$, project smoothly
onto $T \backslash T_1$ (and, at every step, the collection 
of exceptional divisors together with the centre, or together
with the final strict transform of $X$, projects smoothly onto
$T \backslash T_1$). Thus the family of fibres admits
simultaneous (embedded) desingularization over $T \backslash T_1$.
The result follows by induction.
\end{proof}

\begin{remark}
We have not stated this theorem with the precisions
that are evidently needed to cover all categories. For
example, in the complex-analytic category, we should either 
apply the statement to a relatively compact open
subset of $T$, or use a locally finite
filtration; in the real-analytic case, we should either
assume that $p$ admits a proper complexification, or use
a semianalytic filtration of $T$.
\end{remark}

Our main purpose in this section is to give a more precise
version of Theorem 6.11 by extending Theorem 1.14 to a
``desingularization principle for families'' (Theorem 6.15
below). The Hilbert-Samuel function of $X$, as well as the
order of the ideal $\cI_X$ satisfy the Hypotheses 6.13 below
that are needed for Theorem 6.15. (See Lemma 6.14 and
Theorem 6.17). The theorem of \cite{ENV} uses
the order of $\cI_X$, so Theorems 6.15 and 6.17 extend the latter
to the stronger form of embedded desingularization
(cf. Examples 1.19).
\medskip

\noindent {\it Notation}.
Let $X \to T$ be a family of spaces, as above. Let $t \in T$.
We let $X_t$ denote the fibre over $t$, and set
$$
X_{(t)} \ := \ X_t \times (\mbox{germ of } T \mbox{ at } t) \ .
$$
(It is more convenient to use $X_{(t)}$ rather than $X_t$ in
comparing the Hilbert-Samuel functions (or other invariants)
of $X$ and of $X_t$ at a given point in $X_t$.) If $a \in X$,
let $t(a)$ denote the image of $a$ in $T$.
\medskip

If $X \to T$, then, locally in $X$, there is an embedding
$X \hookrightarrow M$ in a manifold $M$, together with a
smooth morphism $M \to T$ that restricts to the given 
projection $X \to T$. In Hypotheses 6.13 and Theorem 6.15
following, we will therefore assume that $X \hookrightarrow M$
and that $X \to T$ is induced by a smooth projection $p:\ M \to T$.
As in \S6.1, however, our results will be independent of the
embedding.

Let $\ucH(a) = (N(a), \cH(a), \cE(a))$ denote a local presentation
of codimension $q$ at $a \in M$. If $\cE(a)$ and $N(a)$ together
map smoothly to $T$, then we can define the {\it restriction of}
$\ucH(a)$ {\it to the fibre} $M_{t(a)}$ in an obvious way: Let
$t = t(a)$ and let $\ucH(a)_t := (N(a)_t, \cH(a)_t, \cE(a)_t)$, 
where $N(a)_t$ is the fibre of $N(a) \to T$, $\cE(a)_t :=
\{H_t: H \in \cE(a)\}$, where $H_t$ denotes the fibre of $H \to T$,
and 
$$
\cH(a)_t \ = \ \{ (h_t, \mu_h):\ (h, \mu_h) \in \cH(a) \} \ ,
$$
where $h_t := h \big|_{N(a)_t}$.

\begin{hypotheses}
We will assume that $\io_X$ satisfies Hypotheses 1.10 (1) and (2),
together with the following additional property:
\begin{enumerate}
\item[(3)] 
Let $a \in X$ and let $t = t(a)$. Then:
\begin{enumerate}
\item 
$\io_{X_{(t)}}(a) \geq \io_X(a)$ .
\item 
If $\io_{X_{(t)}}(a) = \io_X(a)$, then there is
a semicoherent presentation 
$$
\ucH(a) \ = \ (N(a), \cH(a), \emptyset)
$$
of $\io_X$ at $a$, of codimension $q = q(a)$, say, such that
$N(a) \to T$ is smooth, $\ucH(a)_t$ is a semicoherent
presentation (of codimension $q$) of $\io_{X_{(t)}}$ at $a$,
and $\mu_{\ucH(a)_t} = \mu_{\ucH(a)}$.
\item
If there is a (germ of a) smooth subspace $S \subset S_{\io_X}$
such that $S \to T$ is smooth, then $\io_{X_{(t)}}(a) = \io_X(a)$.
\end{enumerate}
\end{enumerate}
\end{hypotheses}

\begin{lemma}
let $\io_X(a)$ denote the order 
$\mu_a(\cJ)$ of
$\cJ := \cI_X$ (cf. Definitions 1.4 and Examples 1.8). Then
$\io_X$ satisfies Hypotheses 6.13.
\end{lemma}

\begin{proof}
By Example 1.13(1), $\io_X$ satisfies Hypotheses 1.10; 
we have to verify 6.13(3).
Choose coordinates $(x_1, \ldots, x_n)$ for $M$ in a neighbourhood
of $a = 0$, in which $p$ is a projection $(x_1, \ldots, x_n)
\mapsto (x_1, \ldots, x_{n-r})$; write $u := (x_1, \ldots, x_{n-r}),\ 
v := (x_{n-r+1}, \ldots, x_n)$. If $\{g_i(u,v):\ i=1,\ldots,q\}$ 
is a set of
generators of $\cJ_a$, then $\cJ_{(t),a}$ is generated by
$\{g_i(0,v)\}$. Therefore, properties (a) and (b) are obvious. 

To prove (c): Consider a smooth germ $S \subset 
S_{\mu_{\cdot}(\cJ)}(a)$
such that $S \to T$ is smooth. Then $\mu_{S,a}(\cJ) = \mu_a(\cJ)$
(because $\mu_{\cdot}(\cJ)$ is constant on $S$; cf. Definitions
1.4). Let $d = \mu_a(\cJ)$ and write
$$
g_i(u,v) \ = \ \sum_{\al \in \IN^r,\ |\al| \geq d} g_{i,\al}(u)v^\al,
\quad i=1,\ldots,q \ .
$$
Then $g_{i,\al}(0) \neq 0$, for some $i$ and some $\al$ such
that $|\al| = d$. Therefore, $\mu_a(\cJ_{(t)}) = d$.
\end{proof}

\begin{theorem}
Suppose that $\io_X$ satisfies Hypotheses 6.13.
Consider any sequence of $\inv_X$-admissible local blowings-up
(1.1). Let $a \in X_j$ and let $a_i$ denote the image of $a$
in $X_i$, $i \leq j$. Assume that, for all $i < j$, $E_i$
together with $C_i$ maps smoothly to $T$ at $a_i$, and
$\inv_{X_{(t)}}(a_i) = \inv_X(a_i)$. (In particular, the
sequence up to year $j$ induces an $\inv_{X_{(t)}}$-admissible
sequence (locally at the points $a_i$) on the fibres over $t$.)
Then
$$
\inv_{X_{(t)}}(a) \ \geq \ \inv_X(a) \ ,
$$ 
and the following conditions are equivalent:
\begin{enumerate}
\item
$\inv_{X_{(t)}}(a) = \inv_X(a)$.
\item
There is a smooth subspace (germ) $S$ of $S_{\inv_X}(a)$ such
that $S \to T$ is smooth.
\item
$S_{\inv_X}(a)$ maps smoothly to $T$.
\end{enumerate}
\end{theorem}

By $\inv_{X_{(t)}}$, we mean the invariant for the fibre $X_t$
given by the desingularization principle Theorem 1.14, beginning
with $\inv_{1/2} = \io_{X_{(t)}}$.
The assumption that $\inv_{X_{(t)}}(a_i) = \inv_X(a_i)$, $i < j$,
is, in fact, redundant -- it follows by induction from the
conclusion of the theorem.

\begin{proof}[Proof of Theorem 6.15]
We will follow the local inductive construction of \S3.2
(or, more generally, \S6.1), {\it but we will use Villamayor's
variant of this construction} (as described
in \S3.4), to extend the properties of
Hypotheses 6.13 to the truncated invariants $\inv_{r - 1/2}(a)$
and $\inv_r(a)$, for each successive $r$. (See Remark 6.16.)

Note first that, by the hypotheses, $E(a) \to T$ is smooth and,
if $\io_{X_{(t)}}(a) = \io_X(a)$, then there is a semicoherent
presentation 
$$
\ucH(a) = (N(a), \cH(a), \cE(a))
$$
of $\io_X$ at $a$, satisfying the obvious generalization of
property 6.13(3)(b): $\ucH(a)_t = (N(a)_t, \cH(a)_t,
\cE(a)_t)$ is a presentation of $\inv_{1/2} = \io_{X_{(t)}}$
at $a$ (for the fibre $X_{j,t}$).

Now consider $r \geq 1$, and assume that $\inv_{r - 1/2}$
satisfies the analogues of Hypotheses
1.10(1) and 6.13(3) (where, in (3)(b), there is a semicoherent
presentation 
$$
\ucH_r(a) = (N_r(a), \cH_r(a), \cE_r(a))
$$
of $\inv_{r - 1/2}$ at $a$, that restricts to a semicoherent
presentation 
$$
\ucH_r(a)_t = (N_r(a)_t, \cH_r(a)_t, \cE_r(a)_t)
$$
of $\inv_{r - 1/2,\, t}$ (the invariant for the fibre) at $a$).
We then have to show that $\inv_r$ satisfies the analogues
of 6.13(3)(a)-(c). (For this, it is enough to assume that the
previous centres of blowing up are $(r - 1/2)$-admissible
and map smoothly to $T$.)

First consider (a) and (b). If $\inv_{r - 1/2,\, t}(a) >
\inv_{r - 1/2}(a)$, then there is nothing more to do. 
Assume $\inv_{r - 1/2,\, t}(a) = \inv_{r - 1/2}(a)$. Then
automatically $E^r_t(a) = E^r(a)$ (more precisely, $E^r_t(a)
= \{H\big|_{X_t}: \ H \in E^r(a)\}$, where $E^r_t(a)$ denotes
the analogue of $E^r(a)$ for the fibre $X_t$), and $s_{r,t}(a)
= s_r(a)$ (where $s_{r,t}(a)$ is again the analogue of $s_r(a)$
for the fibre). (In fact, these equalities hold at $a_i$, for
all $i \leq j$.) It follows that $\inv_{r,t}(a) = \inv_r(a)$,
and there is a semicoherent
presentation
$$
\ucC_r(a) = (N_r(a), \cC_r(a), \cE_r(a))
$$
of $\inv_r$ at $a$, with the properties required for (b).

Property (c) is obvious because, if $S \subset S_{\inv_r}(a)$
is a smooth germ such that $S \to T$ is smooth, then
$\inv_{r - 1/2,\, t}(a) = \inv_{r - 1/2}(a)$, by property (c) for
$\inv_{r - 1/2}$, so that $s_{r,t}(a) = s_r(a)$, as above.
This completes the step from $\inv_{r - 1/2}$ to $\inv_r$.

Now assume that $\inv_r$ satisfies the analogues of Hypotheses
1.10(1) and 6.13(3); in particular, in (b), there is a semicoherent
presentation $\ucC_r(a)$ of $\inv_r$ at $a$, as above, that
restricts to a semicoherent presentation of $\inv_{r,t}$ at $a$.
It is enough to assume that the previous centres of blowing up
are $\inv_{r + 1/2}$-admissible and map smoothly to $T$.

If $\inv_{r,t}(a) > \inv_r(a)$, then again there is nothing to
do. Assume that $\inv_{r,t}(a) = \inv_r(a)$. Then
$$
\mu_{r+1,\, t}(a) \ \geq \ \mu_{r+1}(a)\ ,
$$
by property (b) for $\inv_r$, because $\mu_{r+1}(a)$ and
$\mu_{r+1,\, t}(a)$ are realized as multiplicities of a given
function, before and after restriction to the fibre $X_t$
(cf. Lemma 6.14). On the other hand, for every $H \in \cE_r(a)$,
\begin{equation}
\mu_{r+1,H,\, t}(a) \ = \ \mu_{r+1,H}(a)\ ,
\end{equation}
by induction over the sequence of blowings-up, using Lemma 3.7
and the assumption on the previous centres. Therefore,
\begin{equation}
\nu_{r+1,\, t}(a) \ \geq \ \nu_{r+1}(a)\ ,
\end{equation}
so we have proved (a).

If $\inv_{r,t}(a) = \inv_r(a)$ and $\nu_{r+1,\, t}(a) =
\nu_{r+1}(a)$, then the presentation $\ucG_{r+1}(a)$
of $\inv_{r+1/2}$ at $a$ constructed as in \S\S3.2, 3.4, restricts
to a presentation $\ucG_{r+1, t}(a)$ of $\inv_{r+1/2,\, t}(a)$.
Thus (b) is satisfied.

To prove (c), let $S \subset S_{\inv_{r+1/2}}(a)$ be a smooth
germ such that $S \to T$ is smooth. By property (c) for 
$\inv_r$, $\inv_{r,t}(a) = \inv_r(a)$. Then, by (6.5), in order
to prove that $\nu_{r+1,\, t}(a) = \nu_{r+1}(a)$, it is enough
to prove that $\mu_{r+1,\, t}(a) = \mu_{r+1}(a)$. The latter is a
statement about the order of a function (or an ideal), established 
by Lemma 6.14. This completes the step from $\inv_r$ to 
$\inv_{r+1/2}$.

Finally, the full invariant $\inv_X$ satisfies the analogue
of property 6.13(3), and the conclusion of the theorem follows
immediately.
\end{proof}

\begin{remark}
Recall, from \S\S3.2, 3.3, that, in order to define
$\nu_{r+1}$ and prove the semicontinuity properties of
$\inv_{r + 1/2}$, it is enough to assume that the previous
centres of blowing up are $(r - 1/2)$-admissible.
We need the stronger $(r + 1/2)$-admissibility assumption
on the previous centres in the step from $\inv_r$ to $\inv_{r+1/2}$
in the proof of Theorem 6.15 in order to prove the
semicontinuity condition (6.6) using (6.5). This is the reason
for using Villamayor's smaller block of exceptional divisors
$\cE^r(a)$ -- Lemma 3.7 shows that, for $H$ in the smaller block,
$\mu_{r+1,H}(a)$ can be computed using the values in previous
years. It would seem that (6.5) need not hold for the additional
exceptional divisors that are factored out to define the
residual multiplicities according to the inductive construction
in the Bierstone-Milman algorithm.
\end{remark}

\begin{corollary}
Suppose that $\io_X$ satisfies Hypotheses 6.13. Consider
a desingularization of $X$ by a finite sequence of $\inv_X$-admissible
blowings-up (1.1). Then the following conditions are equivalent:
\begin{enumerate}
\item
All centres of blowing up $C_i$ (as well as the final strict
transform of $X$) project smoothly to $T$.
\item
$\inv_{X_{(t(a))}}(a) = \inv_X(a)$ for all $a \in M_i$,
$i = 0,1,\ldots\ $.
\end{enumerate}
\end{corollary}

This is an immediate consequence of Theorem 6.15.
Corollary 6.17 provides a more precise version of Theorem 6.11
(by using the generic smoothness theorem as in the proof of the
latter).

\begin{theorem}
Let $\io_X(a)$ denote the Hilbert-Samuel function
$H_{X,a}$ (cf. Examples 1.8). Then $\io_X$ satisfies
Hypotheses 6.13.
\end{theorem}

\begin{proof}
By Examples 1.13(2), $\io_X$ satisfies Hypotheses 1.10.
We have to verify properties (3)(a)-(c) of Hypotheses 6.13.
Let $a \in X$ and let $t = t(a)$.
\medskip

\noindent
(a) $H_{X_{(t)},a} \geq H_{X,a}$, by an elementary lemma
\cite[Lemma 7.5]{BMjams}.
\medskip

We will establish properties (b) and (c) (in fact, we will
first prove (c)) using a semicoherent presentation of the
Hilbert-Samuel function constructed in 
\cite[Theorems 9.4, 9.6]{BMinv} (so an understanding of the
following argument requires some familiarity with the
latter). We will use the notation of
\cite[(7.1)]{BMinv}; $\frN(\cdot)$ will denote the {\it
diagram of initial exponents} of an ideal in a ring of
formal power series, as defined, for example, in 
\cite[Sect. 3]{BMinv}. Let $\IK$ denote the underlying
ground field of our space. (In the following arguments in
the case of schemes, $\IK$ should really be understood as
the residue field of the point $a$; see \cite[Remark 3.8]{BMinv}
for an explanation of this point.)
\medskip

\noindent
(c) Let $S \subset S_{H_{X,\cdot}}(a)$ denote a germ of a
smooth subset such that $S \to T$ is smooth. 
Consider a presentation $\ucH(a) = (N(a), \cH(a), \emptyset)$
of $H_{X,\cdot}$ at $a$, satisfying the hypotheses of
\cite[Theorem 9.4]{BMinv}. Then $N(a) \to T$ is smooth, since
$S \subset N(a)$ and $S \to T$ is smooth.
\medskip

{\it Claim.} If $N(a) \to T$ is smooth, then $H_{X_{(t)},a}
= H_{X,a}.$
\medskip

To prove this claim: We use the notation of 
\cite[Theorem 9.4]{BMinv}. In particular, $w = (w_1, \ldots, w_{n-r})$
represents a coordinate system on $N(a)$, $\wcI_{X,a} \subset
\wcO_{M,a} = \IK\lbr W,Z \rbr$, where $W = (W_1, \ldots, W_{n-r})$
and $Z = (Z_1, \ldots, Z_r)$, and the vertices of $\frN (\wcI_{X,a})
\subset \IN^n$ depend only on the Z-coordinates; i.e.,
$\frN (\wcI_{X,a}) = \IN^{n-r} \times \frN^*$, where $\frN^*
\subset \IN^r$. Moreover, since $N(a) \to T$ is smooth, we
can assume that the coordinates $w$ split as $w = (u,v)$,
where $(u,v) \mapsto u$ is the projection $N(a) \to T$,
and that $\wcI_{X,a} \subset \IK\lbr W,Z \rbr = \IK\lbr U,V,Z\rbr$.
Since $\wcI_{X_{(t)},a}$ is obtained from $\wcI_{X,a}$
by setting $U = 0$,
$$
\frN (\wcI_{X_{(t)},a}) \ = \ \frN (\wcI_{X,a}) \ .
$$
Then, by \cite[Lemma 7.5]{BMinv},
$$
H_{X_{(t)},a} \ = \ H_{X,a} \ .
$$
\medskip

\noindent
(b) Assume that $H_{X_{(t)},a} = H_{X,a}$. Choose local
coordinates $(u,v) = (u_1, \ldots, u_{n-s}, v_1, \ldots, v_s)$
for $M$ at $a$ such that $(u,v) \mapsto u$ is the projection
to $T$. Set $J_1 := \wcI_{X,a} \subset \IK \lbr u,v \rbr$ and
$J_0 := \wcI_{X_{(t)},a}$; thus $J_0$ is the ideal in 
$\IK \lbr u,v \rbr$ generated by the evaluations at $u = 0$
of the elements of $J_1$. Write $H_{J_1}$ and $H_{J_0}$ for
the Hilbert-Samuel functions of the quotients of $\IK \lbr u,v \rbr$
by $J_1$ and $J_0$ (respectively); i.e., for $H_{X,a}$ and
$H_{X_{(t)},a}$ (respectively). Let $\um$ denote the maximal ideal
of $\IK \lbr u,v \rbr$. 

We will first show that
\begin{equation}
\frN (\wcI_{X_{(t)},a}) \ = \ \frN (\wcI_{X,a}) \ .
\end{equation}

Let $\frN^* \subset \IN^s$ denote the diagram of initial
exponents of $\wcI_{X_t,a} \subset \IK\lbr v \rbr$, so that
the diagram $\frN = \frN(J_0)$ is 
$$
\frN \ = \ \IN^{n-s} \times \frN^* \ \subset \IN^n \ .
$$
Then 
\begin{equation}
\IK \lbr u,v \rbr \ = \ J_0 \oplus \IK \lbr u,v \rbr^{\frN} \ ,
\end{equation}
where $\IK \lbr u,v \rbr^{\frN}$ denotes the subset of
formal power series
supported in the complement of $\frN$, by Hironaka's formal
division theorem \cite[Theorem 3.17]{BMinv}, and
\begin{equation}
\IK \lbr u,v \rbr \ = \ J_1 + \IK \lbr u,v \rbr^{\frN} 
\end{equation}
(where the sum is not necessarily direct),
by the formal division theorem with parameters. (See 
\cite[Theorem 3.1]{BMpisa}, \cite{galligo}.) By the assumption,
\begin{equation}
H_{J_1} \ = \ H_{J_0} \ = \ H_\frN \ ,
\end{equation}
where 
$$
H_\frN (k) \ := \ \#\{\al \in \IN^n: \ \al \notin \frN, 
\ |\al| \leq k\}, \quad k \in \IN \ .
$$
\medskip

{\it Claim.} $J_1 \cap \IK \lbr u,v \rbr^{\frN} = \{0\}$.
Moreover, if $f = g + h$, where $f \in \um^k$, $g \in J_1$
and $h \in \IK \lbr u,v \rbr^{\frN}$, then $g \in \um^k$ 
and $h \in \um^k$.
\medskip

The claim implies (6.7); i.e., $\frN(J_1) = \frN(J_0)$:
Consider any vertex $\al$ of $\frN^*$. Using (6.9) and the
claim, we can write
\begin{equation}
v^\al \ = \ f_\al (u,v) + r_\al (u,v) \ ,
\end{equation}
where $f_\al (u,v) \in J_1$, $r_\al (u,v) \in 
\IK \lbr u,v \rbr^{\frN}$, and the order of any monomial
in $r_\al (u,v)$ is at least $|\al|$. Set $u = 0$ in (6.11).
Then
$$
v^\al \ = \ f_\al (0,v) + r_\al (0,v) \ .
$$
Thus, $f_\al (0,v)$ is the element of the {\it standard basis}
of $J_0$ representing the vertex $\al$; in particular, every
monomial in $r_\al (0,v)$ has exponent $> (0,\al)$, 
according to the formal division theorem (where the elements
$\delta \in \IN^n$ are ordered according to the lexicographic
order of $(|\delta|, \delta)$. See 
\cite[Corollary 3.19]{BMinv}.) It follows that any monomial
in $r_\al (u,v)$ has exponent $> (0,\al)$: Consider a monomial
$u^\beta v^\gamma$ in $r_\al (u,v)$, where $|\beta| + |\gamma|
= |\al|$. If $\beta \neq 0$, then $(\beta, \gamma) > (0, \al)$,
by the definition of the ordering. On the other hand, if
$\beta = 0$, then $(0, \gamma) > (0, \al)$, since all exponents
of $r_\al (0,v)$ are $> (0,\al)$. Therefore, the initial
exponent (the smallest exponent) of $f_\al (u,v)$ is $(0, \al)$.
Hence
$$
\frN(J_1) \ \supset \ \frN(J_0) \ = \ \frN \ ,
$$
and it follows from (6.10) that $\frN(J_1) = \frN(J_0)$.
\medskip

Proof of the preceding claim. For each $k \in \IN$, there is
a surjective homomorphism
\begin{equation}
\frac{\IK \lbr u,v \rbr^{\frN}}{\IK \lbr u,v \rbr^{\frN} \cap
(J_1 + \um^{k+1})} \ \to 
\ \frac{\IK \lbr u,v \rbr}{J_1 + \um^{k+1}} \ .
\end{equation}
Therefore,
\begin{eqnarray*}
H_\frN(k) & = & \dim \frac{\IK \lbr u,v \rbr^{\frN}}
		{\IK \lbr u,v \rbr^{\frN} \cap \um^{k+1}}\\
          & \geq & \dim \frac{\IK \lbr u,v \rbr^{\frN}}
		   {\IK \lbr u,v \rbr^{\frN} \cap
		   (J_1 + \um^{k+1})}\\
          & \geq & \dim \frac{\IK \lbr u,v \rbr}{J_1 + \um^{k+1}}\\
	  & = & H_{J_1}(k) \ = \ H_\frN(k) \ ,
\end{eqnarray*}
so all terms are equal.	Hence, for every $k$, (6.12) is an 
isomorphism, and 
\begin{equation}
\IK \lbr u,v \rbr^{\frN} \cap (J_1 + \um^{k+1}) \ = 
\ \IK \lbr u,v \rbr^{\frN} \cap \um^{k+1} \ .
\end{equation}
So
$$
J_1 \cap \IK \lbr u,v \rbr^{\frN} \ \subset 
\ \IK \lbr u,v \rbr^{\frN} \cap \um^{k+1}\ ,
$$
for all $k$; therefore $J_1 \cap \IK \lbr u,v \rbr^{\frN} = 0$;
i.e.,
$$
\IK \lbr u,v \rbr \ = \ J_1 \oplus \IK \lbr u,v \rbr^{\frN}\ .
$$
Moreover, if $f = g + h$, where $f \in \um^k$, $g \in J_1$ and
$h \in \IK \lbr u,v \rbr^{\frN}$, then $h \in 
\IK \lbr u,v \rbr^{\frN} \cap \um^k$, by (6.13); i.e.,
$h \in \um^k$, so $g \in \um^k$. This proves the claim.
\medskip

We have proved that the standard basis of $J_1$ gives the
standard basis of $J_0$ when we set $u = 0$. This gives the
variant of (b) where the presentation is merely {\it formal}.
We can, in fact, prove that the semicoherent presentation of the
Hilbert-Samuel function constructed in \cite[Theorem 9.6]{BMinv}
satisfies property (b) as stated :

We have shown that $\frN(\wcI_{X,a}) = \frN$ is a product
$\IN^{n-s} \times \frN^*$, where $\frN^* \subset \IN^s$ (corresponding
to the $v$-coordinates). Let $\ucH(a) = (N(a), \cH(a), \emptyset)$
denote the semicoherent presentation constructed in 
\cite[Theorem 9.6]{BMinv}; then the maximal contact submanifold
$N(a)$ has coordinates $w = (u,w_1)$ such that $(u,w_1) \mapsto u$
is the mapping $N(a) \to T$ (in particular, this mapping is smooth).
Because of this and the structure of \cite[Theorem 9.4]{BMinv},
we get a semicoherent presentation of $H_{X_{(t)},\cdot}$ at $a$
by setting $u = 0$. To see this, it is enough to observe that
the properties of \cite[Theorem 9.4]{BMinv} survive on setting certain
of the $w$-coordinates in the latter equal to $0$. (In particular,
the formal properties \cite[(7.2)(1)-(5)]{BMinv} survive on setting
certain of the formal $W$-coordinates in \cite[Theorem 9.4]{BMinv}
equal to $0$.) This completes the proof.
\end{proof}

\bibliographystyle{ams}

\end{document}